\documentclass[11pt]{article}
\usepackage{eurosym}
\usepackage{amssymb,amsmath,amsfonts,amsthm,enumerate}
\usepackage{setspace, natbib,graphicx,color,caption,subcaption,float,relsize}
\usepackage{url, hyperref}
\usepackage[nottoc,notlof,notlot]{tocbibind}
\usepackage{sectsty}
\usepackage{rotating}

\sectionfont{\large}
\subsectionfont{\large}
\subsubsectionfont{\large}
\makeatletter
\def\maketag@@@#1{\hbox{\m@th\normalfont\normalsize#1}}
\makeatother
\hypersetup{
  colorlinks = true, urlcolor = blue, linkcolor = blue, citecolor = magenta }

\DeclareMathOperator*{\argmin}{arg\,min}
\DeclareMathOperator*{\arginf}{arg\,inf}
\allowdisplaybreaks
\newtheorem{theorem}{Theorem}[section]
\newtheorem*{theorem32*}{Theorem 3.2}
\newtheorem*{theorem41*}{Theorem 4.1}
\newtheorem*{theorem62*}{Theorem 6.2}
\newtheorem{assumption}{Assumption}

\newtheorem{corollary}[theorem]{Corollary}

\newtheorem{example}{Example}

\newtheorem{lemma}[theorem]{Lemma}

\newtheorem{remark}{Remark}

\let\oldbibliography\thebibliography
\renewcommand{\thebibliography}[1]{  \oldbibliography{#1}  \setlength{\itemsep}{2pt}}
\begin{document}

\title{Weak-Identification Robust Wild Bootstrap applied to a Consistent
Model Specification Test}
\author{Jonathan B.~Hill\thanks{%
Dept. of Economics, University of North Carolina, Chapel Hill; www.unc.edu/$%
\sim $jbhill; jbhill@email.unc.edu.\medskip \newline
This paper was previously circulated under the title "\textit{Inference When
There is a Nuisance Parameter under the Alternative and Some Parameters are
Possibly Weakly Identified}". We thank two referees and Co-Editor Michael
Jansson for helpful comments and suggestions.} \\
%EndAName
University of North Carolina -- Chapel Hill}
\date{{\normalsize \today}}
\maketitle

\begin{abstract}
We present a new robust bootstrap method for a test when there is a nuisance
parameter under the alternative, and some parameters are possibly weakly or
non-identified. We focus on a \cite{Bierens1990}-type conditional moment
test of omitted nonlinearity for convenience. Existing methods include the
supremum p-value which promotes a conservative test that is generally not
consistent, and test statistic transforms like the supremum and average for
which bootstrap methods are not valid under weak identification. We propose
a new wild bootstrap method for p-value computation by targeting specific
identification cases. We then combine bootstrapped p-values across polar
identification cases to form an asymptotically valid p-value approximation
that is robust to any identification case. Our wild bootstrap procedure does
not require knowledge of the covariance structure of the bootstrapped
processes, whereas Andrews and Cheng's (%
\citeyear{AndrewsCheng2012,AndrewsCheng2013,AndrewsCheng2014}) simulation
approach generally does. Our method allows for robust bootstrap critical
value computation as well. Our bootstrap method (like conventional ones)
does not lead to a consistent p-value approximation for test statistic
functions like the supremum and average. We therefore smooth over the robust
bootstrapped p-value as the basis for several tests which achieve the
correct asymptotic level, and are consistent, for any degree of
identification. They also achieve uniform size control. A simulation study
reveals possibly large empirical size distortions in non-robust tests when
weak or non-identification arises. One of our smoothed p-value tests,
however, dominates all other tests by delivering accurate empirical size and
comparatively high power.\bigskip \newline
\end{abstract}

\setstretch{1}

\setstretch{1.3}

\section{INTRODUCTION\label{sec:intro}}

We present a new bootstrap procedure for non-standard tests where some
regression model parameters may be weakly or non-identified. We focus ideas
at the expense of greater generality by working with a regression model that
has additive nonlinearity:%
\begin{equation}
y_{t}=\zeta _{0}^{\prime }x_{t}+\beta _{0}^{\prime }g(x_{t},\pi
_{0})+\epsilon _{t}=f(\theta _{0},x_{t})+\epsilon _{t}\text{ where }x_{t}\in
\mathbb{R}^{k_{x}}\text{ and }\theta \equiv \left[ \zeta ^{\prime },\beta
^{\prime },\pi ^{\prime }\right] ^{\prime }\text{.}  \label{model_add}
\end{equation}%
The variable $y_{t}$ is a scalar, $x_{t}$ $\in $ $\mathbb{R}^{k_{x}}$ are
covariates with a constant term and finite $k_{x}$ $\geq $ $2$, $g$ $:$ $%
\mathbb{R}^{k_{x}}$ $\times $ $\Pi $ $\rightarrow $ $\mathbb{R}^{k_{\beta }}$
is a known function, and $\zeta _{0}$ $\in \mathcal{Z}$, $\beta _{0}$ $\in $
$\mathcal{B}$ and $\pi _{0}$ $\in $ $\Pi $, where $\mathcal{B}$, $\mathcal{Z}
$ and $\Pi $\ are compact subsets of $\mathbb{R}^{k_{\beta }}$, $\mathbb{R}%
^{k_{x}}$ and $\mathbb{R}^{k_{\pi }}$\ respectively for finite $(k_{\beta
},k_{\pi })$ $\geq $ $1.$ $x_{t}$ includes a constant term and at least one
stochastic regressor, and let $E[\epsilon _{t}]$ $=$ $0$ and $E[\epsilon
_{t}^{2}]$ $\in $ $\left( 0,\infty \right) $ for some unique $\theta _{0}$ $%
\in $ $\Theta $ $\equiv $ $\mathcal{Z}$ $\times $ $\mathcal{B}$ $\times $ $%
\Pi .$

We want to test $H_{0}$ $:$ $E[y_{t}|x_{t}]$ $=$ $f(\theta _{0},x_{t})$ $%
a.s. $ for some unique $\theta _{0}$ against a general alternative $%
H_{1}:\sup_{\theta \in \Theta }P(E[y_{t}|x_{t}]$ $=$ $f(\theta ,x_{t}))$ $<$
$1$. We assume that the (pseudo) true value $\theta _{0}$ minimizes a
standard criterion function
\citep[see,
e.g.,][]{KullbackLeibler1951,Sawa1978,White1982}. In order to test $H_{0}$,
we work with the \cite{Bierens1990} type conditional moment [CM] test of
omitted nonlinearity for convenience.

Under $H_{1}$ it is known that $E[\epsilon _{t}F(\lambda ^{\prime }x_{t})]$ $%
\neq $ $0$ for a large class of weight functions $F$ $:$ $\mathbb{R}$ $%
\rightarrow $ $\mathbb{R}$, and $\forall \lambda $ $\in $ $\Lambda
/S_{\Lambda }$ where $\Lambda $ is any compact subset of $\mathbb{R}^{k_{x}}$
and $S_{\Lambda }$ $\subset $ $\Lambda $ has measure zero ($S_{\Lambda }$
depends on $F$ and the underlying distribution). Examples of $F$ include the
exponential \citep{Bierens1982,Bierens1990,deJong1996}, logistic %
\citep{White1989}, and the covering class of non-polynomial real analytic
functions \citep{StinchWhite1998}. A CM test operates on a normalized sample
version of $E[\epsilon _{t}F(\lambda ^{\prime }x_{t})]$, cf. \cite{Newey1985}
and \cite{Tauchen1985}, and therefore has the nuisance parameter $\lambda $
under $H_{1}$. As an example, consider testing whether $y_{t}$ is governed
by a Logistic Smooth Transition AR($p$) process with a single transition
function. Under the alternative the true process may be LSTAR with multiple
transition functions, or STAR with a different transition function (e.g.
exponential, normal), or may not be in the STAR class at all (e.g. Self
Exciting Threshold Autoregression, cf. \cite{TongLim1980}). STAR\ model
estimation generally involves the possibility of weakly or non-identified
parameters that is routinely assumed away in the STAR literature %
\citep[cf.][]{Terasvirta1994}. See Example \ref{Example STAR} in Section \ref%
{sec:assum} for further details.

If $\beta _{0}$ $=$ $0$ then $\pi _{0}$ is not identified. In fact, if $n$
is the sample size and there is local drift $\beta _{0}$ $=$ $\beta _{n}$ $%
\rightarrow $ $0$ with $\sqrt{n}||\beta _{n}||$ $\rightarrow $ $[0,\infty )$%
, then estimators for $\theta _{0}$ have nonstandard limit distributions,
and estimators of $\pi _{0}$ have a random probability limit. See %
\citet{AndrewsCheng2012,AndrewsCheng2013,AndrewsCheng2014} for a broad
literature review, and for results on estimation and classic inference
generally under the assumption of model correctness $E[\epsilon _{t}|x_{t}]$
$=$ $0$ $a.s.$ See also \cite{Cheng2015}. We assume throughout that $\pi
_{0} $ is identified when $\beta _{0}$ $\neq $ $0$. Otherwise an approach
similar to \cite{Cheng2015} would be appropriate, leading to more intense
asymptotics. In the weak identification literature in which a regression
model forms the basis of study, correct model specification in the sense
that $\epsilon _{t}$ is iid, a martingale difference or $E[\epsilon
_{t}|x_{t}]$ $=$ $0$ $a.s.$ is typically assumed. Thus, $H_{0}$ above is
assumed to be true. Consider, e.g., \cite{Sargan1983}, \cite{Phillips1990},
\cite{ChoiPhillips1992}, \cite{StockWright2000}, \cite%
{AndrewsMoreiraStock2006}, \citet[Example 1, Section
6]{AndrewsCheng2012}, \citet[Examples 1 and 2, Section 7]{AndrewsCheng2013}, %
\citet[Example 2]{AndrewsCheng2014}, \cite{Cheng2015}, and
\citet[Example
(1)]{McCloskey2017} amongst others. A broad literature exists on weak
identification related to weak instruments %
\citep[e.g.][]{Dufour1997,StockWright2000,Moreira2003,AndrewsMoreiraStock2006}%
. This is not our primary focus since the source of weak identification is a
specific feature of the regression model.\footnote{%
Many treatments in the weak identification literature do not focus on a
regression model, but work with unconditional moment conditions, or on a
(non-stochastic) parametric function. See \cite{Dufour1997} and \cite%
{Caner2010}. See also \cite{ElliotMullerWatson2015}.}

Conversely, in the omitted nonlinearity test literature strong
identification is universally assumed or implied. This translates here to
assuming $\beta _{0}$ $\neq $ $0$, or simply testing whether $y_{t}$ $=$ $%
\zeta ^{\prime }x_{t}$ $+$ $\epsilon _{t}$ is the correct specification.
This literature is equally massive: see, e.g., %
\citet{Bierens1982,Bierens1990}, \cite{White1989}, \cite{HongWhite1995},
\cite{Hansen1996}, \cite{BierensPloberger1997}, \cite{deJong1997}, \cite%
{StinchWhite1998}, \cite{Dette1999}, \cite{Li_1999}, \cite{Whang2000}, \cite%
{DelgadoDominguezLavergne2006}, \cite{Hill2008}, \cite{DavidsonHalyunga2014}
and \cite{LiLiLiu2016}.

Our contributions are twofold. First, we deliver a first-time bridge between
these literatures: a consistent model specification test of $H_{0}$ that is
robust to the full sweep of identification cases %
\citep[cf.][]{AndrewsCheng2012}. Thus, since we test for correct model
specification $E[y_{t}|x_{t}]$ $=$ $f(\theta _{0},x_{t})$ $a.s.$ we do not
assume it a priori, contrary to many offerings in the weak identification
literature. We must, however, make some assumptions on the model error $%
\epsilon _{t}$ in order to identify the (possibly pseudo-true) model
parameters. In a similar vein \citet[Example 1]{AndrewsCheng2014} only
impose a weak orthogonality condition because a regressor may be endogenous,
but only treat iid data in a linear model. We allow for a nonlinear time
series setting and only require a weak orthogonality condition (under the
alternative).

Second, we provide a bootstrap procedure that is robust to the degree of
identification. This topic has apparently been ignored to date. Our method
broadly applies to other tests, including t-, Wald, Lagrange Multiplier, and
QLR tests, as well as other model specification tests including
nonparametric tests, although we restrict attention to a CM statistic for
brevity.

The presence of $\lambda $ prompts a test statistic transform detailed
below. This promotes a nonstandard asymptotic theory and therefore requires
a bootstrap method \citep[e.g.]{Hansen1996}. But the possibility of weak
identification alone leads to nonstandard asymptotics. Thus, even if a
nonparametric model specification test is explored which bypasses a nuisance
parameter, including \cite{HardleMammen1993}, \cite{HongWhite1995},and \cite%
{Zheng1996}, a nonstandard approach is required to handle allowing for any
degree of identification. This paper proposes a new bootstrap method
suitable for the nuisance parametric approach, that is robust to weak
identification. The method is general, and can therefore be extended in
principle to any other model specification test approach.

Let $\mathcal{T}_{n}(\lambda )$ $\geq $ $0$ be the proposed CM test
statistic. In the setting of (\ref{model_add}) and mild regularity
conditions, $\mathcal{T}_{n}(\lambda )$\ has a chi-squared limit law when $%
\sqrt{n}||\beta _{n}||$ $\rightarrow $ $\infty $, and otherwise has a
non-standard limit. This represents polar cases of semi-strong or strong
identification, and weak or non-identification %
\citep[cf.][]{AndrewsCheng2012}. %
\citet{AndrewsCheng2012,AndrewsCheng2013,AndrewsCheng2014} propose a robust
critical value for t, Wald and QLR statistics, where simulated data is used
to approximate the limit distribution.

In the following, unless confusion cannot be avoided, we say "weak
identification" to mean non- or weak cases, and "strong identification" to
mean semi-strong or strong cases.

The simulation approach of %
\citet{AndrewsCheng2012,AndrewsCheng2013,AndrewsCheng2014} requires
knowledge of the covariance kernel of the simulated stochastic process. This
may be intractable when the weakly identified parameter is non-scalar.
Further, simulating an asymptotic distribution presumes the latter well
approximates the small sample distribution. This may fail to be true when
there exists conditional heteroskedasticity, when error tails are
leptokurtic, when the parameter dimension is large, and/or when the sample
size is small. The typical solution is a bootstrap or sub-sampling method.
In the case of testing $H_{0}$, a natural method is the wild (or multiplier)
bootstrap as in \cite{Hansen1996}, cf. \cite{Wu1986} and \cite{Liu1988}. A
major advantage of the wild bootstrap premise over the simulation method is
that knowledge of the covariance kernel of the bootstrapped process is not
required. Bootstrap methods applied to $\mathcal{T}_{n}(\lambda )$ that do
not take into consideration the possibility of weak identification, however,
are asymptotically invalid because the weak limit of $\mathcal{T}%
_{n}(\lambda )$ exhibits discontinuities with respect to $\beta _{n}$, in
which case uniform asymptotics fail. See \cite{BickelFreedman1981}, \cite%
{Romano1989}, \cite{SheehyWellner1992} and \cite{AndrewsGuggenberger2010}.
See \cite{GineZinn1990} for discussion on types of uniformity in bootstrap
environments.

We solve the problem of non-uniformity by targeting the wild bootstrap to
identification category specific first order expansions of $\mathcal{T}%
_{n}(\lambda )$ under the null. Once bootstrapped p-values are computed for
polar identification cases, we combine them as in \cite{AndrewsCheng2012}
using their notions of \textit{Least Favorable} and \textit{Identification
Category Selection} constructions (they develop critical values). The result
is an asymptotically valid p-value approximation $\hat{p}_{n}(\lambda )$,
irrespective of the degree of identification.

The Bonferroni-based size correction approach of \cite{McCloskey2017} is a
plausible alternative to the LF and ICS methods of \cite{AndrewsCheng2012}.
The theory there is presented for a test of a fixed parameter value, where
other (model-based nuisance) parameters are also present and cause the test
statistic limit distribution discontinuity. In McCloskey's (%
\citeyear{McCloskey2017}: eq. (1)) example, the nuisance parameter is part
of the data generating process. In our setting we test whether a chosen
model is correct, where a nuisance parameter arises that is not part of the
data generating process, and is due solely to the construction of a test
statistic. Discontinuity of the limit distribution is not caused by the
nuisance parameter, but by a parameter subset from the model. We leave for
another venue a consideration of generalizing the Bonferroni-based size
correction to our setting.

In order to handle the nuisance parameter, we randomize $\lambda $, use the
classic sup-transform $\sup_{\lambda \in \Lambda }\hat{p}_{n}(\lambda )$
\citep[see, e.g.,][Chapter
3.1]{Lehmann1994}, and use the P-Value Occupation Time [PVOT] $\mathcal{\hat{%
P}}_{n}(\alpha )$ $\equiv $ $\int_{\Lambda }I(\hat{p}_{n}(\lambda )$ $<$ $%
\alpha )d\lambda $ where $I(A)$ $=$ $1$ if $A$ is true, and $\alpha $ $\in $
$(0,1)$ is the nominal level \citep{Hill2018}. Randomizing $\lambda $
sacrifices power \citep[e.g.][]{White1989}. $\sup_{\lambda \in \Lambda }\hat{%
p}_{n}(\lambda )$ by construction promotes a conservative test that is
generally not consistent since a \cite{Bierens1990}-type CM test is not
known to be consistent for all $\lambda $ \citep{Bierens1990,StinchWhite1998}%
. We also present conditions under which our tests achieve uniform size
control for any degree of identification.

The challenge of constructing valid tests in the presence of nuisance
parameters under $H_{1}$ dates at least to \cite{ChernoffZacks1964} and %
\citet{Davies77,Davies87}. Nuisance parameters that are not identified under
$H_{1}$ are either chosen at random \citep[e.g.][]{White1989}; or $\mathcal{T%
}_{n}(\lambda )$ is smoothed over $\Lambda $, resulting in a non-standard
limit distribution %
\citep[e.g.][]{ChernoffZacks1964,Davies77,AndrewsPloberger1994}; or a
computed p-value like $\hat{p}_{n}(\lambda )$ is smoothed. Examples of
transforms are the average $\int_{\Lambda }\mathcal{T}_{n}(\lambda )\mu
(d\lambda )$ and supremum $\sup_{\lambda \in \Lambda }\mathcal{T}%
_{n}(\lambda )$, where $\mu (\lambda )$ is a measure on $\Lambda $ that is
absolutely continuous with respect to Lebesgue measure %
\citep{ChernoffZacks1964,Davies77,AndrewsPloberger1994}. \cite%
{BierensPloberger1997} integrate the squared numerator from a conventional
CM statistic $\mathcal{T}_{n}(\lambda )$, resulting in the Integrated
Conditional Moment [ICM] test, cf. \cite{Bierens1982}.

Unless strong identification is assumed, then $\int_{\Lambda }\mathcal{T}%
_{n}(\lambda )\mu (d\lambda )$, $\sup_{\lambda \in \Lambda }\mathcal{T}%
_{n}(\lambda )$ and the ICM\ cannot be consistently bootstrapped by
conventional methods or our method. The intuition is simple. We can write $%
\mathcal{T}_{n}(\lambda )$ $=$ $\mathcal{Z}_{n}^{2}(\lambda )$ for some
sample process $\{\mathcal{Z}_{n}^{2}(\lambda )\}$, e.g.
\citet[eq.
(18)]{Bierens1990}. Under strong identification and fairly general
assumptions $\{\mathcal{Z}_{n}(\lambda )$ $:$ $\lambda $ $\in $ $\Lambda \}$
converges to a Gaussian process $\{\mathcal{Z}(\lambda )$ $:$ $\lambda $ $%
\in $ $\Lambda \}$ with covariance kernel $E[\mathcal{Z}(\lambda )\mathcal{Z}%
(\tilde{\lambda})]$ that generally depends on $\theta _{0}$ and $E[\epsilon
_{t}^{2}]$. See Theorem \ref{th:CM_weak} below. In order to use $\mathcal{Z}%
_{n}(\lambda )$ to obtain bootstrap draws from the process $\{\mathcal{Z}%
(\lambda )$ $:$ $\lambda $ $\in $ $\Lambda \}$ we therefore need consistent
estimators for $\theta _{0}$ and $E[\epsilon _{t}^{2}]$. That is impossible
if $\pi _{0}$ is truly only weakly identified \citep[cf.][]{AndrewsCheng2012}%
. The same problem applies to Bierens and Ploberger's (%
\citeyear{BierensPloberger1997}) ICM test.

Our setting is decidedly different from Hansen's (\citeyear{Hansen1996}) who
tests $\beta _{0}$ $=$ $0$ and treats $\pi _{0}$ as an unidentified nuisance
parameter under the null. We do not require $F(\lambda ^{\prime }x_{t})$ to
be part of the true data generating process under $H_{0}$, and we estimate
all parameters $\theta _{0}$ allowing for weak identification. Moreover,
\cite{Hansen1996} delivers a valid bootstrap method for test statistic
transforms like $\int_{\Lambda }\mathcal{T}_{n}(\lambda )\mu (d\lambda )$
and $\sup_{\lambda \in \Lambda }\mathcal{T}_{n}(\lambda )$ under strong
identification. If any identification category is allowed. then neither his
nor our bootstrap methods are valid for such \textit{test statistic}
transforms, and we are unaware of any bootstrap method that is valid. The
PVOT\textit{\ p-value} transform, however, does lend itself to weak
identification robust inference \citep[cf.][]{Hill2018}.

We work with p-values due to their convenience of interpretation: one
p-value can be used to test $H_{0}$ at any desired level of significance,
although our bootstrap method can also be used for robust critical value
approximations $\hat{c}_{1-\alpha ,n}(\lambda )$
\citep[see][Appendix
E]{Supp_Mat_2020}. An unavoidable difference in theory, however, is $\hat{c}%
_{1-\alpha ,n}(\lambda )$ leads to an asymptotically correctly \textit{sized}
test, while $\hat{p}_{n}(\lambda )$ only promotes a test with correct
asymptotic \textit{level} \footnote{%
Let $\alpha $ be the desired significance level, and let $AsySz$ be the
asymptotic size of a test. The asymptotic \textit{level} of the test is $%
\alpha $ if $AsySz$ $\leq $ $\alpha $.}. Ultimately this is due to weak
identification and the way parameters enter $\hat{p}_{n}(\lambda )$: see
Theorem \ref{th:pv_weak} and its proof. In simulation experiments not
reported here, however, robust critical and p-values perform essentially
identically.

Our tests are consistent irrespective of the choice of $\Lambda $, although
for a given $\Lambda $ power in small samples is naturally amplified in
certain directions away from the null. These issues are well known and not
dealt with in this paper.

A simulation experiment reveals tests based on $\mathcal{T}_{n}(\lambda
^{\ast })$ with randomly selected $\lambda ^{\ast }$, $\sup_{\lambda \in
\Lambda }\mathcal{T}_{n}(\lambda )$, \linebreak $\int_{\Lambda }\mathcal{T}%
_{n}(\lambda )\mu (d\lambda )$ and the PVOT with a conventional wild
bootstrapped p-value $p_{n}(\lambda )$ are all strongly over-sized under
weak-identification. Somewhat ironically, the conservative test based on $%
\sup_{\lambda \in \Lambda }p_{n}(\lambda )$ counters the large size
distortion under weak identification, but results in low power.\footnote{$%
\sup_{\lambda \in \Lambda }p_{n}(\lambda )$ is not robust to identification
category: its conservativeness merely tempers the degree of size distortion
under weak identification.} The test based on our robust $p$-value $\hat{p}%
_{n}(\lambda ^{\ast })$, however, achieves the correct level, but has
comparatively low power, while $\sup_{\lambda \in \Lambda }\hat{p}%
_{n}(\lambda )$ is conservative with low power. The PVOT test with $\hat{p}%
_{n}(\lambda )$ in simulation experiments has the correct size, and under
weak identification achieves the highest power.

Our approach is parametric, while in the nonparametric literature weak
identification robust methods are increasingly popular. Nevertheless, robust
bootstrap procedures have not apparently been treated. See, for example,
\cite{AndrewsMikusheva2016}, \cite{Cox2016}, \cite{HanMcCloskey2016}, and
\cite{McCloskey2017} and the references provided there.

The remainder of the paper is organized as follows. Section \ref{sec:cm_test}
presents the CM statistic and its transforms. Assumptions and main results
are presented in Sections \ref{sec:assum} and \ref{sec:CM_main}. In Sections %
\ref{sec:CM_pv_construct} and \ref{sec:CM_pv} we present robust p-values and
develop a method for bootstrapping the robust p-values. A simulation study
is contained in Section \ref{sec:sim} and concluding remarks follow in
Section \ref{sec:conclusion}. Proofs are given in Appendix \ref{app:A}%
.\medskip

We use the following notation. $[z]$ rounds $z$ to the nearest integer. $%
I(\cdot )$ is the indicator function: $I(A)$ $=$ $1$ if $A$ is true,
otherwise $I(A)$ $=$ $0$. $a_{n}/b_{n}$ $\sim $ $c$ implies $a_{n}/b_{n}$ $%
\rightarrow $ $c$ as $n$ $\rightarrow $ $\infty $. $|\cdot |$ is the $l_{1}$%
-matrix norm; $||\cdot ||$ is the Euclidean norm; $||\cdot ||_{p}$ is the $%
L_{p}$-norm. $K$ $>$ $0$ is a finite constant whose value may change from
place to place. $0_{a\times b}$ is an $a\times b$ dimensional matrix of
zeros. \emph{a.e.} denotes \emph{almost everywhere}. $\Rightarrow ^{\ast }$
denotes weak convergence on $l_{\infty }$, the space of bounded functions
with sup-norm topology, in the sense of \citet{HoffJorg1984,HoffJorg1991},
cf. \cite{Dudley1978} and \citet{Pollard1984,Pollard1990}.

\section{TEST STATISTIC CONSTRUCTION\label{sec:cm_test}}

Let $\{\beta _{n}\}$ be the drifting sequence such that $\lim_{n\rightarrow
\infty }\beta _{n}$ $=$ $\beta _{0}$. As in %
\citet{AndrewsCheng2012,AndrewsCheng2013}, technical results are derived
under two overlapping cases which align with the following three categories:
I.a. $\beta _{n}$ $=$ $\beta _{0}$ $=$ $0$ $\forall n$ $\geq $ $1$ ($\pi
_{0} $ is \textit{unidentified}); I.b. $n^{1/2}\beta _{n}$ $\rightarrow $ $b$
$\in $ $\mathbb{R}/0$ hence $\beta _{0}$ $=$ $0$ ($\pi _{0}$ is \textit{%
weakly identified}); II. $n^{1/2}||\beta _{n}||$ $\rightarrow $ $\infty $
hence $\beta _{0}\in $ $\mathbb{R}$ ($\pi _{0}$ is \textit{semi-strongly
identified}); and III. $\beta _{n}$ $\rightarrow $ $\beta _{0}$ $\neq $ $0$ (%
$\pi _{0}$ is \textit{strongly identified}).

The two key over-lapping cases for all asymptotic results are denoted as %
\citep[see][eq. (2.7)]{AndrewsCheng2012}:%
\begin{eqnarray*}
&&\mathcal{C}(i,b)\text{. }\beta _{n}\rightarrow \beta _{0}=0\text{ and }%
\sqrt{n}\beta _{n}\rightarrow b\text{ where }b\in (\mathbb{R}\cup \{\pm
\infty \})^{k_{\beta }} \\
&&\mathcal{C}(ii,\omega _{0})\text{. }\beta _{n}\rightarrow \beta _{0}\text{
where }\beta _{0}\gtreqless 0\text{, }\sqrt{n}\left\Vert \beta
_{n}\right\Vert \rightarrow \infty ,\text{ and }\beta _{n}/\left\Vert \beta
_{n}\right\Vert \rightarrow \omega _{0}\text{ where }\left\Vert \omega
_{0}\right\Vert =1.
\end{eqnarray*}%
Case $\mathcal{C}(i,b)$ contains sequences $\beta _{n}$ close to zero, and
when $||b||$ $<$ $\infty $ then $\pi _{0}$ is either weakly or
non-identified. Case $\mathcal{C}(ii,\omega _{0})$ contains sequences $\beta
_{n}$ farther from zero, covering semi-strong ($\beta _{0}$ $=$ $0$ and $%
\sqrt{n}||\beta _{n}||$ $\rightarrow $ $\infty $) and strong ($\beta _{0}$ $%
\neq $ $0$) identification for $\pi _{0}$. Notice $b$ and $\omega _{0}$
represent two different limits: $\sqrt{n}\beta _{n}$ $\rightarrow $ $b$\
versus $\beta _{n}/||\beta _{n}||$ $\rightarrow $ $\omega _{0}$.

Let $F$ $:$ $\mathbb{R}$ $\rightarrow $ $\mathbb{R}$ be a real analytic and
non-polynomial function, and $\mathcal{W}$ $:$ $\mathbb{R}^{k_{x}}$ $%
\rightarrow $ $\mathbb{R}^{k_{x}}$ is a one-to-one and bounded function.
Under $H_{0}$, $E[\epsilon _{t}|x_{t}]$ $=$ $0$ $a.s.,$ hence $E[\epsilon
_{t}F(\lambda ^{\prime }\mathcal{W}(x_{t}))]$ $=$ $0$. Under $H_{1}$, $%
E[\epsilon _{t}F(\lambda ^{\prime }\mathcal{W}(x_{t}))]$ $\neq $ $0$ $%
\forall \lambda $ $\in $ $\Lambda /\mathcal{S}_{\Lambda }$ where $\mathcal{S}%
_{\Lambda }$ has Lebesgue measure zero. See Lemma 1 in \cite{Bierens1990}
for iid data and exponential $F(\cdot )$, see \cite{StinchWhite1998} for
broad theory treating the analytic class, and see \cite{deJong1996} and \cite%
{Hill2008} for the time series case.\footnote{\cite{deJong1996} also allows
for an infinite dimensional conditioning set, e.g. $y_{t-1},y_{t-2},...$ in
an ARMA model. This leads to more nuanced results for identifying whether a
regression model is mis-specified.}

We first require an estimation setting on a chosen estimation parameter
space $\Theta $. Let $y_{t}$ exist on the probability measure space $(\Omega
,\mathcal{P},\mathcal{F})$, where $\mathcal{F}$ $\equiv $ $\sigma (\cup
_{t\in \mathbb{Z}}\mathcal{F}_{t})$ and $\mathcal{F}_{t}$ $\equiv $ $\sigma
(y_{\tau }$ $:$ $\tau $ $\leq $ $t)$. Assume $\Theta $ has the form $%
\{\theta $ $\equiv $ $[\beta ^{\prime },\zeta ^{\prime },\pi ^{\prime
}]^{\prime }$ $:$ $\beta $ $\in $ $\mathcal{B},\zeta $ $\in $ $\mathcal{Z}%
(\beta ),\pi $ $\in $ $\Pi \},$ where $\mathcal{B}$, $\mathcal{Z}(\beta )$
for each $\beta ,$ and $\Pi $ are compact subsets. $\mathcal{Z}(\beta )$
depends on $\beta $ because parameter restrictions may be imposed to ensure
a stationary solution. Define the parameter subset and space%
\begin{equation*}
\psi \equiv \left[ \beta ^{\prime },\zeta ^{\prime }\right] ^{\prime }\in
\Psi \equiv \{(\beta ,\zeta ):\beta \in \mathcal{B},\zeta \in \mathcal{Z}%
(\beta )\}.
\end{equation*}%
The true parameter space $\Theta ^{\ast }$ $=$ $\Psi ^{\ast }$ $\times $ $%
\Pi ^{\ast }$ $=$ $\{\theta $ $\equiv $ $[\beta ^{\prime },\zeta ^{\prime
},\pi ^{\prime }]^{\prime }$ $:$ $\beta $ $\in $ $\mathcal{B}^{\ast },\zeta $
$\in $ $\mathcal{Z}^{\ast }(\mathcal{\beta }),\pi $ $\in $ $\Pi ^{\ast }\}$
lies in the interior of $\Theta $, it contains $\theta _{0}$ $\equiv $ $%
[\beta _{0}^{\prime },\zeta _{0}^{\prime },\pi _{0}^{\prime }]^{\prime }$ ,
and $0$ $\in $ $\mathcal{B}^{\ast }$. The dependence of $\mathcal{Z}^{\ast }(%
\mathcal{\beta })$ on $\mathcal{\beta }$ ensures $\Theta ^{\ast }$\ contains
points consistent with stationarity and moment and memory properties imposed
under Assumption \ref{assum:dgp} below. The spaces $\Theta ^{\ast }$ $%
\subset $ $\Theta $ are assumed different to ensure the true value $\theta
_{0}$ does not lie on the boundary of $\Theta $ for convenience of focus.

The sample is $\{(y_{t},x_{t})\}_{t=1}^{n}$. We work with least squares to
reduce notation, but an extension to a broad class of extremum estimators is
straightforward. Define $\epsilon _{t}(\theta )$ $\equiv $ $y_{t}$ $-$ $%
\zeta ^{\prime }x_{t}$ $-$ $\beta ^{\prime }g(x_{t},\pi )$, and define the
least squares criterion and estimator:%
\begin{equation*}
Q_{n}(\theta )=Q_{n}(\psi ,\pi )\equiv \frac{1}{2}\frac{1}{n}%
\sum_{t=1}^{n}\epsilon _{t}^{2}(\theta )\text{ \ and \ }\hat{\theta}%
_{n}\equiv \arginf_{\theta \in \Theta }Q_{n}(\theta ).
\end{equation*}

The criterion ensures we can express $\hat{\theta}_{n}$ as a concentrated
estimator $\hat{\theta}_{n}$ $=$ $[\hat{\psi}_{n}^{\prime }(\hat{\pi}_{n}),%
\hat{\pi}_{n}^{\prime }]^{\prime },$ where%
\begin{equation*}
\hat{\psi}_{n}(\pi )=\arginf_{\psi \in \Psi }Q_{n}(\psi ,\pi )\text{ and }%
\hat{\pi}_{n}=\arginf_{\pi \in \Pi }Q_{n}(\hat{\psi}_{n}(\pi ),\pi ).
\end{equation*}%
Under weak identification, a suitably normalized $\hat{\theta}_{n}$ has a
non-standard limit distribution, and must be partitioned into $[\sqrt{n}(%
\hat{\psi}_{n}(\hat{\pi}_{n})$ $-$ $\psi _{n})^{\prime },\hat{\pi}%
_{n}^{\prime }]^{\prime }$ since $\hat{\pi}_{n}$ has a stochastic
probability limit, cf. \cite{AndrewsCheng2012}. Thus, we cannot work with
Bierens' (\citeyear{Bierens1990}) original test statistic, nor the
environments of \cite{White1989}, \cite{BierensPloberger1997}, and many
others cited in Section \ref{sec:intro}, since these rely on a first order
expansion of $\sqrt{n}(\hat{\theta}_{n}$ $-$ $\theta _{0})$ in order to
characterize a suitable normalizing scale, implicitly ignoring weak
identification \citep[see, e.g.,][p.
1446]{Bierens1990}.

The robust test statistic is constructed as follows. Define%
\begin{eqnarray}
&&d_{\psi ,t}(\pi )\equiv \left[ g(x_{t},\pi )^{\prime },x_{t}^{\prime }%
\right] ^{\prime }\text{ \ and \ }d_{\theta ,t}(\omega ,\pi )\equiv \left[
g(x_{t},\pi )^{\prime },x_{t}^{\prime },\omega ^{\prime }\frac{\partial }{%
\partial \pi }g(x_{t},\pi )\right] ^{\prime }  \notag \\
&&\widehat{\mathcal{H}}_{n}=\frac{1}{n}\sum_{t=1}^{n}d_{\theta ,t}(\omega (%
\hat{\beta}_{n}),\hat{\pi}_{n})d_{\theta ,t}(\omega (\hat{\beta}_{n}),\hat{%
\pi}_{n})^{\prime }\text{ where }\omega (\beta )\equiv \left\{
\begin{array}{ll}
\beta /\left\Vert \beta \right\Vert & \text{if }\beta \neq 0 \\
1_{k_{\beta }}/\left\Vert 1_{k_{\beta }}\right\Vert & \text{if }\beta =0%
\end{array}%
\right.  \label{Hn} \\
&&\mathfrak{\hat{b}}_{\theta ,n}(\omega ,\pi ,\lambda )\equiv \frac{1}{n}%
\sum_{t=1}^{n}F\left( \lambda ^{\prime }\mathcal{W}(x_{t})\right) d_{\theta
,t}(\omega ,\pi )  \notag \\
&&\hat{v}_{n}^{2}(\hat{\theta}_{n},\lambda )\equiv \frac{1}{n}%
\sum_{t=1}^{n}\epsilon _{t}^{2}(\hat{\theta}_{n})\left\{ F\left( \lambda
^{\prime }\mathcal{W}(x_{t})\right) -\mathfrak{\hat{b}}_{\theta ,n}(\omega (%
\hat{\beta}_{n}),\hat{\pi}_{n},\lambda )^{\prime }\widehat{\mathcal{H}}%
_{n}^{-1}d_{\theta ,t}(\omega (\hat{\beta}_{n}),\hat{\pi}_{n})\right\} ^{2}.
\notag
\end{eqnarray}%
See \citet[p. 2175]{AndrewsCheng2012} and \citet[p. 40]{AndrewsCheng2013}\
for discussions on rescaling with $||\beta ||$ to avoid a singular Hessian
matrix under semi-strong identification, i.e. $\beta _{0}$ $=$ $0$ and $%
\sqrt{n}||\beta _{n}||$ $\rightarrow $ $\infty $. The CM statistic we use
is:
\begin{equation*}
\mathcal{T}_{n}(\lambda )\equiv \left( \frac{1}{\hat{v}_{n}(\hat{\theta}%
_{n},\lambda )}\frac{1}{\sqrt{n}}\sum_{t=1}^{n}\epsilon _{t}(\hat{\theta}%
_{n})F\left( \lambda ^{\prime }\mathcal{W}(x_{t})\right) \right) ^{2}.
\end{equation*}

As discussed in Section \ref{sec:intro}, transforms like $\sup_{\lambda \in
\Lambda }\mathcal{T}_{n}(\lambda )$, $\int_{\Lambda }\mathcal{T}_{n}(\lambda
)\mu (d\lambda )$, and the ICM statistic $\int_{\Lambda }\{1/\sqrt{n}%
\sum_{t=1}^{n}\epsilon _{t}(\hat{\theta}_{n})\times $ $F(\lambda ^{\prime }%
\mathcal{W}(x_{t}))\}^{2}\mu (d\lambda )$ cannot be consistently
bootstrapped by our or apparently any other bootstrap method. We therefore
focus on p-value smoothing since we \textit{can} consistently bootstrap a
p-value approximation for $\mathcal{T}_{n}(\lambda )$ for a given $\lambda $
(see Sections \ref{sec:CM_pv_construct} and \ref{sec:CM_pv}). Let $\hat{p}%
_{n}(\lambda )$\ be a bootstrapped p-value. Along with $\hat{p}_{n}(\lambda
^{\ast })$ with randomly selected $\lambda ^{\ast }$, and $\sup_{\lambda \in
\Lambda }\hat{p}_{n}(\lambda )$, we use Hill's (\citeyear{Hill2018}) P-Value
Occupation Time [PVOT]:
\begin{equation*}
\mathcal{\hat{P}}_{n}(\alpha )\equiv \int_{\Lambda }I\left( \hat{p}%
_{n}(\lambda )<\alpha \right) d\lambda \text{ where }\int_{\Lambda }d\lambda
=1\text{ is assumed.}
\end{equation*}%
If $\int_{\Lambda }d\lambda $ $\neq $ $1$ then we use $\int_{\Lambda }I(\hat{%
p}_{n}(\lambda )$ $<$ $\alpha )d\lambda /\int_{\Lambda }d\lambda $. \cite%
{Hill2018} shows under general conditions that are verified here that $%
\lim_{n\rightarrow \infty }P(\mathcal{\hat{P}}_{n}(\alpha )$ $<$ $\alpha )$ $%
\leq $ $\alpha $ such that the PVOT test has correct asymptotic level.%
\footnote{%
Simulation experiments here and in \cite{Hill2018}\ reveal sharp size for a
PVOT tests of functional form, GARCH effects and a one time structural
break, suggesting $\lim_{n\rightarrow \infty }P(\mathcal{\hat{P}}_{n}(\alpha
)$ $<$ $\alpha )$ $=$ $\alpha $ likely holds in these and similar cases.} We
prove below that the PVOT statistic with the identification robust p-value
also achieves uniform size control.

\section{ASSUMPTIONS\label{sec:assum}}

Recall $\psi _{0}$ $=$ $[\beta _{0}^{\prime },\zeta _{0}^{\prime }]^{\prime
} $ $\in $ $\mathbb{R}^{k_{\psi }}$ where $k_{\psi }$ $=$ $k_{\beta }$ $+$ $%
k_{x}$, and $\pi _{0}$ $\in $ $\mathbb{R}^{k_{\pi }}$. The following matrix
is used to standardize the criterion gradient process below, ensuring a
non-degenerate limit under weak identification
\citep[e.g.][Sect.
3]{AndrewsCheng2012}:%
\begin{equation}
\mathfrak{B}(\beta )=\left[
\begin{array}{ll}
I_{k_{\psi }} & 0_{k_{\psi }\times k_{\pi }} \\
0_{k_{\pi }\times k_{\psi }} & \left\Vert \beta \right\Vert \times I_{k_{\pi
}}%
\end{array}%
\right] .  \label{B(b)}
\end{equation}%
The definition of $\mathfrak{B}(\beta )$ assumes $\beta $\ is a vector; in
the scalar case replace $||\beta ||$ with $\beta $
\citep[see the discussion
following eq. (3.11)]{AndrewsCheng2012}. Now write $\epsilon _{t}(\theta )$ $%
=$ $\epsilon _{t}(\psi ,\pi )$, and define gradient processes:
\begin{eqnarray}
\mathcal{G}_{\psi ,n}(\theta ) &=&\sqrt{n}\left\{ \frac{\partial }{\partial
\psi }Q_{n}(\theta )-E\left[ \frac{\partial }{\partial \psi }Q_{n}(\theta )%
\right] \right\} =-\frac{1}{\sqrt{n}}\sum_{t=1}^{n}\left\{ \epsilon
_{t}(\theta )d_{\psi ,t}(\pi )-E\left[ \epsilon _{t}(\theta )d_{\psi ,t}(\pi
)\right] \right\}  \label{GG} \\
&&  \notag \\
\mathcal{G}_{\theta ,n}(\theta ) &=&\mathfrak{B}(\beta _{n})^{-1}\sqrt{n}%
\left\{ \frac{\partial }{\partial \theta }Q_{n}(\theta )-E\left[ \frac{%
\partial }{\partial \theta }Q_{n}(\theta )\right] \right\} =-\frac{1}{\sqrt{n%
}}\sum_{t=1}^{n}\left\{ \epsilon _{t}(\theta )d_{\theta ,t}(\omega (\beta
),\pi )-E\left[ \epsilon _{t}(\theta )d_{\theta ,t}(\omega (\beta ),\pi )%
\right] \right\} .  \notag
\end{eqnarray}%
In order to make $\psi $ $\equiv $ $[\beta ^{\prime },\zeta ^{\prime
}]^{\prime }$ explicit, we write interchangeably%
\begin{equation*}
\mathcal{G}_{\psi ,n}(\psi ,\pi )=\mathcal{G}_{\psi ,n}(\theta )\text{, etc.}
\end{equation*}

Define:%
\begin{eqnarray}
&&\mathfrak{b}_{\psi }(\pi ,\lambda )=E\left[ F\left( \lambda ^{\prime }%
\mathcal{W}(x_{t})\right) d_{\psi ,t}(\pi )\right]  \label{bHK} \\
&&\mathfrak{b}_{\theta }(\omega ,\pi ,\lambda )\equiv E\left[ F\left(
\lambda ^{\prime }\mathcal{W}(x_{t})\right) d_{\theta ,t}(\omega ,\pi )%
\right] \text{ \ and \ }\mathfrak{b}_{\theta }(\lambda )\equiv E\left[
F\left( \lambda ^{\prime }\mathcal{W}(x_{t})\right) d_{\theta ,t}\right]
\notag \\
&&\mathcal{H}_{\psi }(\pi )\equiv E\left[ d_{\psi ,t}(\pi )d_{\psi ,t}(\pi
)^{\prime }\right] \text{ \ and \ }\mathcal{H}_{\theta }(\omega ,\pi )\equiv
E\left[ d_{\theta ,t}(\omega ,\pi )d_{\theta ,t}^{\prime }(\omega ,\pi )%
\right] \text{ and }\mathcal{H}_{\theta }\equiv \mathcal{H}_{\theta }(\omega
_{0},\pi _{0})  \notag \\
&&\mathcal{K}_{\psi ,t}(\pi ,\lambda )\equiv F\left( \lambda ^{\prime
}W(x_{t})\right) -\mathfrak{b}_{\psi }(\pi ,\lambda )^{\prime }\mathcal{H}%
_{\psi }^{-1}(\pi )d_{\psi ,t}(\pi )  \notag \\
&&\mathcal{K}_{\theta ,t}(\lambda )\equiv F\left( \lambda ^{\prime }\mathcal{%
W}(x_{t})\right) -\mathfrak{b}_{\theta }(\lambda )^{\prime }\mathcal{H}%
_{\theta }^{-1}d_{\theta ,t}(\beta _{n}/\left\Vert \beta _{n}\right\Vert
,\pi _{0})\text{ \ and \ }\mathcal{K}_{\theta ,t}(\lambda ;a,m)\equiv
\sum\nolimits_{i=1}^{m}\alpha _{i}\mathcal{K}_{\theta ,t}(\lambda _{i}).
\notag
\end{eqnarray}

\begin{assumption}[data generating process, test weight]
\label{assum:dgp}\ \ \medskip \newline
a. \emph{Identification}:

$(i)$ Under $H_{0}$, $E[\epsilon _{t}|x_{t}]$ $=$ $0$ $a.s.$ and $E[\epsilon
_{t}^{2}|x_{t}]$ $=$ $\sigma _{0}^{2}$ $a.s.$, a finite positive constant$.$

$(ii)$ Under $\mathcal{C}(i,b)$: $E[(y_{t}$ $-$ $\zeta _{0}^{\prime
}x_{t})d_{\psi ,t}(\pi )]$ $=$ $0$ for unique $\psi _{0}$ $=$ $[0_{k_{\beta
}}^{\prime },\zeta _{0}^{\prime }]^{\prime }$ in the interior of $\Psi
^{\ast }$. Under $\mathcal{C}(ii,\omega _{0})$: $E[\epsilon _{t}(\theta
_{0}) $ $\times $ $d_{\theta ,t}(\omega _{0},\pi _{0})]$ $=$ $0$ for unique $%
\theta _{0}$ $=$ $[\beta _{0}^{\prime },\zeta _{0}^{\prime },\pi
_{0}^{\prime }]^{\prime }$ in the interior of $\Theta ^{\ast }=\Psi ^{\ast }$
$\times $ $\Pi ^{\ast }$.\medskip \newline
b. \emph{Memory and Moments}: $\{\epsilon _{t},x_{t}\}$ are $L_{p}$-bounded
for some $p$ $>$ $6$, strictly stationary, and $\beta $-mixing with mixing
coefficients $\beta _{l}=O(l^{-qp(q-p)-\iota })$ for some $q$ $>$ $p$ and
tiny $\iota $ $>$ $0$.\medskip \newline
c. \emph{Response }$g(x,\pi )$\emph{\ and Test Weight }$F(\lambda ^{\prime }%
\mathcal{W}(x))$:

$(i)$ $g(\cdot ,\pi )$ is Borel measurable for each $\pi $; $g(\cdot ,\pi )$
is twice continuously differentiable in $\pi $ $\in $ $\mathbb{R}^{k_{\pi }}$%
; $g(x_{t},\pi )$ is a non-degenerate random variable for each $\pi $ $\in $
$\Pi $.

$(ii)$ $F$ $:$ $\mathbb{R}$ $\rightarrow $ $\mathbb{R}$ is analytic,
non-polynomial, and $\mathcal{W}$ is one-to-one and bounded.

$(iii)$ $E[\sup_{\pi \in \Pi }|(\partial /\partial \pi )^{i}g(x_{t},\pi
)|^{6}]$ $<$ $\infty $ and $E[\sup_{\lambda \in \Lambda }|(\partial
/\partial \lambda )^{j}F(\lambda ^{\prime }\mathcal{W}(x_{t}))|^{6}]$ $<$ $%
\infty $ for $i$ $=$ $0,1,2$ and $j=0,1$.
\end{assumption}

\begin{remark}
\normalfont Condition (a) imposes $E[\epsilon _{t}|x_{t}]$ $=$ $0$ $a.s.$
under the null. It otherwise requires $\zeta _{0}^{\prime }x_{t}+\beta
_{0}^{\prime }g(x_{t},\pi _{0})$ to be a pseudo true representation for $%
y_{t}$ in the sense of being the minimum mean squared error predictor.
Notice under weak identification $\mathcal{C}(i,b)$ we only need to consider
$\zeta _{0}^{\prime }x_{t}$ as a best predictor in some sense. The only
place where $E[\epsilon _{t}^{2}|x_{t}]$ $=$ $\sigma _{0}^{2}$ $a.s$.\ is
used is to simplify the construction of critical values and p-values in
practice. In principle the assumption can be replaced with $E[\epsilon
_{t}^{2}|x_{t}]$ $=$ $\sigma _{t}^{2}(\cdot )$, a known parametric function.
a(ii) identifies parameter values under either hypothesis.
\end{remark}

\begin{remark}
\normalfont$\beta $-mixing under (b) allows us to exploit a probability
inequality due to \cite{Eberlein1984} and arguments in \cite{ArconesYu1994}
in order to prove a stochastic equicontinuity condition.
\end{remark}

\begin{remark}
\normalfont(c.i) ensures measurability. By assuming $g(x_{t},\pi )$ is a
non-degenerate random variable for each $\pi $ $\in $ $\Pi $, it also
focuses the identification issue to $\pi $ alone %
\citep[cf.][STAR(iv)]{AndrewsCheng2012}.\footnote{%
Consider a Logistic STAR model $y_{t}$ $=$ $\zeta _{1}$ $+$ $\zeta
_{2}y_{t-1}$ $+$ $\beta (1$ $+$ $\exp \{\pi _{1}(y_{t-d}$ $-$ $\pi
_{2})\})^{-1}$ $+$ $\epsilon _{t}$ with two covariates $x_{t}$ $=$ $%
[1,y_{t-1}]^{\prime }$, $\pi _{1}$ $>$ $0$
\citep[see,
e.g.,][]{Terasvirta1994}, and assume $E[\epsilon _{t}^{2}]$ $>$ $0$. If $\Pi
^{\ast }$ $=$ $\{\pi $ $\in $ $\mathbb{R}^{2}$ $:$ $\pi _{1}$ $\geq $ $%
\varepsilon )$ for some $\varepsilon $ $>$ $0$ then $g(x_{t},\pi )$ is a
non-degenerate random variable on $\Pi $. Otherwise $g(x,\pi )$ $=$ $1/2$ $%
\forall x$ at $\pi _{1}$ $=$ $0$, the model reduces to a linear structure,
and $\beta $ cannot be identified. We rule out such cases in order to focus
the identification issue to just $\pi $.} The envelope bounds in (c.iii) are
used to prove consistency of criterion derivatives and a sample variance. If
$g(x,\pi )$ $=$ $xh(x,\pi )$ and $h$ and $F$ are exponential, logistic,
normal, or trigonometric then we need only assume $E|x_{t}|^{6}$ $<$ $\infty
$.\footnote{$F(u)$ $=$ $1/(1$ $-$ $u)$ on $(-\infty ,1)$ is analytic and
non-polynomial. However, $F(\lambda ^{\prime }\mathcal{W}(x_{t}))$ need not
satisfy the envelope bound depending on $\mathcal{W}(\cdot )$ and properties
of $x_{t}$.} See, e.g., the proofs of Lemmas \ref{lm:H_ulln_weak} and \ref%
{lm:JV}.
\end{remark}

\begin{remark}
\normalfont We require additional technical details on the existence of
certain long-run variances, and the true and estimation parameter spaces.
Since these are lengthy, and standard, we place them in Assumptions \ref%
{assum:dgp}.d,e,f in Appendix \ref{append:Assum1def}.
\end{remark}

\begin{remark}
\normalfont It is understood that the true parameter space $\Theta ^{\ast }$
contains points consistent with the moment and memory properties of (a)-(e).
\end{remark}

\begin{example}[STAR Model]
\label{Example STAR}\normalfont Model (\ref{model_add}) contains the Smooth
Transition Autoregression [STAR] class where $x_{t}$ $=$ $%
[1,y_{t-1},...,y_{t-p}]^{\prime }$, $g(x_{t},\pi )$ $=$ $x_{t}h(z_{t},\pi )$%
, $z_{t}$ $=$ $y_{t-d}$ for some $1$ $\leq $ $d$ $\leq $ $p$, and $h$ is a
scalar function (typically exponential, logistic, or the normal distribution
function: see \cite{ChanTong1986}, \cite{Luukkonen_etal_1988}, \cite%
{GrangerTerasvirta1993} and \cite{Terasvirta1994}, and \cite{Hill2008}).

We work with a STAR($p$) model with one logistic transition function:
\begin{equation}
y_{t}=\zeta _{0}^{\prime }x_{t}+\beta _{0}^{\prime }x_{t}\frac{1}{1+\exp
\left\{ -\pi _{0}y_{t-1}\right\} }+\epsilon _{t}=f(\theta
_{0},x_{t})+\epsilon _{t}\text{ where }\pi _{0}>0.  \label{lstar}
\end{equation}%
In the STAR literature both $E[\epsilon _{t}|x_{t}]$ $=$ $0$ $a.s.$ and $%
\beta _{0}$ $\neq $ $0$ are simply assumed %
\citep[e.g.][]{Luukkonen_etal_1988,Terasvirta1994,Hill2008,AndrewsCheng2013}%
. However, neither $E[\epsilon _{t}|x_{t}]$ $=$ $0$ $a.s.$ nor $\beta _{0}$ $%
\neq $ $0$ may be true in practice. If $\beta _{0}$ $=$ $0$ then $\pi _{0}$
is not identified: this occurs when a STAR model is estimated, but the true
data generating process is a linear autoregression.\footnote{%
Obviously a pre-test for omitted STAR effects in a linear model may be
performed \citep[e.g.][]{Hill2008}, but rejecting the linear AR null
hypothesis may be an error, and estimating (\ref{lstar}) can have an
unidentified parameter.}

Further, if $E[\epsilon _{t}|x_{t}]$ $=$ $0$ $a.s.$ is false, then by a(ii)
we still require $f(\theta _{0},x_{t})$ to be pseudo-true in the sense that $%
\epsilon _{t}$ satisfies the orthogonality conditions $E[\epsilon _{t}x_{t}]$
$=$ $0$ and $E[\epsilon _{t}x_{t}(1$ $+$ $\exp \left\{ -\pi
_{0}y_{t-1}\right\} )^{-1}]$ $=$ $0$, and additionally under $\mathcal{C}%
(ii,\omega _{0})$:%
\begin{equation}
E\left[ \epsilon _{t}\beta _{0}^{\prime }x_{t}\frac{y_{t-1}\exp \left\{ -\pi
_{0}y_{t-1}\right\} }{\left( 1+\exp \left\{ -\pi _{0}y_{t-1}\right\} \right)
^{2}}\right] =0.  \label{star_a(ii)}
\end{equation}%
Under these conditions $E[(y_{t}$ $-$ $f(\theta ,x_{t}))^{2}]$ is minimized
at some unique $\theta _{0}$ $=$ $[\zeta _{0}^{\prime },\beta _{0}^{\prime
},\pi _{0}^{\prime }]^{\prime }$ $\in $ $\Theta $ and compact $\Theta $. In
this paper, under $H_{1}$ $:$ $\sup_{\theta \in \Theta }P(E[y_{t}|x_{t}]$ $=$
$f(\theta ,x_{t}))$ $<$ $1$ we are agnostic about what the true data
generating process is. Examples are an LSTAR with \emph{two} transition
functions:%
\begin{equation*}
y_{t}=\zeta _{0}^{\prime }x_{t}+\beta _{1,0}^{\prime }x_{t}\frac{1}{1+\exp
\left\{ -\pi _{1,0}y_{t-1}\right\} }+\beta _{2,0}^{\prime }x_{t}\frac{1}{%
1+\exp \left\{ -\pi _{2,0}y_{t-1}\right\} }+\epsilon _{t};
\end{equation*}%
or ESTAR $y_{t}$ $=$ $\zeta _{0}^{\prime }x_{t}$ $+$ $\beta _{0}^{\prime
}x_{t}\exp \{-\pi _{0}y_{t-1}^{2}\}$ $+$ $\epsilon _{t};$ or $y_{t}$ may be
governed by some other class of processes, like a Self-Exciting Threshold
Autoregression $y_{t}$ $=$ $\zeta _{0}^{\prime }x_{t}$ $+$ $\beta
_{0}^{\prime }x_{t}I(y_{t-1}$ $>$ $c)$ $+$ $\epsilon _{t}$, etc.
\end{example}

Next, we discuss the possible limit process for $\hat{\pi}_{n}$ under weak
identification. Define the value of $\psi $ for the non-identification case:%
\begin{equation*}
\psi _{0,n}\equiv \left[ 0_{k_{\beta }}^{\prime },\zeta _{0}^{\prime }\right]
^{\prime }.
\end{equation*}%
By Lemma \ref{lm:G_uclt_weak}, $\mathcal{G}_{\psi ,n}(\psi _{0,n},\pi )$ $=$
$-1/\sqrt{n}\sum_{t=1}^{n}\{\epsilon _{t}(\psi _{0,n})d_{\psi ,t}(\pi )$ $-$
$E[\epsilon _{t}(\psi _{0,n})d_{\psi ,t}(\pi )]\}$ satisfies $\{\mathcal{G}%
_{\psi ,n}(\psi _{0,n},\pi )$ $:$ $\pi $ $\in $ $\Pi \}$ $\Rightarrow ^{\ast
}$ $\{\mathcal{G}_{\psi }(\pi )$ $:$ $\pi $ $\in $ $\Pi $ $\}$, a zero mean
Gaussian process. Define:%
\begin{equation}
\mathcal{D}_{\psi }(\pi )\equiv -\frac{\partial }{\partial \beta
_{0}^{\prime }}E\left[ \epsilon _{t}(\theta )d_{\psi ,t}(\pi )\right] =-E%
\left[ d_{\psi ,t}(\pi )g(x_{t},\pi _{0})^{\prime }\right] ,  \label{D_psi}
\end{equation}%
and processes $\{\xi _{\psi }(\pi ,\cdot )$ $:$ $\pi $ $\in $ $\Pi \}$ and $%
\{\tau _{\beta }(\pi ,\cdot )$ $:$ $\pi $ $\in $ $\Pi \}$:
\begin{eqnarray}
&&\xi _{\psi }(\pi ,b)\equiv -\frac{1}{2}\left\{ \mathcal{G}_{\psi }(\pi )+%
\mathcal{D}_{\psi }(\pi )b\right\} ^{\prime }\mathcal{H}_{\psi }^{-1}(\pi
)\left\{ \mathcal{G}_{\psi }(\pi )+\mathcal{D}_{\psi }(\pi )b\right\}
\label{z_psi} \\
&&\tau _{\beta }(\pi ,b)\equiv -\mathcal{S}_{\beta }\mathcal{H}_{\psi
}^{-1}(\pi )\left\{ \mathcal{G}_{\psi }(\pi )+\mathcal{D}_{\psi }(\pi
)b\right\} \text{ where }\mathcal{S}_{\beta }\equiv \left[ I_{k_{\beta
}}:0_{k_{x}\times k_{x}}\right] .  \notag
\end{eqnarray}%
Notice $\mathcal{S}_{\beta }$ is the $\beta $ selection matrix. $\mathcal{H}%
_{\psi }(\pi )$ is positive definite uniformly on $\Pi $ by Assumption \ref%
{assum:dgp}.d(iii). The following identifies the limit process of $\pi
^{\ast }(b)$ of $\hat{\pi}_{n}$ under weak identification, and ensures $\tau
_{\beta }(\pi ^{\ast }(b),b)$ is not degenerate (cf. Andrews and Cheng 2012:
Assumption C6, and Andrews and Cheng 2013: STAR3). See
\citet[Assumption
C6']{AndrewsCheng2013} for primitive conditions that ensure the next
assumption when $\beta $ is a scalar.

\begin{assumption}[identification of $\protect\pi $]
\label{assum:pi}Let drift case $\mathcal{C}(i,b)$ hold with $||b||$ $<$ $%
\infty $. $(a)$ Each sample path of the process $\{\xi _{\psi }(\pi ,b)$ $:$
$\pi $ $\in $ $\Pi \}$ in some set $\mathfrak{A}(b)$\ with $P(\mathfrak{A}%
(b))$ $=$ $1$ is minimized over $\Pi $ at a unique point $\pi ^{\ast }(b)$
that may depend on the sample path. $(b)$ $P(\tau _{\beta }(\pi ^{\ast
}(b),b)$ $=$ $0)$ $=$ $0$.
\end{assumption}

The test statistic scale $\hat{v}_{n}^{2}(\hat{\theta}_{n},\lambda )$ can
have a degenerate probability limit for some $\lambda $ under strong
identification cases, depending on the variation of $x_{t}$
\citep[see][p.
1447]{Bierens1990}. Indeed, because we include a constant term, at $\lambda $
$=$ $0$\ with least squares $1/\sqrt{n}\sum_{t=1}^{n}\epsilon _{t}(\hat{%
\theta}_{n})F(0^{\prime }\mathcal{W}(x_{t}))$ $=$ $0$, hence $\hat{v}%
_{n}^{2}(\hat{\theta}_{n},0)$ $\overset{p}{\rightarrow }$ $0$. The following
rules this out on $\Lambda $-\emph{a.e. }by effectively ensuring the
parameters $(\alpha ,\beta ,\zeta ,\pi )$ in the augmented model $y_{t}$ $=$
$\zeta ^{\prime }x_{t}+\beta ^{\prime }g(x_{t},\pi )$ $+$ $\alpha \mu
(x_{t}) $ $+$ $u_{t}$ can be locally identified for some Borel measurable $%
\mu $ $:$ $\mathbb{R}^{k_{x}}$ $\rightarrow $ $\mathbb{R}$. See
\citet[p.
1447]{Bierens1990} for discussion, and see Lemma \ref{lm:v2_nondegen} in
Appendix \ref{app:lemmas}. Further details in the context of a STAR model
are provided online in the supplemental material
\citet[Appendix
F]{Supp_Mat_2020}.

\begin{assumption}[non-degenerate scale on $\Lambda $-\emph{a.e.}]
\label{assum:nondegen_v2_ae} $\ \ \ \medskip $\newline
$a.$ Let $\mathcal{C}(i,b)$ with $||b||$ $<$ $\infty $ hold. Then $%
P(E[\inf_{\pi \in \Pi }\{\epsilon _{t}^{2}(\psi _{0},\pi )\}|x_{t}]$ $>$ $0)$
$=$ $1$. There exists a Borel measurable function $\mu $ $:$ $\mathbb{R}%
^{k_{x}}$ $\rightarrow $ $\mathbb{R}$ such that $\kappa _{t}(\omega ,\pi )$ $%
\equiv $ $[\mu (x_{t}),d_{\theta ,t}(\omega ,\pi )^{\prime }]^{\prime }$ has
nonsingular $E[\kappa _{t}(\omega ,\pi )\kappa _{t}(\omega ,\pi )^{\prime }]$
uniformly on $\{\omega $ $\in $ $\mathbb{R}^{k_{\beta }}$ $:$ $\omega
^{\prime }\omega $ $=$ $1\}$ $\times $ $\Pi $.$\medskip $\newline
$b.$ Let $\mathcal{C}(ii,\omega _{0})$ hold. There exists a Borel measurable
function $\mu $ $:$ $\mathbb{R}^{k_{x}}$ $\rightarrow $ $\mathbb{R}$ such
that $\kappa _{t}$ $\equiv $ $[\mu (x_{t}),d_{\theta ,t}]^{\prime }$ has a
nonsingular $E[\kappa _{t}\kappa _{t}^{\prime }]$.
\end{assumption}

\begin{remark}
\normalfont By Lemma 2 in \cite{Bierens1990}, $E[\epsilon _{t}^{2}|x_{t}]$ $%
= $ $\sigma _{0}^{2}$ $>$ $0$ $a.s.$ under Assumption \ref{assum:dgp}.a(i),
and Assumption \ref{assum:nondegen_v2_ae}.b, imply under strong
identification $E[\epsilon _{t}^{2}\{F\left( \lambda ^{\prime }\mathcal{W}%
(x_{t})\right) $ $-$ $\mathfrak{b}_{\theta }(\lambda )^{\prime }\mathcal{H}%
_{\theta }^{-1}d_{\theta ,t}\}^{2}]$ $>$ $0$ on $\Lambda $-\emph{a.e.} since
$F(\lambda ^{\prime }\mathcal{W}(x_{t}))$ $-$ $\mathfrak{b}_{\theta
}(\lambda )^{\prime }\mathcal{H}_{\theta }^{-1}d_{\theta ,t}$ $\neq $ $0$ $%
a.s.$ on $\Lambda $-\emph{a.e}. In Lemma \ref{lm:v2_nondegen} in Appendix %
\ref{app:lemmas} we prove Assumption \ref{assum:nondegen_v2_ae}.a plays a
similar role under weak identification. This suffices to ensure $\hat{v}%
_{n}^{2}(\hat{\theta}_{n},\lambda )$ $>$ $0$ asymptotically with probability
approaching one, on $\Lambda $-\emph{a.e}.
\end{remark}

Unfortunately, Assumption \ref{assum:nondegen_v2_ae} does not rule out $\hat{%
v}_{n}^{2}(\hat{\theta}_{n},\lambda )$ $\overset{p}{\rightarrow }$ $0$
\textit{everywhere} on $\Lambda $. In order to avoid deviant cases, we make
the following assumption which is mild in view of Lemma \ref{lm:v2_nondegen}
\citep[see
also][p. 1449]{Bierens1990}. We discuss the requirement in the supplemental
material \citet[Appendix
F]{Supp_Mat_2020}.

Define the augmented parameter set $\theta ^{+}$ $\equiv $ $[||\beta
||,\omega (\beta )^{\prime },\zeta ^{\prime },\pi ^{\prime }]^{\prime }$
which is useful under weak identification cases when $\beta $ is a vector.
Let $\theta ^{+}$ $\in $ $\Theta ^{+}$ $\equiv $ $\{\theta ^{+}$ $\in $ $%
\mathbb{R}^{k_{\beta }+k_{x}+k_{\pi }+1}$ $:$ $\theta ^{+}$ $=$ $[||\beta
||,\omega (\beta ),\zeta ,\pi ]^{\prime }$ $:$ $\beta $ $\in $ $\mathcal{B},$
$\zeta $ $\in $ $\mathcal{Z}(\beta ),$ $\pi $ $\in $ $\Pi \}$. Define $%
\epsilon _{t}(\theta ^{+})$ $\equiv $ $y_{t}$ $-$ $\zeta ^{\prime }x_{t}$ $-$
$||\beta ||\omega (\beta )^{\prime }g(x_{t},\pi )$ and%
\begin{eqnarray*}
v^{2}(\theta ^{+},\lambda ) &=&E\left[ \epsilon _{t}^{2}(\theta ^{+})\left\{
F\left( \lambda ^{\prime }\mathcal{W}(x_{t})\right) -\mathfrak{b}_{\theta
}(\omega ,\pi ,\lambda )^{\prime }\mathcal{H}_{\theta }^{-1}(\omega ,\pi
)d_{\theta ,t}(\omega ,\pi )\right\} ^{2}\right] \\
\hat{v}_{n}^{2}(\theta ^{+},\lambda ) &=&\frac{1}{n}\sum_{t=1}^{n}\epsilon
_{t}^{2}(\theta ^{+})\left\{ F\left( \lambda ^{\prime }\mathcal{W}%
(x_{t})\right) -\mathfrak{\hat{b}}_{\theta ,n}(\omega ,\pi ,\lambda
)^{\prime }\widehat{\mathcal{H}}_{n}^{-1}d_{\theta ,t}(\omega ,\pi )\right\}
^{2}.
\end{eqnarray*}%
Since $\sup_{\theta ^{+}\in \Theta ^{+},\lambda \in \Lambda }||\hat{v}%
_{n}^{2}(\theta ^{+},\lambda )$ $-$ $v^{2}(\theta ^{+},\lambda )||$ $\overset%
{p}{\rightarrow }$ $0$ by Lemma \ref{lm:vn}, it suffices to bound $%
v^{2}(\theta ^{+},\lambda )$. In the scalar case $\beta $ $\in $ $\mathcal{B}
$ $\subset $ $\mathbb{R}$\ we use $\epsilon _{t}(\theta )$ $\equiv $ $y_{t}$
$-$ $\zeta ^{\prime }x_{t}$ $-$ $\beta ^{\prime }g(x_{t},\pi )$ and%
\begin{eqnarray*}
v^{2}(\theta ,\lambda ) &=&E\left[ \epsilon _{t}^{2}(\theta )\left\{ F\left(
\lambda ^{\prime }\mathcal{W}(x_{t})\right) -\mathfrak{b}_{\theta }(\omega
(\beta ),\pi ,\lambda )^{\prime }\mathcal{H}_{\theta }^{-1}(\omega (\beta
),\pi )d_{\theta ,t}(\omega (\beta ),\pi )\right\} ^{2}\right] \\
\hat{v}_{n}^{2}(\theta ,\lambda ) &=&\frac{1}{n}\sum_{t=1}^{n}\epsilon
_{t}^{2}(\theta )\left\{ F\left( \lambda ^{\prime }\mathcal{W}(x_{t})\right)
-\mathfrak{\hat{b}}_{\theta ,n}(\omega (\beta ),\pi ,\lambda )^{\prime }%
\widehat{\mathcal{H}}_{n}^{-1}d_{\theta ,t}(\omega (\beta ),\pi )\right\}
^{2}.
\end{eqnarray*}%
Since $\sup_{\theta ^{+}\in \Theta ^{+},\lambda \in \Lambda }||\hat{v}%
_{n}^{2}(\theta ^{+},\lambda )$ $-$ $v^{2}(\theta ^{+},\lambda )||$ $\overset%
{p}{\rightarrow }$ $0$ and $\sup_{\theta \in \Theta ^{+},\lambda \in \Lambda
}||\hat{v}_{n}^{2}(\theta ,\lambda )$ $-$ $v^{2}(\theta ,\lambda )||$ $%
\overset{p}{\rightarrow }$ $0$ by Lemma \ref{lm:vn}, it suffices to bound $%
v^{2}(\theta ^{+},\lambda )$ and $v^{2}(\theta ,\lambda )$.

See \citet[Appendix F]{Supp_Mat_2020} for discussion of the following
assumption.

\begin{assumption}[non-degnerate scale]
\label{assum:nondegen_v2_ae_all} $\ \ \ \medskip $\newline
$a.$ Let $\beta $ be a scalar. Let $\inf_{\pi \in \Pi }v^{2}((\beta
_{0},\zeta _{0},\pi ),\lambda )$ $>$ $0$ $\forall \lambda $ $\in $ $\Lambda $
under identification case $\mathcal{C}(i,b)$ with $|b|$ $<$ $\infty $, and
under $\mathcal{C}(ii,\omega _{0})$ let $v^{2}(\theta _{0},\lambda )$ $>$ $0$
$\forall \lambda $ $\in $ $\Lambda $.\medskip \newline
$b.$ Let $\beta $ be a vector. Let $\inf_{\omega \in \mathbb{R}^{k_{\beta
}}:\omega ^{\prime }\omega =1,\pi \in \Pi }v^{2}((||\beta _{0}||,\omega
,\zeta _{0},\pi ),\lambda )$ $>$ $0$ $\forall \lambda $ $\in $ $\Lambda $
under identification case $\mathcal{C}(i,b)$ with $||b||$ $<$ $\infty $, and
under $\mathcal{C}(ii,\omega _{0})$ let $v^{2}(\theta _{0}^{+},\lambda )$ $>$
$0$ $\forall \lambda $ $\in $ $\Lambda $.
\end{assumption}

\section{TEST STATISTIC LIMIT THEORY\label{sec:CM_main}}

We begin by deriving the weak limit of $\hat{\theta}_{n}$, since it strongly
influences the limit properties of $\mathcal{T}_{n}(\lambda )$. Recall $\hat{%
\theta}_{n}$ $=$ $[\hat{\psi}_{n}(\hat{\pi}_{n})^{\prime },\hat{\pi}%
_{n}^{\prime }]^{\prime }$ where $\hat{\psi}_{n}(\pi )$ $=$ $\arginf_{\psi
\in \Psi }Q_{n}(\psi ,\pi )$ and $\hat{\pi}_{n}$ $=$ $\arginf_{\pi \in \Pi
}Q_{n}(\hat{\psi}_{n}(\pi ),\pi ).$

The limit process of a suitably normalized $\hat{\theta}_{n}$ requires the
following constructions. Recall $\psi _{0,n}\equiv \lbrack 0_{k_{\beta
}}^{\prime },\zeta _{0}^{\prime }]$. We need two Gaussian processes $(%
\mathcal{G}_{\psi }(\pi ),\mathcal{G}_{\theta })$ defined as the weak limits
of $\mathcal{G}_{\psi ,n}(\psi _{0,n},\pi )$ and $\sqrt{n}\mathfrak{B}(\beta
_{n})^{-1}(\partial /\partial \theta )Q_{n}(\theta _{n})$, cf. Lemma \ref%
{lm:G_uclt_weak} and Corollary \ref{cor:G_uclt_strong} in Appendix \ref%
{app:lemmas}. Define:%
\begin{equation*}
\tau (\pi ,b)\equiv -\mathcal{H}_{\psi }^{-1}(\pi )\left\{ \mathcal{G}_{\psi
}(\psi _{0,n},\pi )+\mathcal{D}_{\psi }(\pi )b\right\} -\left[ b^{\prime
},0_{k_{\beta }}^{\prime }\right] ^{\prime }.
\end{equation*}%
Denote the true value $\psi _{n}$ $\equiv $ $[\beta _{n}^{\prime },\zeta
_{0}^{\prime }]^{\prime }$ under drifting sequence $\{\beta _{n}\}$.

The following allows for model mis-specification in the sense of Assumption %
\ref{assum:dgp}.a(ii). This effectively generalizes the low level extensions
of \citet[Theorems 3.1 and 3.2]{AndrewsCheng2012} developed in
\citet[Example
2]{AndrewsCheng2012} for the ARMA model with an iid error, and %
\citet[Section 6]{AndrewsCheng2013} for the STAR model with an mds error.
The proof is similar to arguments in \cite{AndrewsCheng2012}: see %
\citet[Appendix C]{Supp_Mat_2020}.

\begin{theorem}
\label{th:ols}Let Assumptions \ref{assum:dgp} and \ref{assum:pi}
hold.\medskip \newline
$a.$ Under drift case $\mathcal{C}(i,b)$ with $||b||$ $<$ $\infty $, $(\sqrt{%
n}(\hat{\psi}_{n}(\hat{\pi}_{n})$ $-$ $\psi _{n}),\hat{\pi}_{n})$ $\overset{d%
}{\rightarrow }$ $(\tau (\pi ^{\ast }(b),b),\pi ^{\ast }(b))$.$\medskip $%
\newline
$b.$ Under drift case $\mathcal{C}(ii,\omega _{0})$, $\sqrt{n}\mathfrak{B}(%
\hat{\beta}_{n})(\hat{\theta}_{n}$ $-$ $\theta _{n})$ $\overset{d}{%
\rightarrow }$ $-\mathcal{H}_{\theta }^{-1}\mathcal{G}_{\theta }$.
\end{theorem}

\begin{remark}
\label{rm:sig_hat}\normalfont Under any degree of (non)identification, $\hat{%
\sigma}_{n}^{2}$ $\equiv $ $1/n\sum_{t=1}^{n}(y_{t}$ $-$ $f(\hat{\theta}%
_{n},x_{t}))^{2}$ $\overset{p}{\rightarrow }$ $\sigma _{0}^{2}$ is easily
verified. The only issue is case $\mathcal{C}(i,b)$, but\ $\beta _{n}$ $%
\rightarrow $ $0$ ensures the non-standard asymptotic properties of $\hat{\pi%
}_{n}$ are irrelevant asymptotically.
\end{remark}

Now turn to the test statistic $\mathcal{T}_{n}(\lambda )$ $=$ $\{1/\sqrt{n}%
\sum_{t=1}^{n}\epsilon _{t}(\hat{\theta}_{n})F(\lambda ^{\prime }\mathcal{W}%
(x_{t}))/\hat{v}_{n}(\hat{\theta}_{n},\lambda )\}^{2}$. The required limit
processes under the null are constructed as follows. First, consider the
weak identification case $\mathcal{C}(i,b)$ with $||b||$ $<$ $\infty $.
Write $\epsilon _{t}(\psi ,\pi )$ $\equiv $ $y_{t}$ $-$ $\zeta ^{\prime
}x_{t}$ $-$ $\beta _{n}^{\prime }g(x_{t},\pi )$. By Lemma \ref{lm:CM_weak}%
.a, under the null we have weak convergence:%
\begin{equation*}
\left\{ \frac{1}{\sqrt{n}}\sum_{t=1}^{n}\epsilon _{t}\left( F\left( \lambda
^{\prime }\mathcal{W}(x_{t})\right) -\mathfrak{b}_{\psi }(\pi ,\lambda
)^{\prime }\mathcal{H}_{\psi }^{-1}(\pi )d_{\psi ,t}(\pi )\right) :\Pi
,\Lambda \right\} \Rightarrow ^{\ast }\left\{ \mathfrak{Z}_{\psi }(\pi
,\lambda ):\Pi ,\Lambda \right\} ,
\end{equation*}%
where $\{\mathfrak{Z}_{\psi }(\pi ,\lambda )$ $:$ $\pi $ $\in $ $\Pi
,\lambda $ $\in $ $\Lambda \}$ is a zero mean Gaussian process with
covariance kernel $\sigma _{0}^{2}E[\mathcal{K}_{\psi ,t}(\pi ,\lambda )$ $%
\times $ $\mathcal{K}_{\psi ,t}(\tilde{\pi},\tilde{\lambda})]$. The
numerator of the test statistic, $(1/\sqrt{n}\sum_{t=1}^{n}\epsilon _{t}(%
\hat{\theta}_{n})F(\lambda ^{\prime }\mathcal{W}(x_{t})))^{2}$, therefore
converges under $H_{0}$ to $\mathfrak{T}_{\psi }^{2}(\pi ^{\ast }(b),\lambda
,b)$, where:%
\begin{eqnarray}
\mathfrak{T}_{\psi }(\pi ,\lambda ,b) &\equiv &\mathfrak{Z}_{\psi }(\pi
,\lambda )+\mathfrak{b}_{\psi }(\pi ,\lambda )^{\prime }\left\{ \mathcal{H}%
_{\psi }^{-1}(\pi )\mathcal{D}_{\psi }(\pi )b+\left[ b,0_{k_{\beta
}}^{\prime }\right] ^{\prime }\right\}  \label{T(pi,lam,b)} \\
&&+\mathfrak{b}_{\psi }(\pi ,\lambda )^{\prime }\mathcal{H}_{\psi }^{-1}(\pi
)E\left[ d_{\psi ,t}(\pi )\left\{ g(x_{t},\pi _{0})-g(x_{t},\pi )\right\}
^{\prime }\right] b  \notag \\
&&+E\left[ \mathcal{K}_{\psi ,t}(\pi ,\lambda )\left\{ g(x_{t},\pi
_{0})-g(x_{t},\pi )\right\} ^{\prime }\right] b.  \notag
\end{eqnarray}%
The form reflects that ($i$) we only expand around the (possibly drifting)
true $\psi _{n}$ because $\hat{\pi}_{n}$ has a nonstandard limit law, hence $%
\mathfrak{Z}_{\psi }(\pi ,\lambda )$; ($ii$) weak identification with $b$ $%
\neq $ $0$ adds asymptotic bias in $\sqrt{n}(\hat{\psi}_{n}(\hat{\pi}_{n})$ $%
-$ $\psi _{n})$ in the first order expansion, hence $\mathfrak{b}_{\psi
}(\pi ,\lambda )^{\prime }\{\mathcal{H}_{\psi }^{-1}(\pi )\mathcal{D}_{\psi
}(\pi )b$ $+$ $[b,0_{k_{\beta }}^{\prime }]^{\prime }\}$; and ($iii$) bias
subsequently arises through $\epsilon _{t}(\psi _{n},\hat{\pi}_{n})$ which
does not have a zero mean in general, hence the remaining two terms. See the
proof of Theorem \ref{th:CM_weak} for details.

The scale $\hat{v}_{n}^{2}(\hat{\theta}_{n},\lambda )$ limit under weak
identification is constructed from:%
\begin{eqnarray}
&&v^{2}(\omega ,\pi ,\lambda )\equiv E\left[ \epsilon _{t}^{2}(\psi _{0},\pi
)\left\{ F\left( \lambda ^{\prime }\mathcal{W}(x_{t})\right) -\mathfrak{b}%
_{\theta }(\omega ,\pi ,\lambda )^{\prime }\mathcal{H}_{\theta }^{-1}(\omega
,\pi )d_{\theta ,t}(\omega ,\pi )\right\} ^{2}\right]  \label{v2} \\
&&\bar{v}^{2}(\pi ,\lambda ,b)\equiv v^{2}(\omega ^{\ast }(\pi ,b),\pi
,\lambda )\text{ where }\omega ^{\ast }(\pi ,b)\equiv \tau _{\beta }(\pi
,b)/\left\Vert \tau _{\beta }(\pi ,b)\right\Vert .  \notag
\end{eqnarray}%
The null limit process of the test statistic under weak identification is
therefore:%
\begin{equation*}
\mathcal{T}_{\psi }(\pi ,\lambda ,b)\equiv \frac{\mathfrak{T}_{\psi
}^{2}(\pi ,\lambda ,b)}{\bar{v}^{2}(\pi ,\lambda ,b)}\text{\ and }\mathcal{T}%
_{\psi }(\lambda ,b)\equiv \mathcal{T}_{\psi }(\pi ^{\ast }(b),\lambda ,b).
\end{equation*}

Now consider strong identification $\mathcal{C}(ii,\omega _{0})$. By Lemma %
\ref{lm:CM_weak}.b we have the weak limit under the null:%
\begin{equation*}
\left\{ \frac{1}{\sqrt{n}}\sum_{t=1}^{n}\epsilon _{t}\left( F\left( \lambda
^{\prime }\mathcal{W}(x_{t})\right) -\mathfrak{b}_{\theta }(\lambda
)^{\prime }\mathcal{H}_{\theta }^{-1}d_{\theta ,t}\right) :\lambda \in
\Lambda \right\} \Rightarrow ^{\ast }\left\{ \mathfrak{Z}_{\theta }(\lambda
):\lambda \in \Lambda \right\} ,
\end{equation*}%
where $\{\mathfrak{Z}_{\theta }(\lambda )$ $:$ $\lambda $ $\in $ $\Lambda \}$
is a zero mean Gaussian process with variance $v^{2}(\lambda )$ $\equiv $ $%
E[\epsilon _{t}^{2}\{F\left( \lambda ^{\prime }\mathcal{W}(x_{t})\right) $ $%
- $ $\mathfrak{b}_{\theta }(\lambda )^{\prime }\mathcal{H}_{\theta
}^{-1}d_{\theta ,t}\}^{2}]$. The limit process $\mathcal{T}(\lambda )$ $%
\equiv $ $\mathfrak{Z}_{\theta }^{2}(\lambda )/v^{2}(\lambda )$ is therefore
chi-squared with one degree of freedom, as in \cite{Bierens1990}, cf. \cite%
{deJong1996} and \cite{Hill2008}.

\begin{theorem}
\label{th:CM_weak}Let Assumptions \ref{assum:dgp}, \ref{assum:pi} and \ref%
{assum:nondegen_v2_ae_all}, and $H_{0}$, hold.\medskip \newline
$a.$ If drift case $\mathcal{C}(i,b)$ holds with $||b||$ $<$ $\infty $, then
$\{\mathcal{T}_{n}(\lambda )$ $:$ $\lambda $ $\in $ $\Lambda \}$ $%
\Rightarrow ^{\ast }$ $\{\mathcal{T}_{\psi }(\lambda ,b)$ $:$ $\lambda $ $%
\in $ $\Lambda \}$. Further, $\inf_{\pi \in \Pi }\bar{v}^{2}(\pi ,\lambda
,b) $ $>$ $0$ $\forall \lambda $ $\in $ $\Lambda $, $\mathcal{T}_{\psi
}(\lambda ,b)$ $\geq $ $0$ $a.s.$, and $\sup_{\lambda \in \Lambda }\{%
\mathcal{T}_{\psi }(\lambda ,b)\}$ $<$ $\infty $ $a.s$.\medskip \newline
$b.$ If drift case $\mathcal{C}(ii,\omega _{0})$ holds then $\{\mathcal{T}%
_{n}(\lambda )$ $:$ $\lambda $ $\in $ $\Lambda \}$ $\Rightarrow ^{\ast }$ $\{%
\mathcal{T}(\lambda )$ $:$ $\lambda $ $\in $ $\Lambda \}$, a chi-squared
process with one degree of freedom, with a version that has \emph{almost
surely} uniformly continuous sample paths. Further, $v^{2}(\lambda )$ $>$ $0$
$\forall \lambda $ $\in $ $\Lambda $, $\mathcal{T}(\lambda )$ $\geq $ $0$ $%
a.s.$, $\sup_{\lambda \in \Lambda }\{\mathcal{T}(\lambda )\}$ $<$ $\infty $ $%
a.s$., and $\mathcal{T}(\lambda )$ has an absolutely continuous distribution
for each $\lambda $.
\end{theorem}

The CM test is consistent on $\Lambda $-\emph{a.e}. under any identification
category. Recall the alternative is non-local to null: $H_{1}$ $:$ $%
\sup_{\theta \in \Theta }P(E[y_{t}|x_{t}]$ $=$ $f(\theta ,x_{t}))$ $<$ $1$.

\begin{theorem}
\label{th:CM_H1}Let Assumptions \ref{assum:dgp}, \ref{assum:pi} and \ref%
{assum:nondegen_v2_ae_all} hold, and let drift case $\mathcal{C}(i,b)$ with $%
||b||$ $<$ $\infty $, or $\mathcal{C}(ii,\omega _{0})$, apply. Under $H_{1}$%
, $\mathcal{T}_{n}(\lambda )$ $\overset{p}{\rightarrow }$ $\infty $ for all $%
\lambda $ $\in $ $\Lambda /S$ where $S$ $\subset $ $\Lambda $ has Lebesgue
measure zero.
\end{theorem}

\section{ROBUST P-VALUE CONSTRUCTIONS\label{sec:CM_pv_construct}}

We now develop identification category robust p-values, while their
bootstrapped approximations are handled in Section \ref{sec:CM_pv}. Critical
value computation is presented in \citet[Appendix E]{Supp_Mat_2020}.

Recall the Theorem \ref{th:CM_weak}.a null limit process $\{\mathcal{T}%
_{\psi }(\lambda ,b)$ $:$ $\lambda $ $\in $ $\Lambda \}$ of $\mathcal{T}%
_{n}(\lambda )$ under weak identification $\sqrt{n}\beta _{n}$ $\rightarrow $
$b$ with $||b||$ $<$ $\infty $. Note $\theta =[\zeta ^{\prime },\beta
^{\prime },\pi ^{\prime }]^{\prime }$ may not fully parameterize the
distribution of $W_{t}$ $\equiv $ $[y_{t},x_{t}^{\prime }]^{\prime }$, hence
$\mathcal{T}_{\psi }(\lambda ,b)$\ may not reveal \textit{all} distribution
based nuisance parameters. Let $\phi _{0}$ index all remaining (nuisance)
parameters such that the distribution of $W_{t}$ is determined by $\gamma
_{0}$ $\equiv $ $(\theta _{0},\phi _{0})$ where
\citep[see][p.
2161-2162]{AndrewsCheng2012}:
\begin{equation}
\gamma _{0}\equiv (\theta _{0},\phi _{0})\in \Gamma ^{\ast }\equiv \left\{
\theta \in \Theta ^{\ast },\phi \in \Phi ^{\ast }(\theta )\right\} .
\label{gamma*}
\end{equation}%
In our nonlinear regression setting (\ref{model_add}), $\phi _{0}$ is a
possibly infinite dimensional parameter that indexes all remaining
characteristics of the error distribution not represented in (\ref{model_add}%
).

Assume $\Phi ^{\ast }(\theta )$ $\subset $ $\Phi ^{\ast }$ $\forall \theta $
$\in $ $\Theta ^{\ast }$, where $\Phi ^{\ast }$ is a compact metric space
with a metric that induces weak convergence for $\{W_{t},W_{t+m}\}$ with
respect to drifting $\gamma $ $\rightarrow $ $\gamma _{0}$
\citep[eq.
(2.3)]{AndrewsCheng2012}.\footnote{%
Let $d_{\phi }(\cdot ,\cdot )$ be the metric on $\Phi ^{\ast }$. Under these
assumptions $\Gamma ^{\ast }$ is a metric space with metric $d_{\Gamma
}(\gamma _{1},\gamma _{2})$ $\equiv $ $||\theta _{1}$ $-$ $\theta _{2}||$ $+$
$d_{\phi }(\phi _{1},\phi _{1})$. Now let $P_{\gamma }$ denote the joint
probability function of $(W_{t},W_{t+m})$ induced by the measure $\mathcal{P}
$ under $\gamma $. The metric $d_{\Gamma }(\gamma _{1},\gamma _{2})$
therefore satisfies the following weak convergence: if $\gamma $ $%
\rightarrow $ $\gamma _{0}$ then $P_{\gamma }(W_{t}$ $\leq $ $a,W_{t+m}$ $%
\leq $ $b)$ $\rightarrow $ $P_{\gamma _{0}}(W_{t}$ $\leq $ $a,W_{t+m}$ $\leq
$ $b)$ $\forall a,b$ $\in $ $\mathbb{R}^{p+1}$, $\forall t$, $\forall m$ $%
\geq $ $1$. A key implicit use of these ideas arises under local drift $%
\gamma _{n}$ $\rightarrow $ $\gamma _{0}$, in which case only $\beta $
matters under our assumptions, hence $d_{\Gamma }(\gamma _{n},\gamma _{0})$ $%
\equiv $ $||\beta _{n}$ $-$ $\beta _{0}||$. In the sequel we therefore do
not explicitly state such weak convergence in order to avoid redundancy. See
also \citet[Footnote 21]{AndrewsCheng2012}.} The space $\Phi ^{\ast }(\theta
)$ generally depends on $\theta $ $\in $ $\Theta ^{\ast }$ because we
implicitly assume $\Gamma ^{\ast }$ indexes only those distributions that
satisfy the maintained assumptions, which ensures uniform asymptotics.%
\footnote{%
See \citet[p.
119]{AndrewsCheng2012_supp} for an example.} Under drift $\theta _{n}$ $=$ $%
[\zeta _{0}^{\prime },\beta _{n}^{\prime },\pi _{0}^{\prime }]^{\prime }$
the parameter set becomes $\gamma _{n}$ $\equiv $ $(\theta _{n},\phi _{0})$ $%
\rightarrow $ $\gamma _{0}$. We only let $\beta $ exhibit drift to ease
notation, and since that parameter governs identification cases for $\pi $.

Now define the total parameter set that characterizes data generating
processes under weak identification $\beta _{n}$ $\rightarrow $ $\beta _{0}$
$=$ $0$, and $\sqrt{n}\beta _{n}$ $\rightarrow b$ with $||b||$ $<$ $\infty $%
:
\begin{equation}
h\equiv \left( \gamma _{0},b\right) \in \mathfrak{H}\equiv \left\{ h:\gamma
_{0}\in \Gamma ^{\ast }\text{, and }\left\Vert b\right\Vert <\infty \text{,
with }\beta _{0}=0\right\} .  \label{hH}
\end{equation}

\subsection{P-VALUES FOR $\mathcal{T}_{n}(\protect\lambda )$}

Operate under $H_{0}$. Define $\mathcal{F}_{\infty }(c)$ $\equiv $ $P(%
\mathcal{T}(\lambda )$ $\leq $ $c)$ where $\{\mathcal{T(}\lambda \mathcal{)}$
$:$ $\lambda $ $\in $ $\Lambda \mathcal{\}}$ is the asymptotic null
chi-squared process under strong identification. Similarly, $\mathcal{F}%
_{\lambda ,h}(c)$ $\equiv $ $P(\mathcal{T}_{\psi }(\lambda ,h)$ $\leq $ $c)$
where $\{\mathcal{T}_{\psi }(\lambda ,h)$ $:$ $\lambda $ $\in $ $\Lambda \}$
is the asymptotic null process under weak identification, and we now reveal
all nuisance parameters $h$. The case specific asymptotic p-values are%
\begin{equation*}
p_{n}^{\infty }(\lambda )\equiv 1-\mathcal{F}_{\infty }(\mathcal{T}%
_{n}(\lambda ))=\mathcal{\bar{F}}_{\infty }(\mathcal{T}_{n}(\lambda ))\text{
\ and \ }p_{n}(\lambda ,h)\equiv 1-\mathcal{F}_{\lambda ,h}(\mathcal{T}%
_{n}(\lambda ))=\mathcal{\bar{F}}_{\lambda ,h}(\mathcal{T}_{n}(\lambda )).
\end{equation*}

The following summarizes and extends ideas developed in
\citet[Section
5]{AndrewsCheng2012}. The\textit{\ Least Favorable} [LF] p-value is defined
as $p_{n}^{(LF)}(\lambda )$ $\equiv $ $\max \{\sup_{h\in \mathfrak{H}%
}\{p_{n}(\lambda ,h)\},p_{n}^{\infty }(\lambda )\}$. A better p-value in
terms of power uses the fact that $(\zeta _{0},\beta _{n},\sigma _{0}^{2})$
are consistently estimated by $(\hat{\zeta}_{n},\hat{\beta}_{n},\hat{\sigma}%
_{n}^{2})$ under any degree of (non)identification. The\textit{\ plug-in} LF
p-value $\hat{p}_{n}^{(LF)}(\lambda )$ uses $\widehat{\mathfrak{H}}$ $\equiv
$ $\{h$ $\in $ $\mathfrak{H}$ $:$ $\theta $ $=$ $[\hat{\zeta}_{n}^{\prime },%
\hat{\beta}_{n}^{^{\prime }},\pi ^{\prime }]^{\prime },\sigma ^{2}$ $=$ $%
\hat{\sigma}_{n}^{2}\}$ in place of $\mathfrak{H}$.\footnote{%
The null hypothesis is tested by using a sample version of $E[\epsilon
_{t}F(\lambda ^{\prime }\mathcal{W}(x_{t}))]$. Thus, so-called parametric
\textit{null imposed} p-values, similar to null imposed critical values in
\cite{AndrewsCheng2012} for t-, Quasi-Likelihood Ratio and Wald statistics,
do not play a role here.}

The LF p-value does not exploit information that may point toward a
particular identification case. The \textit{identification category selection%
} [ICS] procedure uses the sample to choose between $\sqrt{n}\beta _{n}$ $%
\rightarrow $ $b$ when $||b||$ $<$ $\infty $ (weak and non-identification)
and $||b||$ $=$ $\infty $ (semi-strong and strong identification). The
statistic used to determine whether $b$ is finite is%
\begin{equation}
\mathcal{A}_{n}\equiv \left( \frac{1}{k_{\beta }}n\hat{\beta}_{n}^{\prime }%
\hat{\Sigma}_{\beta ,\beta ,n}^{-1}\hat{\beta}_{n}\right) ^{1/2}
\label{A_ICS}
\end{equation}%
where $\hat{\Sigma}_{\beta ,\beta ,n}$ is the upper $k_{\beta }$ $\times $ $%
k_{\beta }$ block of $\hat{\Sigma}_{n}$ $\equiv $ $\widehat{\mathcal{H}}%
_{n}^{-1}\mathcal{\hat{V}}_{n}\widehat{\mathcal{H}}_{n}^{-1}$, and
\begin{equation}
\mathcal{\hat{V}}_{n}\equiv \frac{1}{n}\sum_{t=1}^{n}\epsilon _{t}^{2}(\hat{%
\theta}_{n})d_{\theta ,t}(\omega (\hat{\beta}_{n}),\hat{\pi}_{n})d_{\theta
,t}(\omega (\hat{\beta}_{n}),\hat{\pi}_{n})^{\prime }.  \label{JV_hat}
\end{equation}%
Now let $\{\kappa _{n}\}$ be a sequence of positive constants, with $\kappa
_{n}$ $\rightarrow $ $\infty $ and $\kappa _{n}$ $=$ $o(\sqrt{n})$. The case
$||b||$ $<$ $\infty $ is selected when $\mathcal{A}_{n}$ $\leq $ $\kappa
_{n} $, else $||b||$ $=$ $\infty $ is selected. The type 1 ICS [ICS-1]
p-value is:
\begin{equation*}
p_{n}^{(ICS-1)}(\lambda )=\left\{
\begin{array}{ll}
p_{n}^{(LF)}(\lambda )\text{ if }\mathcal{A}_{n}\leq \kappa _{n}, &
p_{n}^{\infty }(\lambda )\text{ if }\mathcal{A}_{n}>\kappa _{n}%
\end{array}%
\right. .
\end{equation*}%
A plug-in version is similar. Only when $\sqrt{n}||\beta _{n}||$ $%
\rightarrow $ $\infty $ faster than $\kappa _{n}$ $\rightarrow $ $\infty $
will the chi-squared based p-value be chosen asymptotically with probability
approaching one since then $\mathcal{A}_{n}/\kappa _{n}$ $\overset{p}{%
\rightarrow }$ $\infty $. Thus, a high bar must be passed in order to select
the strong identification case. In every other case the LF value is chosen,
which is always asymptotically correct.

The type 2 ICS [ICS-2] p-value involves a subtler comparison for category
selection, cf. \citet[Section 5.3]{AndrewsCheng2012}. Since our simulation
study focuses on LF and ICS-1 p-values due to the added computational
complexity of ICS-2 p-values, and ICS-1 works well, we relegate ICS-2
details to the supplemental material \citet[Appendix D]{Supp_Mat_2020}.

A limit theory for the ICS-1 p-value requires the limit distribution of $%
\mathcal{A}_{n}$. It is easily derived along the lines of Theorems \ref%
{th:ols} and \ref{th:CM_weak}. Recall the augmented parameter set $\theta
^{+}$ $\equiv $ $[||\beta ||,\omega (\beta )^{\prime },\zeta ^{\prime },\pi
^{\prime }]^{\prime }$. Define $\mathcal{H}_{\theta }(\theta ^{+})$ $=$ $%
E[d_{\theta ,t}(\omega (\beta ),\pi )d_{\theta ,t}(\omega (\beta ),\pi
)^{\prime }]$ and $\mathcal{V}(\theta ^{+})$ $=$ $E[\epsilon _{t}^{2}(\theta
^{+})d_{\theta ,t}(\omega (\beta ),\pi )$ $\times $ $d_{\theta ,t}(\omega
(\beta ),\pi )^{\prime }]$, and%
\begin{eqnarray}
&&\Sigma (\theta ^{+})\equiv \mathcal{H}_{\theta }(\theta ^{+})^{-1}\mathcal{%
V}(\theta ^{+})\mathcal{H}_{\theta }(\theta ^{+})^{-1}\text{, }\Sigma \left(
\omega ,\pi \right) \equiv \Sigma \left( \left\Vert \beta _{0}\right\Vert
,\omega ,\zeta _{0},\pi \right)  \label{Sigma_the+} \\
&&\bar{\Sigma}(\pi ,b)=\left[ \bar{\Sigma}_{i,j}(\pi ,b)\right]
_{i,j=1}^{k_{\theta }}\equiv \Sigma \left( \omega ^{\ast }(\pi ,b),\pi
\right) .  \notag
\end{eqnarray}

\begin{theorem}
\label{th:ICS_stat}Let Assumptions \ref{assum:dgp} and \ref{assum:pi}, and $%
H_{0}$, hold.$\medskip $\newline
$a.$ Under drift case $\mathcal{C}(i,b)$ with $||b||$ $<$ $\infty $, $%
\mathcal{A}_{n}\overset{d}{\rightarrow }\mathcal{A}(b)$ $\mathcal{\equiv }$ $%
\{\tau _{\beta }(\pi ^{\ast }(b),b)^{\prime }$ $\times $ $\bar{\Sigma}%
_{\beta ,\beta }^{-1}(\pi ^{\ast }(b),b)$ $\times $ $\tau _{\beta }(\pi
^{\ast }(b),b)/k_{\beta }\}^{1/2},$ where $\bar{\Sigma}_{\beta ,\beta }(\pi
,b)\equiv \lbrack \bar{\Sigma}_{i,j}(\pi ,b)]_{i,j=1}^{k_{\beta }}$.$%
\medskip $\newline
$b.$ Let $\{\kappa _{n}\}$\ be a sequence of positive constants, $\kappa
_{n} $ $\rightarrow $ $\infty $ and $\kappa _{n}$ $=$ $o(\sqrt{n})$. Under
drift case $\mathcal{C}(ii,\omega _{0})$ we have $\mathcal{A}_{n}$ $\overset{%
p}{\rightarrow }$ $\infty $. If $\sqrt{n}||\beta _{n}||/\kappa _{n}$ $=$ $%
O(1)$ then $\kappa _{n}^{-1}\mathcal{A}_{n}$ $\overset{p}{\rightarrow }$ $%
[0,\infty )$. If $\sqrt{n}||\beta _{n}||/\kappa _{n}$ $\rightarrow $ $\infty
$ then $\kappa _{n}^{-1}\mathcal{A}_{n}$ $\overset{p}{\rightarrow }$ $\infty
$, for example when $\beta _{0}$ $\neq $ $0$\ for any sequence $\{\kappa
_{n}\}$ defined above.
\end{theorem}

\begin{remark}
\normalfont Intuitively, $\kappa _{n}^{-1}\mathcal{A}_{n}$ $\overset{p}{%
\rightarrow }$ $\infty $ (and therefore the ICS-1 p-value is based on the
chi-squared distribution) only when there is strong evidence in favor of
strong identification. If $\sqrt{n}||\beta _{n}||$ $\rightarrow $ $\infty $
too slowly, in this case $\sqrt{n}||\beta _{n}||/\kappa _{n}$ $=$ $O(1)$,
then the LF value is selected, which leads always to asymptotically correct
inference.
\end{remark}

\section{P-VALUE ASYMPTOTICS AND COMPUTATION\label{sec:CM_pv}}

Let $p_{n}^{(\cdot )}(\lambda )$\ be the LF or ICS-1 p-value. We first prove
that $p_{n}^{(\cdot )}(\lambda )$\ leads to a test with correct asymptotic
level. Analogous results carry over to plug-in versions. We then show how to
bootstrap the key component $p_{n}(\lambda ,h)$ $\equiv $ $1$ $-$ $\mathcal{F%
}_{\lambda ,h}(\mathcal{T}_{n}(\lambda ))$ where $\mathcal{F}_{\lambda
,h}(c) $ $\equiv $ $P(\mathcal{T}_{\psi }(\lambda ,h)$ $\leq $ $c)$, with
accompanying limit theory. Once we have an asymptotically valid
approximation for $p_{n}(\lambda ,h)$, a robust p-value follows as in
Section \ref{sec:CM_pv_construct}. The same method leads to robust critical
value approximations \citep[see][Appendix
E]{Supp_Mat_2020}.

\subsection{ASYMPOTICS $p_{n}^{(\cdot )}(\protect\lambda )$\label%
{sec:CM_pv_asym}}

Technical arguments are made feasible when $\mathcal{F}_{\lambda ,h}(c)$ $%
\equiv $ $P(\mathcal{T}_{\psi }(\lambda ,h)$ $\leq $ $c)$ is continuous,
because $p_{n}^{(\cdot )}(\lambda )$ contains $\mathcal{F}_{\lambda ,h}(%
\mathcal{T}_{n}(\lambda ))$ which is evaluated by a weak limit theory and
the continuous mapping theorem. See \citet[Appendix F]{Supp_Mat_2020} for
discussion.

\begin{assumption}[p-value]
\label{assum:pvs}$a.$ $\mathcal{F}_{\lambda ,h}(\cdot )$ is continuous \emph{%
a.e.} on $[0,\infty )$, $\forall h$ $\in $ $\mathfrak{H}$.\footnote{%
The bulk of Assumptions LF and V3 in \cite{AndrewsCheng2012} ensure
distribution continuity at non-random critical value points for each $\alpha
$. Since these must hold for \textit{any} nominal level $\alpha $,
Assumption \ref{assum:pvs} is not restrictive by comparison.} $b.$ The ICS-1
threshold sequence $\{\kappa _{n}\}$ satisfies $\kappa _{n}$ $\rightarrow $ $%
\infty $ and $\kappa _{n}$ $=$ $o(\sqrt{n})$.
\end{assumption}

Let $F_{\gamma }$ be the distribution function of $W_{t}$ $=$ $%
[y_{t},x_{t}^{\prime }]^{\prime }$ under some $\gamma $ $\in $ $\Gamma
^{\ast }$, where $\Gamma ^{\ast }$ is the true parameter space in (\ref%
{gamma*}). Let $P_{\gamma }$ denote probability under $F_{\gamma }$. For any
p-value $p_{n}^{(\cdot )}(\lambda )$ and each $\lambda $, the asymptotic
size of the test is the asymptotic maximum rejection probability over $%
\gamma $ such that the null is true: $AsySz(\lambda )$ $=$ $%
\limsup_{n\rightarrow \infty }\sup_{\gamma \in \Gamma ^{\ast }}P_{\gamma
}(p_{n}^{(\cdot )}(\lambda )$ $<$ $\alpha |H_{0})$. Uniform size control
over $\lambda $ is captured by $AsySz$ $\equiv $ $\sup_{\lambda \in \Lambda
}AsySz(\lambda )$.

\begin{theorem}
\label{th:pv_weak}Let Assumptions \ref{assum:dgp}, \ref{assum:pi}, \ref%
{assum:nondegen_v2_ae_all} and \ref{assum:pvs} hold.$\medskip $\newline
$a$. LF and ICS-1 $p_{n}^{(\cdot )}(\lambda )$ satisfy $AsySz$ $\leq $ $%
\alpha $.$\medskip $\newline
$b$. Let $H_{1}$ $:$ $\sup_{\theta \in \Theta }P(E[y_{t}|x_{t}]$ $=$ $%
f(\theta ,x_{t}))$ $<$ $1$ be true. Then $p_{n}^{(\cdot )}(\lambda )$ $%
\overset{p}{\rightarrow }$ $0$ for all $\lambda $ $\in $ $\Lambda /S$ where $%
S$ $\subset $ $\Lambda $ has Lebesgue measure zero.
\end{theorem}

\begin{remark}
\normalfont Under Assumptions K, LF, and V3 in \cite{AndrewsCheng2012}, the
robust ICS-1 critical value leads to a correctly sized test. The assumptions
primarily concern continuity of $\mathcal{F}_{\lambda ,h}$ at the critical
value for a given level $\alpha $. In turn these allow for $AsySz$ to be
reduced by their Lemma 2.1, a key step toward proving $AsySz$ $=$ $\alpha $.
The robust p-value requires that $\mathcal{F}_{\lambda ,h}(c)$ $\equiv $ $P(%
\mathcal{T}_{\psi }(\lambda ,h)$ $\leq $ $c)$ be continuous everywhere, and
force us to exploit probability bounds rather than their Lemma 2.1. Hence,
we can only prove that $p_{n}^{(\cdot )}(\lambda )$\ yields a correct
uniform asymptotic \emph{level} $AsySz$ $\leq $ $\alpha $. That seems
irrelevant in small sample experiments since correct size appears to be
achieved: see Section \ref{sec:sim}.
\end{remark}

\subsection{COMPUTATION OF $p_{n}^{(\cdot )}(\protect\lambda )$\label%
{sec:CM_pv_comp}}

We propose a wild bootstrap method for computing $p_{n}(\lambda ,h)$ $\equiv
$ $1$ $-$ $\mathcal{F}_{\lambda ,h}(\mathcal{T}_{n}(\lambda ))$. A similar
method applies to bootstrapping the chi-squared based p-value $p_{n}^{\infty
}(\lambda )$ $\equiv $ $1$ $-$ $\mathcal{F}_{\infty }(\mathcal{T}%
_{n}(\lambda ))$ under strong identification \citep[see, e.g.,][]{Hansen1996}%
. This may be an attractive option in small samples where the chi-squared
distribution may not well approximate the small sample distribution of $%
\mathcal{T}_{n}(\lambda )$\ under strong identification.

Operate under $H_{0}$, and under weak identification $\sqrt{n}\beta _{n}$ $%
\rightarrow $ $b$ and $||b||$ $<$ $\infty $. Set for the sake of brevity $%
\beta _{n}$ $=$ $b/\sqrt{n}$ where $||b||$ $<$ $\infty $ indexes the true
value $\beta _{n}$. We first give the steps for computing $p_{n}(\lambda ,h)$%
, and then prove its validity. In all that follows, independence is
conditional on the sample.

\subsubsection*{Step 1: Compute components $\mathcal{H}_{\protect\psi }(%
\protect\pi )$, $\mathcal{D}_{\protect\psi }(\protect\pi )$, etc.}

Define $\epsilon _{t}(\psi ,\pi )$ $\equiv $ $y_{t}$ $-$ $\zeta ^{\prime
}x_{t}$ $-$ $\beta ^{\prime }g(x_{t},\pi )$ and estimators:%
\begin{eqnarray*}
&&\widehat{\mathcal{H}}_{\psi ,n}(\pi )\equiv \frac{1}{n}\sum_{t=1}^{n}d_{%
\psi ,t}(\pi )d_{\psi ,t}(\pi )^{\prime }\text{ and }\widehat{\mathcal{H}}%
_{n}(\omega ,\pi )=\frac{1}{n}\sum_{t=1}^{n}d_{\theta ,t}(\omega ,\pi
)d_{\theta ,t}(\omega ,\pi )^{\prime } \\
&&\mathcal{\hat{D}}_{\psi ,n}(\pi ,\pi _{0})\equiv -\frac{1}{n}%
\sum_{t=1}^{n}d_{\psi ,t}(\pi )g(x_{t},\pi _{0})^{\prime }\text{ and }%
\mathcal{\hat{K}}_{\psi ,n,t}(\pi ,\lambda )\equiv F\left( \lambda ^{\prime }%
\mathcal{W}(x_{t})\right) -\mathfrak{\hat{b}}_{\psi ,n}(\pi ,\lambda
)^{\prime }\widehat{\mathcal{H}}_{\psi ,n}^{-1}(\pi )d_{\psi ,t}(\pi ) \\
&&\mathfrak{\hat{b}}_{\psi ,n}(\pi ,\lambda )\equiv \frac{1}{n}%
\sum_{t=1}^{n}F\left( \lambda ^{\prime }\mathcal{W}(x_{t})\right) d_{\psi
,t}(\pi )\text{ and }\mathfrak{\hat{b}}_{\theta ,n}(\omega ,\pi ,\lambda
)\equiv \frac{1}{n}\sum_{t=1}^{n}F\left( \lambda ^{\prime }\mathcal{W}%
(x_{t})\right) d_{\theta ,t}(\omega ,\pi ).
\end{eqnarray*}

\subsubsection*{Step 2: Draw from $\protect\pi ^{\ast }(b)$}

By Assumption \ref{assum:pi}, $\pi ^{\ast }(b)$ $\equiv $ $\arg \inf_{\pi
\in \Pi }\xi _{\psi }(\pi ,b)$ $\equiv $ $-\arg \inf_{\pi \in \Pi }\{%
\mathcal{S}_{\beta }\mathcal{H}_{\psi }^{-1}(\pi )(\mathcal{G}_{\psi }(\pi )$
$+$ $\mathcal{D}_{\psi }(\pi )b)\}$. Under weak identification, Lemma \ref%
{lm:G_uclt_weak} yields that $\{\mathcal{G}_{\psi }(\pi )$ $:$ $\pi $ $\in $
$\Pi \}$ is the weak limit of%
\begin{eqnarray}
\mathcal{G}_{\psi ,n}(\psi _{0,n},\pi ) &=&\frac{1}{\sqrt{n}}%
\sum_{t=1}^{n}\left\{ \epsilon _{t}(\psi _{0,n},\pi )d_{\psi ,t}(\pi )-E%
\left[ \epsilon _{t}(\psi _{0,n},\pi )d_{\psi ,t}(\pi )\right] \right\}
\label{G_psi_n} \\
&=&\frac{1}{\sqrt{n}}\sum_{t=1}^{n}\epsilon _{t}d_{\psi ,t}(\pi )+\frac{1}{n}%
\sum_{t=1}^{n}\left\{ d_{\psi ,t}(\pi )g(x_{t},\pi _{0})^{\prime }-E\left[
d_{\psi ,t}(\pi )g(x_{t},\pi _{0})^{\prime }\right] \right\} \times b.
\notag
\end{eqnarray}%
By the argument used to prove Lemma \ref{lm:H_ulln_weak}, $\sup_{\pi \in \Pi
}|1/n\sum_{t=1}^{n}\{d_{\psi ,t}(\pi )g(x_{t},\pi _{0})^{\prime }$ $-$ $%
E[d_{\psi ,t}(\pi )g(x_{t},\pi _{0})^{\prime }]\}|$ $\overset{p}{\rightarrow
}$ $0$. Hence, by $E[\epsilon _{t}|x_{t}]$ $=$ $0$ $a.s.$ and $E[\epsilon
_{t}^{2}|x_{t}]$ $=$ $\sigma _{0}^{2}$ $\in $ $(0,\infty )$ $a.s.$ under $%
H_{0}$, the covariance for $\mathcal{G}_{\psi }(\pi )$ is
\begin{equation*}
E\left[ \epsilon _{t}^{2}d_{\psi ,t}(\pi )d_{\psi ,t}(\tilde{\pi})^{\prime }%
\right] =E\left[ \epsilon _{t}^{2}\right] \times E\left[ d_{\psi ,t}(\pi
)d_{\psi ,t}(\tilde{\pi})^{\prime }\right] =\sigma _{0}^{2}\times \mathcal{H}%
_{\psi }\left( \pi ,\tilde{\pi}\right) ,
\end{equation*}%
say. Thus $\mathcal{H}_{\psi }^{-1/2}(\pi )\mathcal{G}_{\psi }(\pi )$ is
distributed $N(0,\sigma _{0}^{2})$ with kernel $\sigma _{0}^{2}\mathcal{H}%
_{\psi }^{-1/2}(\pi )$ $\times $ $\mathcal{H}_{\psi }(\pi ,\tilde{\pi})$ $%
\times $ $\mathcal{H}_{\psi }^{-1/2}(\pi )$.

Next, let $\{z_{t}\}_{t=1}^{n}$ be a sequence of independent draws from $%
N(0,1)$, and define $\widehat{\mathcal{G}}_{\psi ,n}^{\ast }(\pi )$ $\equiv $
\linebreak $1/\sqrt{n}\sum_{t=1}^{n}z_{t}\widehat{\mathcal{H}}_{\psi
,n}^{-1/2}(\pi )d_{\psi ,t}(\pi )$. By $\hat{\sigma}_{n}$ $\overset{p}{%
\rightarrow }$ $\sigma _{0}$ and the proof of Theorem \ref{th:boot_pv},
below, $\{\hat{\sigma}_{n}\widehat{\mathcal{G}}_{\psi ,n}^{\ast }(\pi )$ $:$
$\pi $ $\in $ $\Pi \}$ $\mathcal{\Rightarrow }^{p}$ $\{\mathcal{H}_{\psi
}^{-1/2}(\pi )\mathcal{G}_{\psi }(\pi )$ $:$ $\pi $ $\in $ $\Pi \},$ where $%
\mathcal{\Rightarrow }^{p}$ denotes weak convergence in probability defined
in \citet[Section
3]{GineZinn1990}.\footnote{\cite{GineZinn1990} work under weak convergence
in $l_{\infty }$ as in \citet{HoffJorg1991}, which is the same rubric of
weak convergence that we work with. Thus, for example, $\{\sigma \widehat{%
\mathcal{G}}_{\psi ,n}^{\ast }(\pi )$ $:$ $\pi $ $\in $ $\Pi \}$ $\mathcal{%
\Rightarrow }^{p}$ $\{\mathcal{H}_{\psi }^{-1/2}(\pi )\mathcal{G}_{\psi
}(\pi )$ $:$ $\pi $ $\in $ $\Pi \}$ \textit{if and only if} $\{\sigma
\widehat{\mathcal{G}}_{\psi ,n}^{\ast }(\pi )$ $:$ $\pi $ $\in $ $\Pi \}$ $%
\mathcal{\Rightarrow }^{\ast }$ $\{\mathcal{H}_{\psi }^{-1/2}(\pi )\mathcal{G%
}_{\psi }(\pi )$ $:$ $\pi $ $\in $ $\Pi \}$ asymptotically with probability
approaching one with respect to the sample draw.} Thus, $\hat{\sigma}_{n}%
\widehat{\mathcal{G}}_{\psi ,n}^{\ast }(\pi )$ is a draw from $\{\mathcal{H}%
_{\psi }^{-1/2}(\pi )\mathcal{G}_{\psi }(\pi )$ $:$ $\pi $ $\in $ $\Pi \}$
with probability approaching one as $n$ $\rightarrow $ $\infty $.

Now use $\{\widehat{\mathcal{G}}_{\psi ,n}^{\ast }(\pi ),\widehat{\mathcal{H}%
}_{\psi ,n}(\pi ),\mathcal{\hat{D}}_{\psi ,n}(\pi ,\pi _{0})\}$ to compute

\begin{equation*}
\hat{\xi}_{\psi ,n}^{\ast }(\pi ,\pi _{0},b)=-\frac{1}{2}\left\{ \hat{\sigma}%
_{n}\widehat{\mathcal{G}}_{\psi ,n}^{\ast }(\pi )+\widehat{\mathcal{H}}%
_{\psi ,n}^{-1/2}(\pi )\mathcal{\hat{D}}_{\psi ,n}(\pi ,\pi _{0})\times
b\right\} ^{\prime }\left\{ \hat{\sigma}_{n}\widehat{\mathcal{G}}_{\psi
,n}^{\ast }(\pi )+\widehat{\mathcal{H}}_{\psi ,n}^{-1/2}(\pi )\mathcal{\hat{D%
}}_{\psi ,n}(\pi ,\pi _{0})\times b\right\} .
\end{equation*}%
The bootstrapped $\pi ^{\ast }(b)$ is therefore:
\begin{equation}
\hat{\pi}_{n}^{\ast }(\pi _{0},b)=\argmin_{\pi \in \Pi }\left\{ \hat{\xi}%
_{\psi ,n}^{\ast }(\pi ,\pi _{0},b)\right\} .  \label{pi*}
\end{equation}

\subsubsection*{Step 3: Draw from $\mathfrak{Z}_{\protect\psi }(\protect\pi ,%
\protect\lambda )$}

By Lemma \ref{lm:CM_weak}.a, under the null $\mathfrak{Z}_{\psi }(\pi
,\lambda )$ is the zero mean Gaussian limit process of $1/\sqrt{n}%
\sum_{t=1}^{n}\epsilon _{t}\mathcal{K}_{\psi ,t}(\pi ,\lambda )$. Use the
Step 2 draws $\{z_{t}\}_{t=1}^{n}$ to define%
\begin{equation}
\mathfrak{\hat{Z}}_{\psi ,n}^{\ast }(\pi ,\lambda )\equiv \frac{1}{\sqrt{n}}%
\sum_{t=1}^{n}z_{t}\left( F\left( \lambda ^{\prime }\mathcal{W}%
(x_{t})\right) -\mathfrak{\hat{b}}_{\psi ,n}(\pi ,\lambda )^{\prime }%
\widehat{\mathcal{H}}_{\psi ,n}^{-1}(\pi )d_{\psi ,t}(\pi )\right) .
\label{Z*}
\end{equation}%
Then $\{\hat{\sigma}_{n}\mathfrak{\hat{Z}}_{\psi ,n}^{\ast }(\pi ,\lambda )$
$:$ $\pi $ $\in $ $\Pi ,\lambda $ $\in $ $\Lambda \}$ $\mathcal{\Rightarrow }%
^{p}$ $\{\mathfrak{Z}_{\psi }(\pi ,\lambda )$ $:$ $\pi $ $\in $ $\Pi
,\lambda $ $\in $ $\Lambda \}$, hence $\hat{\sigma}_{n}\mathfrak{\hat{Z}}%
_{\psi ,n}^{\ast }(\pi ,\lambda )$ is the bootstrap\ draw from $\mathfrak{Z}%
_{\psi }(\pi ,\lambda )$.

\subsubsection*{Step 4: $\protect\tau _{\protect\beta }(\cdot ),$ $\mathfrak{%
\hat{T}}_{\protect\psi }(\cdot ),$ $\bar{v}^{2}(\cdot )$, $\mathcal{T}_{,%
\protect\psi }(\cdot )$}

We now have all the required components for computing the following key
quantities (recall $\mathcal{S}_{\beta }$ $\equiv $ $[I_{k_{\beta
}}:0_{k_{x}\times k_{x}}]$):%
\begin{eqnarray}
&&\hat{\tau}_{\beta ,n}^{\ast }(\pi _{0},b)\equiv -\mathcal{S}_{\beta }%
\widehat{\mathcal{H}}_{\psi ,n}^{-1}(\hat{\pi}_{n}^{\ast }(\pi
_{0},b))\left\{ \hat{\sigma}_{n}\widehat{\mathcal{G}}_{\psi ,n}^{\ast }(\hat{%
\pi}_{n}^{\ast }(\pi _{0},b))+\mathcal{\hat{D}}_{\psi ,n}(\hat{\pi}%
_{n}^{\ast }(\pi _{0},b),\pi _{0})\times b\right\}  \label{tao_boot} \\
&&\hat{\omega}_{n}^{\ast }(\pi _{0},b)\equiv \frac{\hat{\tau}_{\beta
,n}^{\ast }(\pi _{0},b)}{\left\Vert \hat{\tau}_{\beta ,n}^{\ast }(\pi
_{0},b)\right\Vert }  \notag \\
&&  \notag \\
&&\mathfrak{\hat{T}}_{\psi ,n}^{\ast }(\pi ,\lambda ,\pi _{0},b)  \notag \\
&&\text{ \ \ \ \ \ \ \ \ \ \ \ }\equiv \hat{\sigma}_{n}\mathfrak{\hat{Z}}%
_{\psi ,n}^{\ast }(\pi ,\lambda )+\mathfrak{\hat{b}}_{\psi ,n}(\hat{\psi}%
_{n},\pi ,\lambda )^{\prime }\left( \widehat{\mathcal{H}}_{\psi ,n}^{-1}(\pi
)\mathcal{\hat{D}}_{\psi ,n}(\pi ,\pi _{0})\times b+\left[ b,0_{k_{\beta
}}^{\prime }\right] ^{\prime }\right)  \notag \\
&&\text{ \ \ \ \ \ \ \ \ \ \ \ \ \ \ \ \ }+\mathfrak{\hat{b}}_{\psi ,n}(\hat{%
\psi}_{n},\pi ,\lambda )^{\prime }\widehat{\mathcal{H}}_{\psi ,n}^{-1}(\pi )%
\frac{1}{n}\sum_{t=1}^{n}d_{\psi ,t}(\pi )\left\{ g(x_{t},\pi
_{0})-g(x_{t},\pi )\right\} ^{\prime }b  \notag \\
&&\text{ \ \ \ \ \ \ \ \ \ \ \ \ \ \ \ \ }+\frac{1}{n}\sum_{t=1}^{n}\left\{
F\left( \lambda ^{\prime }\mathcal{W}(x_{t})\right) -\mathfrak{\hat{b}}%
_{\psi ,n}(\pi ,\lambda )^{\prime }\widehat{\mathcal{H}}_{\psi ,n}^{-1}(\pi
)d_{\psi ,t}(\pi )\right\} \left\{ g(x_{t},\pi _{0})-g(x_{t},\pi )\right\}
^{\prime }b  \notag \\
&&  \notag \\
&&\hat{v}_{n}^{2}(\omega ,\pi ,\lambda )\equiv \frac{1}{n}%
\sum_{t=1}^{n}\epsilon _{t}^{2}(\hat{\psi}_{n},\pi )\left\{ F\left( \lambda
^{\prime }\mathcal{W}(x_{t})\right) -\mathfrak{\hat{b}}_{\theta ,n}(\omega
,\pi ,\lambda )^{\prime }\widehat{\mathcal{H}}_{n}^{-1}(\omega ,\pi
)d_{\theta ,t}(\omega ,\pi )\right\} ^{2}  \notag \\
&&\widehat{\bar{v}}_{n}^{2}(\pi ,\lambda ,b)\equiv \hat{v}_{n}^{2}(\hat{%
\omega}_{n}^{\ast }(\pi _{0},b),\pi ,\lambda ).  \notag
\end{eqnarray}%
The bootstrap draw from $\mathcal{T}_{\psi }(\pi ^{\ast }(b),\lambda ,b)$ is
$\mathcal{\hat{T}}_{\psi ,n}^{\ast }(\lambda ,h)$ $=$ $\mathcal{\hat{T}}%
_{\psi ,n}^{\ast }(\lambda ,\pi _{0},b)$ $\equiv $ $\mathcal{\hat{T}}_{\psi
,n}^{\ast }(\hat{\pi}_{n}^{\ast }(\pi _{0},b),$ $\lambda ,\pi _{0},b)$ where%
\begin{equation}
\mathcal{\hat{T}}_{\psi ,n}^{\ast }(\pi ,\lambda ,\pi _{0},b)\equiv \left(
\frac{\mathfrak{\hat{T}}_{\psi ,n}^{\ast }(\pi ,\lambda ,\pi _{0},b)}{%
\widehat{\bar{v}}_{n}(\pi ,\lambda ,b)}\right) ^{2}.  \label{Tn_boot}
\end{equation}%
Notice $h$ $=$ $(\pi _{0},b)$ are nuisance parameters that cannot be
consistently estimated under weak identification $\sqrt{n}||\beta _{n}||$ $%
\rightarrow $ $[0,\infty )$.

\subsubsection*{Step 5}

Repeat Steps 1-4 $\mathcal{M}$ times resulting in a sequence of independent
draws $\{\mathcal{\hat{T}}_{\psi ,n,j}^{\ast }(\lambda ,h)\}_{j=1}^{\mathcal{%
M}}$. The p-value approximation is simply:
\begin{equation*}
\hat{p}_{n,\mathcal{M}}^{\ast }(\lambda ,h)\equiv \frac{1}{\mathcal{M}}%
\sum_{j=1}^{\mathcal{M}}I\left( \mathcal{\hat{T}}_{\psi ,n,j}^{\ast
}(\lambda ,h)>\mathcal{T}_{n}(\lambda )\right) .
\end{equation*}

Let $\hat{p}_{n,\mathcal{M}}^{(\cdot )}(\lambda )$\ be the LF or ICS-1
p-value computed with $\hat{p}_{n,\mathcal{M}}^{\ast }(\lambda ,h)$, and the
corresponding asymptotic size $AsySz^{\ast }(\lambda )$ $\equiv $ $%
\limsup_{n\rightarrow \infty }\sup_{\gamma \in \Gamma ^{\ast }}P_{\gamma }(%
\hat{p}_{n,\mathcal{M}}^{(\cdot )}(\lambda )$ $<$ $\alpha |H_{0})$ and $%
AsySz^{\ast }$\ $\equiv $ $\sup_{\lambda \in \Lambda }AsySz^{\ast }(\lambda
) $. $\hat{p}_{n,\mathcal{M}}^{\ast }(\lambda ,h)$ is consistent for the
asymptotic p-value under weak identification $p_{n}(\lambda ,h)$ $\equiv $ $%
1 $ $-$ $\mathcal{F}_{\lambda ,h}(\mathcal{T}_{n}(\lambda ))$, and the
resulting test achieves the correct uniform asymptotic level, $AsySz^{\ast }$
$\leq $ $\alpha $.

In order to demonstrate $AsySz^{\ast }$ $\leq $ $\alpha $ we need to verify
uniform convergence $\sup_{\lambda \in \Lambda }|\hat{p}_{n,\mathcal{M}%
_{n}}^{\ast }(\lambda ,h)$ $-$ $p_{n}(\lambda ,h)|$ $\overset{p}{\rightarrow
}$ $0$. Due to the nonsmooth structure of $\hat{p}_{n,\mathcal{M}_{n}}^{\ast
}(\lambda ,h)$ and how it enters $AsySz^{\ast }$, we need additional
structure on key processes. We exploit properties of the Vapnik-\v{C}%
ervonenkis\textit{\ subgraph }class of functions, denoted $\mathcal{V}(%
\mathcal{C})$. The $\mathcal{V}(\mathcal{C})$\ class is large: it contains
indicator, monotonic and continuous functions; and $\mathcal{V}(\mathcal{C})$%
\ mappings of $\mathcal{V}(\mathcal{C})$\ functions are in $\mathcal{V}(%
\mathcal{C})$, including linear combinations, minima, maxima, products and
indicator transforms. See, e.g., \citet[Chap. 2.6]{VaartWellner1996} for a
compendium of $\mathcal{V}(\mathcal{C})$ properties.\footnote{%
We exploit the facts that an indicator function of a $\mathcal{V}(\mathcal{C}%
)$\ index function is in $\mathcal{V}(\mathcal{C})$, and a continuous
function evaluated at a $\mathcal{V}(\mathcal{C})$\ function is in $\mathcal{%
V}(\mathcal{C})$.} See \cite{VC1971}, \citet[Section 7]{Dudley1978} and %
\citet[Section 2]{VaartWellner1996}, and see \citet[Chap.
II.4]{Pollard1984} for the closely related \textit{polynomial discrimination}
class.

Write $F_{n,\lambda }(c)$ $\equiv $ $P(\mathcal{T}_{n}(\lambda )$ $\leq $ $%
c) $ and $F_{n,\lambda ,h}^{\ast }(c)$ $\equiv $ $P(\mathcal{\hat{T}}_{\psi
,n,1}^{\ast }(\lambda ,h)$ $\leq $ $c|\mathfrak{W}_{n})$ where $\mathfrak{W}%
_{n}$ $\equiv $ $\{(y_{t},x_{t})\}_{t=1}^{n}$.

\begin{assumption}
\label{assum_pvs_uniform}The test weight $\{F(w)$ $:$ $w$ $\in $ $\mathbb{R}%
\}$ and distribution functions $\{F_{n,\lambda }(c)$ $:$ $\lambda $ $\in $ $%
\Lambda ,c\in \lbrack 0,\infty )\}$ and $\{F_{n,\lambda ,h}^{\ast }(c)$ $:$ $%
\lambda $ $\in $ $\Lambda ,c$ $\in $ $[0,\infty )\}$ belong to the $\mathcal{%
V}(\mathcal{C})$ class.
\end{assumption}

\begin{remark}
\normalfont The popularly used logistic and exponential weight functions $%
F(\cdot )$ are in $\mathcal{V}(\mathcal{C})$ because they are continuous.
Under Assumption \ref{assum:pvs} we impose continuity on the distribution
function $F_{n,\lambda }(\cdot )$, but we need more structure here to handle
uniform asymptotics over $\lambda $ for the bootstrapped p-value. The index
functions $F_{n,\lambda }(c)$ and $F_{n,\lambda ,h}^{\ast }(c)$ need to be
in $\mathcal{V}(\mathcal{C})$\ both in terms of the argument $c$ \emph{and}
the index $\lambda $ since they are evaluated at the test statistics $\{%
\mathcal{\hat{T}}_{\psi ,n,j}^{\ast }(\lambda ,h),\mathcal{T}_{n}(\lambda
)\} $.
\end{remark}

\begin{theorem}
\label{th:boot_pv}Let $\mathcal{M}$ $=$ $\mathcal{M}_{n}$ $\rightarrow $ $%
\infty $ as $n$ $\rightarrow $ $\infty $, and let Assumptions \ref{assum:dgp}%
, \ref{assum:pi}, \ref{assum:nondegen_v2_ae_all} and \ref{assum:pvs}\ hold. $%
a.$ $|\hat{p}_{n,\mathcal{M}}^{\ast }(\lambda ,h)$ $-$ $p_{n}(\lambda ,h)|$ $%
\overset{p}{\rightarrow }$ $0$. $b.$ If additionally Assumption \ref%
{assum_pvs_uniform} holds then $\sup_{\lambda \in \Lambda }|\hat{p}_{n,%
\mathcal{M}}^{\ast }(\lambda ,h)$ $-$ $p_{n}(\lambda ,h)|$ $\overset{p}{%
\rightarrow }$ $0$ and $AsySz^{\ast }$ $\leq $ $\alpha $.
\end{theorem}

Finally, we consider a theory for the PVOT test. Define the LF or ICS-1 PVOT
$\mathcal{\hat{P}}_{n,\mathcal{M}}(\alpha )$ $\equiv $ $\int_{\Lambda }I(%
\hat{p}_{n,\mathcal{M}}^{(\cdot )}(\lambda )$ $<$ $\alpha )d\lambda $. The
test rejects $H_{0}$ when $\mathcal{\hat{P}}_{n,\mathcal{M}}^{(\cdot
)}(\alpha )$ $>$ $\alpha $. The (non-uniform) asymptotic level of the test
is $\alpha $, and the test is consistent.

\begin{theorem}
\label{th:pvot}Let $\mathcal{M}$ $=$ $\mathcal{M}_{n}$ $\rightarrow $ $%
\infty $ as $n$ $\rightarrow $ $\infty $, and let Assumptions \ref{assum:dgp}%
, \ref{assum:pi}, \ref{assum:nondegen_v2_ae_all} and \ref{assum:pvs} hold.
Under $H_{0}$, $\lim_{n\rightarrow \infty }P(\mathcal{\hat{P}}_{n,\mathcal{M}%
}(\alpha )$ $>$ $\alpha )$ $\leq $ $\alpha $. Conversely, $P(\mathcal{\hat{P}%
}_{n,\mathcal{M}}(\alpha )$ $>$ $\alpha )$ $\rightarrow $ $1$ under $H_{1}$ $%
:$ $\sup_{\theta \in \Theta }P(E[y_{t}|x_{t}]$ $=$ $f(\theta ,x_{t}))$ $<$ $%
1 $.
\end{theorem}

The (uniform) asymptotic size of the PVOT test is $AsySz(pvot)$ $=$ $%
\limsup_{n\rightarrow \infty }\sup_{\gamma \in \Gamma ^{\ast }}P_{\gamma }(%
\mathcal{\hat{P}}_{n,\mathcal{M}}(\alpha )$ $>$ $\alpha |H_{0})$. Under the
additional structure of Assumption \ref{assum_pvs_uniform}, $AsySz(pvot)$ $%
\leq $ $\alpha $. Hence the PVOT test controls for size uniformly.

\begin{theorem}
\label{th:pvot_size}Let $\mathcal{M}$ $=$ $\mathcal{M}_{n}$ $\rightarrow $ $%
\infty $ as $n$ $\rightarrow $ $\infty $, and let Assumptions \ref{assum:dgp}%
, \ref{assum:pi}, \ref{assum:nondegen_v2_ae_all}, \ref{assum:pvs} and \ref%
{assum_pvs_uniform} hold. Then $AsySz(pvot)$ $\leq $ $\alpha $.
\end{theorem}

\section{MONTE CARLO STUDY\label{sec:sim}}

We now perform a simulation study in order to assess how well the proposed
bootstrap method works.

\subsection{SET UP}

Throughout $\epsilon _{t}$ is iid $N(0,1)$ distributed, $10,000$ samples are
generated, and sample sizes are $n$ $\in $ $\{100,250,500\}$. The wild
bootstrap used for robust p-value computation, and for the supremum and
average tests, uses $500$ bootstrap samples to reduce computation time.

The data generating process is%
\begin{equation*}
y_{t}=\zeta _{0}y_{t-1}+\beta _{n}y_{t-1}\frac{1}{1+\exp \left\{ -10\left(
y_{t-1}-\pi _{0}\right) \right\} }+\varpi _{0}\frac{1}{1+y_{t-1}^{2}}%
+\epsilon _{t}.
\end{equation*}%
If $\varpi _{0}$ $=$ $0$ then $y_{t}$ is a Logistic STAR model and the null
hypothesis is true. We use a fixed value $10$ for the speed of transition to
reduce computation complexity \citep[see also][]{AndrewsCheng2013}.\footnote{%
In empirical work when the transition function is, e.g., $(1$ $+$ $\exp
\{-\pi _{0,1}(y_{t-d}$ $-$ $\pi _{0,2})\})^{-1}$, often the transition speed
$\pi _{0,1}$ or threshold $\pi _{0,2}$ are fixed to ease computation, e.g.
\citet[p.
592]{LundberghTerasvirta2006} and \citet[Definition 1]{Gonzalez-Rivera1998}.}

We use $\zeta _{0}$ $=$ $.6$, $\pi _{0}$ $=$ $0$ and $\varpi _{0}$ $\in $ $%
\{0,.03,.3\}$. The latter allows us to inspect power against weak and strong
degrees of deviation from a STAR null hypothesis. The key parameter for
identification cases takes values $\beta _{n}$ $\in $ $\{.3,.3/\sqrt{n},0\}$%
, representing strong identification, weak identification with $\sqrt{n}%
\beta _{n}$ $=$ $b$ $=$ $.3$ and $\beta _{n}$ $\rightarrow $ $\beta _{0}$ $=$
$0$, and non-identification with $\beta _{n}$ $=$ $\beta _{0}$ $=$ $0$.
Other values for $(\zeta _{0},\beta _{n})$ lead to similar results.

Let $\iota $ $=$ $10^{-10}$. The true parameter spaces are $\mathcal{B}%
^{\ast }$ $=$ $[-1$ $+$ $2\iota ,1$ $-$ $2\iota ]$, $\mathcal{Z}^{\ast
}(\beta )$ $=$ $[-1-\beta $ $+$ $\iota <\zeta <1-\beta $ $-$ $\iota ]$, and $%
\Pi ^{\ast }$ $=$ $[-1,1]$. The estimation spaces are $\mathcal{B}$ $=$ $[-1$
$+$ $\iota ,1$ $-$ $\iota ]$, $\mathcal{Z}(\beta )$ $=$ $[-1-\beta <\zeta
<1-\beta ]$, and $\Pi $ $=$ $[-2,2]$. Thus $|\zeta $ $+$ $\beta |$ $<$ $1$
on $\Theta $ $\equiv $ $\mathcal{B}$ $\times $ $\mathcal{Z}(\beta )$ $\times
$ $\Pi $.

The estimated model is an LSTAR:%
\begin{equation*}
y_{t}=\zeta _{0}y_{t-1}+\beta _{0}y_{t-1}\frac{1}{1+\exp \left\{ -10\left(
y_{t-1}-\pi _{0}\right) \right\} }+\epsilon _{t}.
\end{equation*}%
We draw $100$ start values $\theta $ from the uniform distribution on $%
\Theta $ and estimate $\theta _{0}=[\zeta _{0},\beta _{0},\pi _{0}]^{\prime
} $ by least squares, resulting in $100$ estimates $\{\hat{\theta}%
_{n,i}\}_{i=1}^{100}$. The final estimate $\hat{\theta}_{n}$ \ minimizes the
least squares criterion over $\{\hat{\theta}_{n,i}\}_{i=1}^{100}$.\footnote{%
An analytic gradient is used for optimization. The criterion tolerance for
ceasing iterations is $1e^{-8}$, and the maximum number of allowed
iterations is $20,000$.} The conditional moment weight is logistic $F(u)$ $=$
$1/(1$ $+$ $\exp \{u\})$, and $F(\lambda ^{\prime }\mathcal{W}(x_{t}))$ uses
the bounded one-to-one transform $\mathcal{W}(x)$ $=$ atan$(x)$ as in
\citet[p. 1445,
1453]{Bierens1990}. The parameter space is $\Lambda $ $=$ $[1,5]$. We use a
discretization $\Lambda _{n}$ with endpoints $\{1,5\}$, and equal increments
with $n$ elements (e.g $\Lambda _{100}$ $=$ $\{1,$ $1.04,$ $1.08,...,$ $5$).

Eleven tests are performed. The first five are not robust to weak
identification: $(i)$ uniformly randomly chosen $\lambda ^{\ast }$ from $%
\Lambda _{n}$, compute $\mathcal{T}_{n}(\lambda ^{\ast })$ and use $\chi
^{2}(1)$ for p-value computation; $(ii)$ $\sup_{\lambda \in \Lambda
_{n}}p_{n}(\lambda )$; $(iii)$ $\sup_{\lambda \in \Lambda _{n}}\mathcal{T}%
_{n}(\lambda )$ and $(iv)$ $\int_{\Lambda _{n}}\mathcal{T}_{n}(\lambda )\mu
(d\lambda )$\ where $\mu $ is the uniform measure on $\Lambda _{n}$, and
p-values are computed by wild bootstrap; and $(v)$ the PVOT test using $%
\Lambda _{n}$, and a p-value computed from the $\chi ^{2}(1)$ distribution
[PVOT-$\chi ^{2}$].

The final six tests are robust. We compute $\mathcal{T}_{n}(\lambda ^{\ast
}) $ using ($vi$) the plug-in LF and ($vii$) plug-in ICS-1 p-values [$%
\mathcal{T}_{n}(\lambda ^{\ast })$-LF, $\mathcal{T}_{n}(\lambda ^{\ast })$%
-ICS]; $\sup_{\lambda \in \Lambda _{n}}p_{n}(\lambda )$ using ($viii$) the
plug-in LF and ($ix$) plug-in ICS-1 p-values [$\sup p_{n}$-LF, $\sup p_{n}$%
-ICS]; and PVOT using ($x$)\ the plug-in LF and ($xi$) plug-in ICS-1
p-values [PVOT-LF, PVOT-ICS]. The bootstrap procedure in Section \ref%
{sec:CM_pv_comp} is used to approximate the p-value under weak
identification $p_{n}(\lambda ,h)$ with $\hat{p}_{n,\mathcal{M}}^{\ast
}(\lambda ,h)$. Then $\hat{p}_{n,\mathcal{M}}^{\ast }(\lambda ,h)$ is used
to compute the plug-in LF and plug-in ICS-1 p-values from Section \ref%
{sec:CM_pv_construct}. Using ICS-2 would naturally lead to an improved
p-value, but the computational cost is too great at this time.

Theorems \ref{th:boot_pv}-\ref{th:pvot_size} provide the theory
demonstrating robustness and correct asymptotic size or level for tests ($vi$%
)-($xi$). Uniform bootstrap p-value convergence, and uniform size control,
require Assumption \ref{assum_pvs_uniform}: logistic $F(\cdot )$ is in the $%
\mathcal{V}(\mathcal{C})$ class. We need to assume $\{F_{n,\lambda }(c)$ $:$
$\lambda $ $\in $ $\Lambda ,c\in \lbrack 0,\infty )\}$ and $\{F_{n,\lambda
,h}^{\ast }(c)$ $:$ $\lambda $ $\in $ $\Lambda ,c$ $\in $ $[0,\infty )\}$
belong to the $\mathcal{V}(\mathcal{C})$ class due to their nonlinear
complexity.

The computation of LF and ICS p-values using $\hat{p}_{n,\mathcal{M}}^{\ast
}(\lambda ,h)$ requires a grid of nuisance parameters $h$ $=$ $(\pi ,b)$. We
use $\pi \in \{-2,-1.5,...,1.5,2\}$ and $b$ $\in $ $%
\{-.5,-.3,-.2,-.1,0,.1,.2,.3,.5\}$. Finer grids lead to significant
increases in computation time. The ICS-1 p-value require the threshold $%
\kappa _{n}$: we use $\kappa _{n}$ $=$ $a\ln (\ln (n))$ with $a$ $=$ $1$.
Values of $a$\ close to $1$ lead to similar results, while under rejection
of the null is more prominent as $a$ increases. Larger rates of increase for
$\kappa _{n}$ like $c(\ln (n))^{\delta }$ for some $c,\delta $ $>$ $0$
generally result in the ICS-1 p-value being nearly equal to the LF p-value,
at least within our chosen design. Indeed, under $\kappa _{n}$ $=$ $(\ln
(n))^{1/2}$ there is little difference between LF and ICS-1 values. Finally,
the selection matrix $\mathcal{S}_{\beta }$ $\equiv $ $[I_{k_{\beta
}}:0_{k_{x}\times k_{x}}]$ reduces to $\mathcal{S}_{\beta }$ $\equiv $ $%
[1,0] $ since $k_{\beta }$ $=$ $k_{x}$ $=$ $1$.

\subsection{RESULTS}

Rejection frequencies are given in Tables \ref{tbl:starn100}-\ref%
{tbl:starn500}.

\subsubsection{Strong Identification}

Consider the strong identification case $\beta _{n}$ $=$ $.3$. Under the
null, $\mathcal{T}_{n}(\lambda ^{\ast })$ and PVOT-$\chi ^{2}$ exhibit
rejection rates close to the nominal levels. The supremum test is over-sized
and exhibits the largest size distortion, while the average test is slightly
over-sized.

The LF p-value leads to under-rejection for $\mathcal{T}_{n}(\lambda ^{\ast
})$-LF and PVOT-LF at each $n$. The ICS-1 p-value with $\kappa _{n}$ $=$ $%
\ln (\ln (n))$ corrects for the size distortion in nearly every case. The
exceptions are when $n$ $=$ $100$ at the $10\%$ level for both $\mathcal{T}%
_{n}(\lambda ^{\ast })$-ICS and PVOT-ICS, and when $n$ $=$ $250$ at the $%
10\% $ level for $\mathcal{T}_{n}(\lambda ^{\ast })$-ICS. In these cases
empirical size is roughly $7\%$. The PVOT-ICS therefore yields sharp size in
nearly every case. The slight advantage of ICS-1 over LF is not surprising
since the LF p-value is generally larger.

The supremum and PVOT-$\chi ^{2}$ tests have the largest size corrected
power under the weak alternative (raw power is displayed). The robust PVOT
tests perform better than the robust tests based on $\mathcal{T}_{n}(\lambda
^{\ast })$; PVOT-ICS performs better than PVOT-LF at $n$ $\in $ $\{100,250\}$%
, but the two are essentially identical at $n$ $=$ $500$; and PVOT-ICS
approaches the size corrected power of supremum and PVOT-$\chi ^{2}$ tests
as $n$ increases.

Under the strong alternative, supremum, average and PVOT-ICS tests are
comparable, although the average test is weaker at $n$ $=$ $100$. Both $%
\mathcal{T}_{n}(\lambda ^{\ast })$-LF and $\mathcal{T}_{n}(\lambda ^{\ast })$%
-ICS yield lower power than PVOT-LF and PVOT-ICS. Generally the LF p-value
results in lower rejection rates than the ICS-1 p-value since the LF p-value
is larger (hence rejection is less likely).

\subsubsection{Weak Identification}

Now consider weak and non-identification $\beta _{n}$ $=$ $3/n^{1/2}$ and $%
\beta _{n}$ $=$ $0$. Each non-robust test is strongly over-sized, up to an
order of $3$ to $5$, depending on the significance level. As an example, at $%
n$ $=$ $100$ under non-identification $\beta _{n}$ $=$ $0$ the rejection
rates for PVOT--$\chi ^{2}$ are $\{.050,.140,.205\}$ at nominal sizes $%
(1\%,5\%,10\%)$, compared to $\{.015,.065,.124\}$ under strong
identification. At $n$ $=$ $500$ the rates are $\{.061,.148,.208\}$ and $%
\{.014,.055,.115\}$ under non- and strong identification. Thus $\mathcal{T}%
_{n}(\lambda )$ is strongly positively skewed relative to the $\chi ^{2}(1)$
distribution. The remaining tests are qualitatively similar. For example,
the average test based on $\int_{\Lambda _{n}}\mathcal{T}_{n}(\lambda )\mu
(d\lambda )$ generates rejection frequencies $\{.057,.146,.219\}$ when $n$ $%
= $ $100$ under weak identification. These drop to $\{.029,.125,.176\}$ when
$n $ $=$ $500$.

LF and ICS-1 p-values lead to correct size for both $\mathcal{T}_{n}(\lambda
^{\ast })$-ICS and PVOT-ICS tests. Under the alternative, however, the ICS-1
p-value leads to a power gain that reaches close to 25\%. Thus, the
construction of the LF p-value works well under the null, but not
surprisingly weakens empirical power. The maximum gain for ICS-1 occurs
under the weak alternative, with weak identification and $n$ $=$ $100$: see
Table \ref{tbl:starn100} (middle panel, fourth column). The typical gain is
5\%-10\% depending on the alternative and sample size.

The PVOT-ICS test generally has the greatest size corrected power, in
particular under $(i)$ the strong alternative at $n$ $=$ $100$, and the $5\%$
and $10\%$ levels; $(ii)$ the strong alternative at $n$ $\geq $ $250$; and $%
(iii)$ the weak alternative at the $5\%$ and $10\%$ levels, when $n$ $\geq $
$250$. Under those alternatives and sample sizes the supremum and PVOT-ICS
tests are comparable at the $1\%$ level. In the remaining cases the supremum
test has the largest size corrected power with a margin of about $5\%$-$10\%$%
.

When empirical size and size corrected power are considered jointly,
PVOT-ICS is the most promising test in this study across cases, and in terms
of robustness against weak and non-identification. Average and supremum
tests with wild bootstrapped p-values exhibit large size distortions under
weak and non-identification, while PVOT-ICS controls for size, and yields
competitive or dominant size-corrected power.

\section{CONCLUSION\label{sec:conclusion}}

We offer a new bootstrap procedure that is robust to any degree of
(non)identification in nonlinear regression models. The procedure targets
case specific degrees of identification, avoiding the breakdown of uniform
bootstrap asymptotic validity over the parameter space. We focus on a
conditional moment test of functional form, but the method extends to a wide
variety of tests. An occupation time smoothed bootstrapped p-value leads to
a test that achieves uniform size control, and is consistent.

The procedure works well in a simulation experiment, in particular when the
proposed bootstrapped p-value is imbedded in the p-value occupation time.
Future work may include expanding the proposed bootstrap to other tests
where identification is a potential problem, including tests of white noise
for model residuals, structural break tests, and so on.

\section{SUPPLEMENTARY MATERIAL}
Hill, J. B. (2020): Supplement to ``Weak-Identification Robust Wild Bootstrap applied to a Consistent Model Specification Test", Econometric Theory Supplementary Material.

\bigskip\bigskip
\bibliographystyle{econometrica}
\bibliography{refs_weak_ident_boot}

\setcounter{equation}{0} \renewcommand{\theequation}{{\thesection}.%
\arabic{equation}} \appendix
\setstretch{1.2}

\section{APPENDIX\label{app:A}}

\subsection{ASSUMPTION \protect\ref{assum:dgp}.d,e,f\label{append:Assum1def}}

Assumptions \ref{assum:dgp}.d,e,f contain technical restrictions on long-run
variances, and parameter space details that are useful when any degree of
identification is allowed.

In order to conserve space below, we use the following notation. Let \emph{%
int}$(\mathcal{M})$ denote the interior of set $\mathcal{M}$. Write
compactly $\inf_{\alpha ,r,\theta }$ $=$ $\inf_{\alpha ^{\prime }\alpha
=1,r^{\prime }r=1,\{\theta _{i}\}_{i=1}^{m}\in \Theta ^{m}}$, $\inf_{\alpha
,\pi ,\lambda }$ $=$ $\inf_{\alpha ^{\prime }\alpha =1,\{\pi _{i},\lambda
_{i}\}_{i=1}^{m}\in (\Pi \times \Lambda )^{m}}$, $\inf_{r,\omega ,\pi }$ $=$
$\inf_{r^{\prime }r=1,\omega \in \mathbb{R}^{k_{\beta }}:||\omega ||=1,\pi
\in \Pi }$, and so on. $m$ $\in $ $\mathbb{N}$ and $a$ $\in $ $\mathbb{R}%
^{m} $\ are arbitrary; $\Theta ^{m}$ $\equiv $ $\Theta \times \cdots \times
\Theta $ $\subset $ $\mathbb{R}^{m}$. $r$ $=$ $[r_{1},r_{2}^{\prime
}]^{\prime }$, with $r_{1}$ $\in $ $\mathbb{R}$, is an arbitrary vector
whose dimension is implicitly defined. We need the following definitions:

\begin{eqnarray*}
&&\mathbb{E}_{\psi ,n}(\pi ;a,r)\equiv \frac{1}{\sqrt{n}}\sum_{t=1}^{n}%
\epsilon _{t}\sum_{i=1}^{m}\alpha _{i}r^{\prime }d_{\psi ,t}(\pi _{i})\text{
and }\mathbb{E}_{\theta ,n}(\omega ,\pi ;a,r)\equiv \frac{1}{\sqrt{n}}%
\sum_{t=1}^{n}\epsilon _{t}\sum_{i=1}^{m}\alpha _{i}r^{\prime }d_{\theta
,t}(\omega _{i},\pi _{i}) \\
&&\mathfrak{EG}_{\psi ,n}(\lambda ;a,r)\equiv r_{1}\frac{1}{\sqrt{n}}%
\sum_{t=1}^{n}\sum_{i=1}^{m}\alpha _{i}\left\{ \epsilon _{t}(\psi _{n},\pi
_{i})\mathcal{K}_{\psi ,t}(\pi _{i},\lambda _{i})-E\left[ \epsilon _{t}(\psi
_{n},\pi _{i})\mathcal{K}_{\psi ,t}(\pi _{i},\lambda _{i})\right] \right\} \\
&&\text{ \ \ \ \ \ \ \ \ \ \ \ \ \ \ \ \ \ \ \ \ \ }+r_{2}^{\prime
}\sum_{i=1}^{m}\alpha _{i}\mathcal{G}_{\psi ,n}(\psi _{n},\pi _{i}).
\end{eqnarray*}%
\textbf{Assumption \ref{assum:dgp} }(data generating process, test
weight)\medskip \newline
d. \emph{Long-Run Variances}:

$(i)$ Under $\mathcal{C}(i,b)$ with $||b||$ $<$ $\infty $\ let $\lim
\inf_{n\rightarrow \infty }E[\inf_{\alpha ,r,\theta }\mathfrak{(}r^{\prime
}\sum_{i=1}^{m}\alpha _{i}\mathcal{G}_{\psi ,n}(\theta _{i}))^{2}]$ $>$ $0$
and \linebreak $\lim \sup_{n\rightarrow \infty }E[\sup_{\alpha ,r,\theta }%
\mathfrak{(}r^{\prime }\sum_{i=1}^{m}\alpha _{i}\mathcal{G}_{\psi ,n}(\theta
_{i}))^{2}]$ $<$ $\infty $.

$(ii)$ Under $\mathcal{C}(ii,\omega _{0})$ let $\lim \inf_{n\rightarrow
\infty }E[\inf_{\alpha ,r,\theta }\mathfrak{(}r^{\prime
}\sum_{i=1}^{m}\alpha _{i}\mathcal{G}_{\theta ,n}(\theta _{i}))^{2}]$ $>$ $0$
and \linebreak $\lim \sup_{n\rightarrow \infty }E[\sup_{\alpha ,r,\theta }%
\mathfrak{(}r^{\prime }\sum_{i=1}^{m}\alpha _{i}\mathcal{G}_{\theta
,n}(\theta _{i}))^{2}]$ $<$ $\infty $.

$(iii)$ $E[\inf_{r,\omega ,\pi }(r^{\prime }d_{\theta ,t}(\omega ,\pi
))^{2}] $ $>$ $0$ and $E[\sup_{r,\omega ,\pi }(r^{\prime }d_{\theta
,t}(\omega ,\pi ))^{2}]$ $<$ $\infty $; $E\left[ \inf_{r,\pi }(r^{\prime
}d_{\psi ,t}(\pi ))^{2}\right] $ $>$ $0$ and $E\left[ \sup_{r,\pi
}(r^{\prime }d_{\psi ,t}(\pi ))^{2}\right] $ $<$ $\infty $.

$(iv)$ $\lim \inf_{n\rightarrow \infty }\inf_{a,r,\pi }E[\mathbb{E}_{\psi
,n}(\pi ;a,r)^{2}]$ $>0$ and $\lim \sup_{n\rightarrow \infty }\sup_{a,r,\pi
}E[\mathbb{E}_{\psi ,n}(\pi ;a,r)^{2}]$ $<$ $\infty $; and $\lim
\inf_{n\rightarrow \infty }\inf_{a,r,\omega ,\pi }E[\mathbb{E}_{\theta
,n}(\omega ,\pi ;a,r)^{2}]$ $>$ $0$ and $\lim \sup_{n\rightarrow \infty
}\sup_{a,r,\omega ,\pi }E[\mathbb{E}_{\theta ,n}(\omega ,\pi ;a,r)^{2}]$ $<$
$\infty $.

$(v)$ Under $\mathcal{C}(i,b)$ with $||b||$ $<$ $\infty $, $\lim
\inf_{n\rightarrow \infty }E[\sup_{\alpha ,r,\lambda }\mathfrak{EG}_{\psi
,n}(\lambda ;a,r)^{2}]$ $<$ $\infty $.

$(vi)$ Under $\mathcal{C}(ii,\omega _{0})$, $E[\sup_{\alpha ,r,\lambda }(1/%
\sqrt{n}\sum_{t=1}^{n}\epsilon _{t}\mathcal{K}_{\theta ,t}(\lambda
;a,m))^{2}]$ $<$ $\infty $\ for each $m$.\medskip \newline
e. \emph{True Parameter Space}:

$(i)$ $\Theta ^{\ast }$ $\equiv $ $\{(\beta ,\zeta ,\pi )$ $:$ $\beta $ $\in
$ $\mathcal{B}^{\ast },$ $\zeta $ $\in $ $\mathcal{Z}^{\ast }(\beta ),$ $\pi
$ $\in $ $\Pi ^{\ast }\}$ is\ compact.

$(ii)$ $0_{k_{\beta }}$ $\in $ \emph{int}$(\mathcal{B}^{\ast })$.

$(iii)$ For some set $\mathcal{Z}_{0}^{\ast }$ and some $\delta $ $>$ $0$, $%
\mathcal{Z}^{\ast }(\beta )$ $=$ $\mathcal{Z}_{0}^{\ast }$ $\forall ||\beta
||$ $<$ $\delta $.\medskip \newline
f. \emph{Optimization Parameter Space}:

$(i)$ $\Theta $ $\equiv $ $\{(\beta ,\zeta ,\pi )$ $:$ $\beta $ $\in $ $%
\mathcal{B},$ $\zeta $ $\in $ $\mathcal{Z}(\beta ),$ $\pi $ $\in $ $\Pi \}$
and $\Theta ^{\ast }$ $\subset $ $\emph{int}(\Theta )$.

$(ii)$ $(\Theta ,\mathcal{B},\Pi )$ are compact, and $\mathcal{Z}(\beta )$
is compact for each $\beta $. $(iii)$ For some set $\mathcal{Z}_{0}$ and
some $\delta $ $>$ $0$, $\mathcal{Z}(\beta )$ $=$ $\mathcal{Z}_{0}$ $\forall
||\beta ||$ $<$ $\delta $ and $\mathcal{Z}_{0}^{\ast }$ $\subset $ $\emph{int%
}(\mathcal{Z}_{0})$.

\begin{remark}
\normalfont(d.i,ii) are standard for non-degenerate finite dimensional
asymptotics for the least squares first order equations, under stationarity.
(d.iii,iv) likewise imply the components of the least squares asymptotic
variance are non-degenerate and positive definite. Each is standard for
ensuring non-degenerate asymptotics. (d.v) promotes a joint weak limit theory
for the test statistic numerator $1/\sqrt{n}\sum_{t=1}^{n}\epsilon _{t}(\psi
_{n},\pi )F(\lambda ^{\prime }\mathcal{W}(x_{t}))$ and $\hat{\pi}_{n}$ under
weak identification, which in turn leads to a limit theory for $1/\sqrt{n}%
\sum_{t=1}^{n}\epsilon _{t}(\hat{\theta}_{n})F(\lambda ^{\prime }\mathcal{W}%
(x_{t}))$. (d.vi) similarly covers $1/\sqrt{n}\sum_{t=1}^{n}\epsilon _{t}(%
\hat{\theta}_{n})F(\lambda ^{\prime }\mathcal{W}(x_{t}))$ under strong
identification. In the latter cases (d.v) and (d.vi) we cannot assume strict
positivity due to possible degeneracy: see the discussion leading to
Assumptions \ref{assum:nondegen_v2_ae} and \ref{assum:nondegen_v2_ae_all}.
\end{remark}

\begin{remark}
\normalfont(e.ii) ensures non-identification ($\beta $ $=$ $0$) and near
non-identification ($\beta $ close to $0$) points are in $\mathcal{B}^{\ast }
$. (e.iii) ensures the true space $\mathcal{Z}^{\ast }(\beta )$ is bounded
from the empty set for values of $\beta $ near the non-identification point $%
\beta $ $=$ $0$, and therefore allows for derivatives of certain moments
with respect to $\beta $ near $\beta $ $=$ $0$
\citep[cf.][comments following Assumptions B2(iii) and
STAR4(iv)]{AndrewsCheng2013}.\footnote{%
In the Logistic STAR model in the preceding footnote it suffices to assume
compact $\mathcal{B}^{\ast }$ $\subset $ $(-\infty ,\infty )$ $\times $ $(-1$
$+$ $\iota ,1$ $-$ $\iota )$ for some infinitesimal $\iota $ $>$ $0$, and $%
\mathcal{Z}^{\ast }(\beta )$ is a compact subset of $\{\zeta $ $\in $ $%
\mathbb{R}^{2}$ $:$ $\zeta _{1}$ $\in $ $(-\infty ,\infty )$, $-1$ $<$ $%
\zeta _{2}$ $+$ $\beta _{2}$ $<$ $1\}$. Now assume $\mathcal{Z}_{0}^{\ast }$
$=$ $\{\zeta $ $\in $ $\mathbb{R}^{2}$ $:$ $\zeta _{1}$ $\in $ $(-\infty
,\infty ),-1$ $+$ $\iota $ $<$ $\zeta _{2}$ $<$ $1$ $-$ $\iota \}$ and pick $%
\delta $ $=$ $\iota $. The same idea extends to $\mathcal{Z}_{0}$.} The
(f.i) property $\Theta ^{\ast }$ $\subset $ $\emph{int}(\Theta )$ implies
the true value does not lie on the boundary of the optimization space. This
is non-essential, but allows a focus on weak identification.
\end{remark}

\subsection{PROOFS OF MAIN RESULTS\label{app:proofs}}

We assume all random variables exist on a complete measure space such that
majorants and integrals over uncountable families of measurable functions
are measurable, and probabilities where applicable are outer probability
measures.\footnote{%
See Pollard's (\citeyear{Pollard1984}: Appendix C) \textit{permissibility}
criteria, and see Dudley's (\citeyear{Dudley1984}: p. 101) \textit{%
admissible Suslin} property.}

In order to conserve space, write processes on compact spaces variously as $%
\{f(a,b,c)$ $:$ $A,B,C\}$ $=$ $\{f(a,b,c)$ $:$ $a$ $\in $ $A,b$ $\in $ $B,c$
$\in $ $C\}$. $A_{n}(\lambda )$ $=$ $o_{p,\lambda }(1)$ implies $%
\sup_{\lambda \in \Lambda }||A_{n}(\lambda )||$ $\overset{p}{\rightarrow }$ $%
0$. All Gaussian processes below have a version that has \emph{almost surely}
continuous and uniformly bounded sample paths, hence we just say \textit{%
Gaussian process}.

The following proofs require supporting results presented in Appendix \ref%
{app:lemmas} and proven in the supplemental material
\citet[Appendix
B]{Supp_Mat_2020}. Recall $\omega (\beta )$ $\equiv $ $\{\beta /||\beta ||$
if $\beta $ $\neq $ $0$, $1_{k_{\beta }}/||1_{k_{\beta }}||$ if $\beta $ $=$
$0$.

Recall the augmented parameter $\theta ^{+}$ $\equiv $ $[||\beta ||,\omega
^{\prime },\zeta ^{\prime },\pi ^{\prime }]^{\prime }$ $\in $ $\Theta ^{+}$
where $\Theta ^{+}$ $\equiv $ $\{\theta ^{+}$ $\in $ $\mathbb{R}^{k_{\beta
}+k_{x}+k_{\pi }+1}$ $:$ $\theta ^{+}$ $=$ $[||\beta ||,\omega (\beta
),\zeta ,\pi ]^{\prime }$ $:$ $\beta $ $\in $ $\mathcal{B},$ $\zeta $ $\in $
$\mathcal{Z}(\beta ),$ $\pi $ $\in $ $\Pi \}$.\medskip \newline
\textbf{Proof of Theorem \ref{th:CM_weak}.}\medskip \newline
\textbf{Claim a.}\qquad Let drift case $\mathcal{C}(i,b)$ hold with $||b||$ $%
<$ $\infty $. Recall $\epsilon _{t}(\theta )$ $=$ $\epsilon _{t}(\psi ,\pi )$%
, and define $f(x_{t},\theta )$ $=$ $\zeta ^{\prime }x_{t}$ $+$ $\beta
^{\prime }g(x_{t},\pi )$.\medskip \newline
\textbf{Step 1.}\qquad We prove the following expansion:
\begin{eqnarray}
&&\frac{1}{\sqrt{n}}\sum_{t=1}^{n}\epsilon _{t}(\hat{\theta}_{n})F\left(
\lambda ^{\prime }\mathcal{W}(x_{t})\right)  \label{eF_expand} \\
&&\text{ \ \ \ }=\frac{1}{\sqrt{n}}\sum_{t=1}^{n}\left\{ \epsilon _{t}(\psi
_{n},\hat{\pi}_{n})\mathcal{K}_{\psi ,t}(\hat{\pi}_{n},\lambda )-E\left[
\epsilon _{t}(\psi _{n},\hat{\pi}_{n})\mathcal{K}_{\psi ,t}(\hat{\pi}%
_{n},\lambda )\right] \right\}  \notag \\
&&\text{ \ \ \ \ \ \ \ \ \ \ \ \ \ \ \ \ \ \ \ }+E\left[ \mathcal{K}_{\psi
,t}(\hat{\pi}_{n},\lambda )\left\{ g(x_{t},\pi _{0})-g(x_{t},\hat{\pi}%
_{n})\right\} ^{\prime }\right] b  \notag \\
&&\text{ \ \ \ \ \ \ \ \ \ \ \ \ \ \ \ \ \ \ \ }+\mathfrak{b}_{\psi }(\hat{%
\pi}_{n},\lambda )^{\prime }\mathcal{H}_{\psi }^{-1}(\hat{\pi}_{n})E\left[
d_{\psi ,t}(\hat{\pi}_{n})\left\{ g(x_{t},\pi _{0})-g(x_{t},\hat{\pi}%
_{n})\right\} ^{\prime }\right] b  \notag \\
\text{\ \ \ \ } &&\text{\ \ \ \ \ \ \ \ \ \ \ \ \ \ \ \ \ \ \ \ }+\mathfrak{b%
}_{\psi }(\hat{\pi}_{n},\lambda )^{\prime }\left\{ \mathcal{H}_{\psi }^{-1}(%
\hat{\pi}_{n})\mathcal{D}_{\psi }(\hat{\pi}_{n})b+\left[ b,0_{k_{\beta
}}^{\prime }\right] ^{\prime }\right\} +o_{p,\lambda }\left( 1\right)  \notag
\\
&&\text{ \ \ \ }=\mathfrak{Z}_{n}(\hat{\pi}_{n},\lambda )+\mathcal{R}(\hat{%
\pi}_{n},\lambda )+o_{p,\lambda }\left( 1\right) ,  \notag
\end{eqnarray}%
where $\mathcal{R}(\pi ,\lambda )$ is implicitly defined, and
\begin{equation*}
\mathfrak{Z}_{n}(\pi ,\lambda )\equiv \frac{1}{\sqrt{n}}\sum_{t=1}^{n}\left%
\{ \epsilon _{t}(\psi _{n},\pi )\mathcal{K}_{\psi ,t}(\pi ,\lambda )-E\left[
\epsilon _{t}(\psi _{n},\pi )\mathcal{K}_{\psi ,t}(\pi ,\lambda )\right]
\right\} .
\end{equation*}

Recall $\hat{\theta}_{n}$ $=$ $[\hat{\psi}_{n}(\hat{\pi}_{n})^{\prime },\hat{%
\pi}_{n}^{\prime }]^{\prime }$ and write $\hat{\psi}_{n}$ $=$ $\hat{\psi}%
_{n}(\hat{\pi}_{n})$. By the mean value theorem, there exists $\psi
_{n}^{\ast }$ $\in $ $\Psi $, $||\psi _{n}^{\ast }$ $-$ $\psi _{n}||$ $\leq $
$||\hat{\psi}_{n}$ $-$ $\psi _{n}||$, such that:%
\begin{eqnarray*}
\frac{1}{\sqrt{n}}\sum_{t=1}^{n}\epsilon _{t}(\hat{\theta}_{n})F\left(
\lambda ^{\prime }\mathcal{W}(x_{t})\right) &=&\frac{1}{\sqrt{n}}%
\sum_{t=1}^{n}\epsilon _{t}(\psi _{n},\hat{\pi}_{n})F\left( \lambda ^{\prime
}\mathcal{W}(x_{t})\right) \\
&&-\frac{1}{n}\sum_{t=1}^{n}F\left( \lambda ^{\prime }\mathcal{W}%
(x_{t})\right) \frac{\partial }{\partial \psi ^{\prime }}f(x_{t,},\left[
\psi _{n}^{\ast },\hat{\pi}_{n}\right] )\sqrt{n}\left( \hat{\psi}_{n}-\psi
_{n}\right) \\
&=&\frac{1}{\sqrt{n}}\sum_{t=1}^{n}\epsilon _{t}(\psi _{n},\hat{\pi}%
_{n})F\left( \lambda ^{\prime }\mathcal{W}(x_{t})\right) -\mathfrak{\hat{b}}%
_{\psi ,n}(\hat{\pi}_{n},\lambda )^{\prime }\sqrt{n}\left( \hat{\psi}%
_{n}-\psi _{n}\right) .
\end{eqnarray*}%
By Lemma \ref{lm:bn}, $\sup_{\pi \in \Pi ,\lambda \in \Lambda }||\mathfrak{%
\hat{b}}_{\psi ,n}(\pi ,\lambda )$ $-$ $\mathfrak{b}_{\psi }(\pi ,\lambda
)|| $ $\overset{p}{\rightarrow }$ $0$. The proof of Theorem \ref{th:ols}
verifies \citep[see (C.19) in][]{Supp_Mat_2020}:
\begin{equation}
\sup_{\pi \in \Pi }\left\Vert \sqrt{n}\left( \hat{\psi}_{n}(\pi )-\psi
_{n}\right) -\left( -\mathcal{H}_{\psi }^{-1}(\pi )\left\{ \mathcal{G}_{\psi
,n}(\psi _{0,n},\pi )+\mathcal{D}_{\psi }(\pi )b\right\} -\left[
b,0_{k_{\beta }}^{\prime }\right] ^{\prime }\right) \right\Vert \overset{p}{%
\rightarrow }0.  \label{psi_expand1}
\end{equation}%
Now apply (\ref{psi_expand1}) for $\sqrt{n}(\hat{\psi}_{n}$ $-$ $\psi _{n})$
and the Lemma \ref{lm:H_ulln_weak} uniform consistency of $\widehat{\mathcal{%
H}}_{\psi ,n}(\pi )$ to yield:%
\begin{eqnarray*}
&&\frac{1}{\sqrt{n}}\sum_{t=1}^{n}\epsilon _{t}(\hat{\theta}_{n})F\left(
\lambda ^{\prime }\mathcal{W}(x_{t})\right) \\
&&\text{ \ \ \ \ \ \ \ }=\frac{1}{\sqrt{n}}\sum_{t=1}^{n}\epsilon _{t}(\psi
_{n},\hat{\pi}_{n})F\left( \lambda ^{\prime }\mathcal{W}(x_{t})\right) \\
&&\text{ \ \ \ \ \ \ \ \ \ \ \ \ \ \ \ }-\mathfrak{b}_{\psi }(\hat{\pi}%
_{n},\lambda )^{\prime }\left\{ -\mathcal{H}_{\psi }^{-1}(\hat{\pi}%
_{n})\left\{ \mathcal{G}_{\psi ,n}(\psi _{0,n},\hat{\pi}_{n})+\mathcal{D}%
_{\psi }(\hat{\pi}_{n})b\right\} -\left[ b,0_{k_{\beta }}^{\prime }\right]
^{\prime }\right\} +o_{p,\lambda }\left( 1\right) \\
&&\text{ \ \ \ \ \ \ \ }=\frac{1}{\sqrt{n}}\sum_{t=1}^{n}\epsilon _{t}(\psi
_{n},\hat{\pi}_{n})F\left( \lambda ^{\prime }\mathcal{W}(x_{t})\right) +%
\mathfrak{b}_{\psi }(\hat{\pi}_{n},\lambda )^{\prime }\mathcal{H}_{\psi
}^{-1}(\hat{\pi}_{n})\mathcal{G}_{\psi ,n}(\psi _{0,n},\hat{\pi}_{n}) \\
&&\text{ \ \ \ \ \ \ \ \ \ \ \ \ \ \ \ }+\mathfrak{b}_{\psi }(\hat{\pi}%
_{n},\lambda )^{\prime }\left\{ \mathcal{H}_{\psi }^{-1}(\hat{\pi}_{n})%
\mathcal{D}_{\psi }(\hat{\pi}_{n})b+\left[ b,0_{k_{\beta }}^{\prime }\right]
^{\prime }\right\} +o_{p,\lambda }\left( 1\right) .
\end{eqnarray*}

Next, by the construction of $\mathcal{G}_{\psi ,n}(\theta )$ in (\ref{GG}):%
\begin{eqnarray*}
-\mathcal{G}_{\psi ,n}(\psi _{0,n},\hat{\pi}_{n}) &=&-\frac{1}{\sqrt{n}}%
\sum_{t=1}^{n}\left\{ \epsilon _{t}(\psi _{0,n},\hat{\pi}_{n})d_{\psi ,t}(%
\hat{\pi}_{n})-E\left[ \epsilon _{t}(\psi _{0,n},\hat{\pi}_{n})d_{\psi ,t}(%
\hat{\pi}_{n})\right] \right\} \\
&=&-\frac{1}{\sqrt{n}}\sum_{t=1}^{n}\left\{ \epsilon _{t}(\psi _{n},\hat{\pi}%
_{n})d_{\psi ,t}(\hat{\pi}_{n})-E\left[ \epsilon _{t}(\psi _{n},\hat{\pi}%
_{n})d_{\psi ,t}(\hat{\pi}_{n})\right] \right\} \\
&&+\sqrt{n}\frac{1}{n}\sum_{t=1}^{n}\left\{ d_{\psi ,t}(\hat{\pi}%
_{n})g(x_{t},\hat{\pi}_{n})^{\prime }-E\left[ d_{\psi ,t}(\hat{\pi}%
_{n})g(x_{t},\hat{\pi}_{n})^{\prime }\right] \right\} \beta _{n}.
\end{eqnarray*}%
Combine $\sqrt{n}||\beta _{n}||$ $\rightarrow $ $[0,\infty )$, Lemma \ref%
{lm:H_ulln_weak}, and Theorem \ref{th:ols} to yield: $\sqrt{n}%
n^{-1}\sum_{t=1}^{n}\{d_{\psi ,t}(\hat{\pi}_{n})g(x_{t},\hat{\pi}%
_{n})^{\prime }$ $-$ $E[d_{\psi ,t}(\hat{\pi}_{n})g(x_{t},\hat{\pi}%
_{n})^{\prime }]\}\beta _{n}$ $\overset{p}{\rightarrow }$ $0$. Therefore%
\begin{eqnarray}
&&\frac{1}{\sqrt{n}}\sum_{t=1}^{n}\epsilon _{t}(\hat{\theta}_{n})F\left(
\lambda ^{\prime }\mathcal{W}(x_{t})\right)  \label{sumEF} \\
&&\text{ \ \ \ }=\frac{1}{\sqrt{n}}\sum_{t=1}^{n}\epsilon _{t}(\psi _{n},%
\hat{\pi}_{n})F\left( \lambda ^{\prime }\mathcal{W}(x_{t})\right)  \notag \\
&&\text{ \ \ \ \ \ \ \ \ \ \ \ \ \ \ \ }-\mathfrak{b}_{\psi }(\hat{\pi}%
_{n},\lambda )^{\prime }\mathcal{H}_{\psi }^{-1}(\hat{\pi}_{n})\frac{1}{%
\sqrt{n}}\sum_{t=1}^{n}\left\{ \epsilon _{t}(\psi _{n},\hat{\pi}_{n})d_{\psi
,t}(\hat{\pi}_{n})-E\left[ \epsilon _{t}(\psi _{n},\hat{\pi}_{n})d_{\psi ,t}(%
\hat{\pi}_{n})\right] \right\}  \notag \\
&&\text{ \ \ \ \ \ \ \ \ \ \ \ \ \ \ \ }+\mathfrak{b}_{\psi }(\hat{\pi}%
_{n},\lambda )^{\prime }\left\{ \mathcal{H}_{\psi }^{-1}(\hat{\pi}_{n})%
\mathcal{D}_{\psi }(\hat{\pi}_{n})b+\left[ b,0_{k_{\beta }}^{\prime }\right]
^{\prime }\right\} +o_{p,\lambda }\left( 1\right)  \notag \\
&&\text{ \ \ \ }=\frac{1}{\sqrt{n}}\sum_{t=1}^{n}\epsilon _{t}(\psi _{n},%
\hat{\pi}_{n})\mathcal{K}_{\psi ,t}(\hat{\pi}_{n},\lambda )  \label{sum} \\
&&\text{ \ \ \ \ \ \ \ \ \ \ \ \ \ \ \ \ \ \ \ }+\mathfrak{b}_{\psi }(\hat{%
\pi}_{n},\lambda )^{\prime }\mathcal{H}_{\psi }^{-1}(\hat{\pi}_{n})\sqrt{n}E%
\left[ \epsilon _{t}(\psi _{n},\hat{\pi}_{n})d_{\psi ,t}(\hat{\pi}_{n})%
\right]  \label{mean} \\
\text{\ \ \ \ } &&\text{\ \ \ \ \ \ \ \ \ \ \ \ \ \ \ \ \ \ \ \ }+\mathfrak{b%
}_{\psi }(\hat{\pi}_{n},\lambda )^{\prime }\left\{ \mathcal{H}_{\psi }^{-1}(%
\hat{\pi}_{n})\mathcal{D}_{\psi }(\hat{\pi}_{n})b+\left[ b,0_{k_{\beta
}}^{\prime }\right] ^{\prime }\right\} +o_{p,\lambda }\left( 1\right) .
\notag
\end{eqnarray}%
By adding and subtracting $E[\epsilon _{t}(\psi _{n},\hat{\pi}_{n})\mathcal{K%
}_{\psi ,t}(\hat{\pi}_{n},\lambda )]$, summand (\ref{sum}) satisfies:
\begin{eqnarray}
&&\frac{1}{\sqrt{n}}\sum_{t=1}^{n}\epsilon _{t}(\psi _{n},\hat{\pi}_{n})%
\mathcal{K}_{\psi ,t}(\hat{\pi}_{n},\lambda )  \label{sum1} \\
&&\text{ \ \ \ }=\frac{1}{\sqrt{n}}\sum_{t=1}^{n}\left\{ \epsilon _{t}(\psi
_{n},\hat{\pi}_{n})\mathcal{K}_{\psi ,t}(\hat{\pi}_{n},\lambda )-E\left[
\epsilon _{t}(\psi _{n},\hat{\pi}_{n})\mathcal{K}_{\psi ,t}(\hat{\pi}%
_{n},\lambda )\right] \right\} +\sqrt{n}E\left[ \epsilon _{t}(\psi _{n},\hat{%
\pi}_{n})\mathcal{K}_{\psi ,t}(\hat{\pi}_{n},\lambda )\right]  \notag \\
&&\text{ \ \ \ }=\frac{1}{\sqrt{n}}\sum_{t=1}^{n}\left\{ \epsilon _{t}(\psi
_{n},\hat{\pi}_{n})\mathcal{K}_{\psi ,t}(\hat{\pi}_{n},\lambda )-E\left[
\epsilon _{t}(\psi _{n},\hat{\pi}_{n})\mathcal{K}_{\psi ,t}(\hat{\pi}%
_{n},\lambda )\right] \right\}  \notag \\
&&\text{ \ \ \ \ \ \ \ \ \ \ \ \ \ \ \ }+\sqrt{n}E\left[ \epsilon _{t}%
\mathcal{K}_{\psi ,t}(\hat{\pi}_{n},\lambda )\right] +E\left[ \mathcal{K}%
_{\psi ,t}(\hat{\pi}_{n},\lambda )\left\{ g(x_{t},\pi _{0})-g(x_{t},\hat{\pi}%
_{n})\right\} ^{\prime }\right] \sqrt{n}\beta _{n}.  \notag
\end{eqnarray}%
Under $H_{0}$, trivially $\sup_{\pi \in \Pi }||E[\epsilon _{t}\mathcal{K}%
_{\psi ,t}(\pi ,\lambda )]||$ $=$ $0$. Turning to the expectations in (\ref%
{mean}):%
\begin{eqnarray}
E\left[ \epsilon _{t}(\psi _{n},\hat{\pi}_{n})d_{\psi ,t}(\hat{\pi}_{n})%
\right] &=&E\left[ \epsilon _{t}d_{\psi ,t}(\hat{\pi}_{n})\right] +E\left[
d_{\psi ,t}(\hat{\pi}_{n})\left\{ g(x_{t},\pi _{0})-g(x_{t},\hat{\pi}%
_{n})\right\} ^{\prime }\right] \beta _{n}  \label{mean1} \\
&=&E\left[ d_{\psi ,t}(\hat{\pi}_{n})\left\{ g(x_{t},\pi _{0})-g(x_{t},\hat{%
\pi}_{n})\right\} ^{\prime }\right] \beta _{n}.  \notag
\end{eqnarray}%
Combine (\ref{sumEF}), (\ref{sum1}) and (\ref{mean1}) with $\sqrt{n}\beta
_{n}$ $\rightarrow $ $b$, $||b||$ $<$ $\infty $, to arrive at (\ref%
{eF_expand}).\medskip \newline
\textbf{Step 2.}\qquad We will show $\{1/\sqrt{n}\sum_{t=1}^{n}\epsilon _{t}(%
\hat{\theta}_{n})F(\lambda ^{\prime }\mathcal{W}(x_{t}))$ $:$ $\Lambda \}$ $%
\Rightarrow ^{\ast }$ $\{\mathfrak{Z}_{\psi }(\pi ^{\ast }(b),\lambda )$ $+$
$\mathcal{R}(\pi ^{\ast }(b),\lambda )$ $:$ $\Lambda \}$. Lemma \ref%
{lm:CM_weak}.a states $\{\mathfrak{Z}_{n}(\pi ,\lambda )$ $:$ $\Pi ,\Lambda
\}$ $\Rightarrow ^{\ast }$ $\{\mathfrak{Z}_{\psi }(\pi ,\lambda )$ $:$ $\Pi
,\Lambda \}$, a zero mean Gaussian process, and by Theorem \ref{th:ols}.a $%
\hat{\pi}_{n}$ $\overset{d}{\rightarrow }$ $\pi ^{\ast }(b)$ where $\pi
^{\ast }(b)$ is defined by Assumption \ref{assum:pi}. Step 3 proves joint
weak convergence%
\begin{equation}
\left\{ \mathfrak{Z}_{n}(\pi ,\lambda ),\hat{\pi}_{n}:\pi \in \Pi ,\lambda
\in \Lambda \right\} \Rightarrow ^{\ast }\left\{ \mathfrak{Z}_{\psi }(\pi
,\lambda ),\pi ^{\ast }(b):\pi \in \Pi ,\lambda \in \Lambda \right\} .
\label{Zpi_joint}
\end{equation}%
The mapping theorem and expansion (\ref{eF_expand}) deliver the desired
result.\medskip \newline
\textbf{Step 3}.\qquad We need to show (\ref{Zpi_joint}). By the proof of
Theorem \ref{th:ols}.a, $\hat{\pi}_{n}$ is a continuous function of $%
\mathcal{G}_{\psi ,n}(\psi _{0,n},\pi )$ and $\widehat{\mathcal{H}}_{\psi
,n}(\pi )$. Further, $\{\mathcal{G}_{\psi ,n}(\psi _{0,n},\pi )$ $:$ $\Pi \}$
$\Rightarrow ^{\ast }$ $\{\mathcal{G}_{\psi }(\pi )$ $:$ $\Pi \}$ by Lemma %
\ref{lm:G_uclt_weak}, and $\widehat{\mathcal{H}}_{\psi ,n}(\pi )$ has a
non-random limit uniformly on $\Pi $ by Lemma \ref{lm:H_ulln_weak}.
Therefore, (\ref{Zpi_joint}) follows from the mapping theorem and Cram\'{e}%
r's theorem provided jointly:
\begin{eqnarray*}
&&\left\{ \left[
\begin{array}{c}
\mathfrak{Z}_{n}(\pi ,\lambda ) \\
\mathcal{G}_{\psi ,n}(\psi _{0,n},\pi )%
\end{array}%
\right] :\Pi ,\Lambda \right\} =\left\{ \left[
\begin{array}{c}
\dfrac{1}{\sqrt{n}}\sum_{t=1}^{n}\left\{ \epsilon _{t}(\pi ,\lambda )%
\mathcal{K}_{\psi ,t}(\pi ,\lambda )-E\left[ \epsilon _{t}(\psi _{n},\pi )%
\mathcal{K}_{\psi ,t}(\pi ,\lambda )\right] \right\} \\
-\dfrac{1}{\sqrt{n}}\sum_{t=1}^{n}\left\{ \epsilon _{t}(\psi ,\pi )d_{\psi
,t}(\pi )-E\left[ \epsilon _{t}(\psi ,\pi )d_{\psi ,t}(\pi )\right] \right\}%
\end{array}%
\right] :\Pi ,\Lambda \right\} \\
&&\text{ \ \ \ \ \ \ \ \ \ \ \ \ \ \ \ \ \ \ \ \ \ \ \ \ \ \ \ \ \ \ \ \ \ \
\ \ }\Rightarrow ^{\ast }\left\{ \left[
\begin{array}{c}
\mathfrak{Z}_{\psi }(\pi ,\lambda ) \\
\mathcal{G}_{\psi }(\pi )%
\end{array}%
\right] :\Pi ,\Lambda \right\} .
\end{eqnarray*}%
The latter holds by the same arguments used to prove Lemmas \ref%
{lm:G_uclt_weak} and \ref{lm:CM_weak}, hence we only provide a sketch of the
proof. First, $[\mathfrak{Z}_{n}(\pi ,\lambda ),\mathcal{G}_{\psi ,n}(\psi
_{0,n},\pi )^{\prime }]^{\prime }$ converges in finite dimensional
distributions over $\Pi $ $\times $ $\Lambda $ to a zero mean, finite
variance Gaussian random vector. This follows because linear combinations $%
\sum_{i=1}^{m}a_{i}\{r_{1}\mathfrak{Z}_{n}(\pi _{i},\lambda _{i})$ $+$ $%
r_{2}^{\prime }\mathcal{G}_{\psi ,n}(\psi _{0,n},\pi _{i})\}$ for any $m$ $%
\in $ $\mathbb{N}$ and $a$ $\in $ $\mathbb{R}^{m}$ with $a^{\prime }a$ $=$ $%
1 $, and any $r$ $=$ $[r_{1},r_{2}^{\prime }]^{\prime }$, $r^{\prime }r$ $=$
$1 $, satisfy a Gaussian central theorem under the moment and memory
properties of Assumption \ref{assum:dgp}.b,c,d(vi). Second, $[\mathfrak{Z}%
_{n}(\pi ,\lambda ),\mathcal{G}_{\psi ,n}(\psi _{0,n},\pi )^{\prime
}]^{\prime }$ is stochastically equicontinuous because, by probability
subadditivity, we require $\mathfrak{Z}_{n}(\pi ,\lambda )$ and each element
$\mathcal{G}_{\psi ,n,i}(\psi _{0,n},\pi )$ of $\mathcal{G}_{\psi ,n}(\psi
_{0,n},\pi )$ $=$ $[\mathcal{G}_{\psi ,n,i}(\psi _{0,n},\pi
)]_{i=1}^{k_{x}+k_{\beta }}$ to be stochastically equicontinuous, and these
properties are established in the proofs of Lemmas \ref{lm:G_uclt_weak} and %
\ref{lm:CM_weak}.$\medskip $\newline
\textbf{Step 4.}\qquad We now tackle $\hat{v}_{n}(\hat{\theta}_{n},\lambda )$
and complete the proof. $\hat{v}_{n}(\hat{\theta}_{n},\lambda )$ is a
function of $\mathfrak{\hat{b}}_{\theta ,n}(\omega (\hat{\beta}_{n}),\hat{\pi%
}_{n},\lambda )$ and $d_{\theta ,t}(\omega (\hat{\beta}_{n}),\hat{\pi}_{n})$%
. By Lemma \ref{lm:vn}.a $\sup_{\theta ^{+}\in \Theta ^{+},\lambda \in
\Lambda }||\hat{v}_{n}^{2}(\theta ^{+},\lambda )-$ $v^{2}(\theta
^{+},\lambda )||$ $\overset{p}{\rightarrow }$ $0$. Furthermore,
\begin{equation}
\omega (\hat{\beta}_{n}(\hat{\pi}_{n}))=\frac{\sqrt{n}\mathcal{S}_{\beta }%
\hat{\psi}_{n}(\hat{\pi}_{n})}{\left\Vert \sqrt{n}\mathcal{S}_{\beta }\hat{%
\psi}_{n}(\hat{\pi}_{n})\right\Vert }=\frac{\sqrt{n}\mathcal{S}_{\beta
}\left( \hat{\psi}_{n}(\hat{\pi}_{n})-\psi _{n}\right) +\sqrt{n}\beta _{n}}{%
\left\Vert \sqrt{n}\mathcal{S}_{\beta }\left( \hat{\psi}_{n}(\hat{\pi}%
_{n})-\psi _{n}\right) +\sqrt{n}\beta _{n}\right\Vert }\equiv \omega _{n}(%
\hat{\pi}_{n}),  \label{wwn}
\end{equation}%
by construction of $\mathcal{S}_{\beta }$. Notice $\omega _{n}$\ is an
implicitly defined stochastic function of $\hat{\pi}_{n}$. Now, by Theorem %
\ref{th:ols}.a, the mapping theorem and $\sqrt{n}\beta _{n}$ $\rightarrow $ $%
b$, $||b||$ $<$ $\infty $:
\begin{equation*}
\omega _{n}(\hat{\pi}_{n})\overset{d}{\rightarrow }\frac{\mathcal{S}_{\beta
}\tau (\pi ^{\ast }(b),b)}{\left\Vert \mathcal{S}_{\beta }\tau (\pi ^{\ast
}(b),b)\right\Vert }=\frac{\tau _{\beta }(\pi ^{\ast }(b),b)}{\left\Vert
\tau _{\beta }(\pi ^{\ast }(b),b)\right\Vert }.
\end{equation*}%
Joint weak convergence for $\sqrt{n}(\hat{\psi}_{n}(\pi )$ $-$ $\psi _{n}),$
$\hat{\pi}_{n}$ and $\omega _{n}(\hat{\pi}_{n})$ follows from arguments in
the proof of Theorem \ref{th:ols}.a because $\omega _{n}(\pi )$ is a
continuous function of $\sqrt{n}(\hat{\psi}_{n}(\pi )$ $-$ $\psi _{n})$, and
$\sqrt{n}(\hat{\psi}_{n}(\pi )$ $-$ $\psi _{n})$ and $\hat{\pi}_{n}$ are
continuous functions of $\mathcal{G}_{\psi ,n}(\psi _{0,n},\pi )$ and $%
\widehat{\mathcal{H}}_{\psi ,n}(\pi )$. Hence:%
\begin{equation}
\left[ \sqrt{n}\left( \hat{\psi}_{n}(\hat{\pi}_{n})-\psi _{n}\right)
^{\prime },\hat{\pi}_{n}^{\prime },\omega _{n}(\hat{\pi}_{n})^{\prime }%
\right] ^{\prime }\overset{d}{\rightarrow }\left[ \tau (\pi ^{\ast
}(b),b)^{\prime },\pi ^{\ast }(b)^{\prime },\frac{\tau _{\beta }(\pi ^{\ast
}(b),b)^{\prime }}{\left\Vert \tau _{\beta }(\pi ^{\ast }(b),b)\right\Vert }%
\right] ^{\prime }.  \label{psi_pi_om}
\end{equation}%
Using $v^{2}(\cdot )$ defined in (\ref{v2}), we may therefore write:
\begin{equation}
\hat{v}_{n}^{2}(\hat{\theta}_{n},\lambda )=v^{2}(\omega _{n}(\hat{\pi}_{n}),%
\hat{\pi}_{n},\lambda )+o_{p,\lambda }(1)\text{ where }v^{2}(\omega ,\pi
,\lambda )=v^{2}(\left[ \beta _{0},\omega ,\zeta _{0},\pi \right] ,\lambda ),
\label{vnv}
\end{equation}%
and $\lim \inf_{n\rightarrow \infty }v^{2}(\omega _{n}(\hat{\pi}_{n}),\hat{%
\pi}_{n},\lambda )$ $>$ $0$ $a.s.$ $\forall \lambda $ $\in $ $\Lambda $\ by
Assumption \ref{assum:nondegen_v2_ae_all}. The claim now follows from Step
2, $\sup_{\theta ^{+}\in \Theta ^{+},\lambda \in \Lambda }||\hat{v}%
_{n}^{2}(\theta ^{+},\lambda )-$ $v^{2}(\theta ^{+},\lambda )||$ $\overset{p}%
{\rightarrow }$ $0$, (\ref{psi_pi_om}), (\ref{vnv}) and the mapping theorem:

\begin{eqnarray}
\left\{ \mathcal{T}_{n}(\lambda ):\Lambda \right\} &=&\left\{ \frac{\left(
\mathfrak{Z}_{n}(\hat{\pi}_{n},\lambda )+\mathcal{R}(\hat{\pi}_{n},\lambda
)\right) ^{2}}{v^{2}(\omega _{n}(\hat{\pi}_{n}),\hat{\pi}_{n},\lambda )}%
+o_{p,\lambda }\left( 1\right) :\Lambda \right\}  \label{Tn_expand} \\
&\Rightarrow &^{\ast }\left\{ \frac{\left( \mathfrak{Z}(\pi ^{\ast
}(b),\lambda )+\mathcal{R}(\pi ^{\ast }(b),\lambda )\right) ^{2}}{%
v^{2}(\omega (\pi ^{\ast }(b)),\pi ^{\ast }(b),\lambda )}:\Lambda \right\} .
\notag
\end{eqnarray}%
\textbf{Claim b.}\qquad Let $\mathcal{C}(ii,\omega _{0})$ apply. A first
order expansion yields for some midpoint $\theta _{n}^{\ast }$, $||\theta
_{n}^{\ast }$ $-$ $\theta _{n}||$ $\leq $ $||\hat{\theta}_{n}$ $-$ $\theta
_{n}||$:
\begin{eqnarray*}
\frac{1}{\sqrt{n}}\sum_{t=1}^{n}\epsilon _{t}(\hat{\theta}_{n})F\left(
\lambda ^{\prime }\mathcal{W}(x_{t})\right) &=&\frac{1}{\sqrt{n}}%
\sum_{t=1}^{n}\epsilon _{t}F\left( \lambda ^{\prime }\mathcal{W}%
(x_{t})\right) -\frac{1}{n}\sum_{t=1}^{n}F\left( \lambda ^{\prime }\mathcal{W%
}(x_{t})\right) \frac{\partial }{\partial \theta ^{\prime }}f(x_{t,},\theta
_{n}^{\ast })\sqrt{n}\left( \hat{\theta}_{n}-\theta _{n}\right) \\
&=&\frac{1}{\sqrt{n}}\sum_{t=1}^{n}\epsilon _{t}F\left( \lambda ^{\prime }%
\mathcal{W}(x_{t})\right) -\mathfrak{\hat{b}}_{\theta ,n}(\beta _{n}^{\ast
}/\left\Vert \beta _{n}^{\ast }\right\Vert ,\pi _{n}^{\ast },\lambda
)^{\prime }\sqrt{n}\mathfrak{B}(\beta _{n})\left( \hat{\theta}_{n}-\theta
_{n}\right) .
\end{eqnarray*}%
The proof of Theorem \ref{th:ols}.b shows
\citep[see (C.21)
in][]{Supp_Mat_2020}:
\begin{equation*}
\sqrt{n}\mathfrak{B}(\beta _{n})\left( \hat{\theta}_{n}-\theta _{n}\right) =%
\widehat{\mathcal{H}}_{n}^{-1}(\omega (\beta _{n}^{\ast }),\pi _{n}^{\ast })%
\mathfrak{B}(\beta _{n})^{-1}\sqrt{n}\frac{\partial }{\partial \theta }%
Q_{n}(\theta _{n})=\widehat{\mathcal{H}}_{n}^{-1}(\omega (\beta _{n}^{\ast
}),\pi _{n}^{\ast })\mathfrak{B}(\beta _{n})^{-1}\mathcal{G}_{\theta
,n}(\theta _{n}).
\end{equation*}%
By Lemma \ref{lm:JV} $\sup_{\omega \in \mathbb{R}^{k_{\beta }}:||\omega
||=1,\pi \in \Pi }||\widehat{\mathcal{H}}_{n}(\omega ,\pi )$ $-$ $\mathcal{H}%
_{\theta }(\omega ,\pi )||$ $\overset{p}{\rightarrow }$ $0$, and $\hat{\theta%
}_{n}$ $\overset{p}{\rightarrow }$ $\theta _{0}$ by the proof of Theorem \ref%
{th:ols}. Thus, by definition of $\mathcal{G}_{\theta ,n}(\cdot )$:
\begin{equation*}
\sqrt{n}\mathfrak{B}(\beta _{n})\left( \hat{\theta}_{n}-\theta _{n}\right) =%
\mathcal{H}_{\theta }^{-1}\mathfrak{B}(\beta _{n})^{-1}\mathcal{G}_{\theta
,n}(\theta _{n})+o_{p}(1)=-\mathcal{H}_{\theta }^{-1}\frac{1}{\sqrt{n}}%
\sum_{t=1}^{n}\epsilon _{t}d_{\theta ,t}(\beta _{n}/\left\Vert \beta
_{n}\right\Vert ,\pi _{0})+o_{p}(1).
\end{equation*}%
Moreover, $\sup_{\lambda \in \Lambda }||\mathfrak{\hat{b}}_{\theta ,n}(\beta
_{n}^{\ast }/||\beta _{n}^{\ast }||,\hat{\pi}_{n},\lambda )$ $-$ $\mathfrak{b%
}_{\theta }(\omega _{0},\pi _{0},\lambda )||$ $\overset{p}{\rightarrow }$ $0$
by Lemma \ref{lm:bn} and $\hat{\theta}_{n}$ $\overset{p}{\rightarrow }$ $%
\theta _{0}$. Combined, we obtain:%
\begin{eqnarray}
&&\frac{1}{\sqrt{n}}\sum_{t=1}^{n}\epsilon _{t}(\hat{\theta}_{n})F\left(
\lambda ^{\prime }\mathcal{W}(x_{t})\right)  \label{eF_expand_strong} \\
&&\text{ \ \ \ \ \ \ \ \ \ }=\frac{1}{\sqrt{n}}\sum_{t=1}^{n}\epsilon
_{t}F\left( \lambda ^{\prime }\mathcal{W}(x_{t})\right) +\mathfrak{b}%
_{\theta }(\omega _{0},\pi _{0},\lambda )^{\prime }\mathcal{H}_{\theta }^{-1}%
\frac{1}{\sqrt{n}}\sum_{t=1}^{n}\epsilon _{t}d_{\theta ,t}(\beta
_{n}/\left\Vert \beta _{n}\right\Vert ,\pi _{0})+o_{p,\lambda }(1)  \notag \\
&&\text{ \ \ \ \ \ \ \ \ \ }=\frac{1}{\sqrt{n}}\sum_{t=1}^{n}\epsilon
_{t}\left( F\left( \lambda ^{\prime }\mathcal{W}(x_{t})\right) +\mathfrak{b}%
_{\theta }(\omega _{0},\pi _{0},\lambda )^{\prime }\mathcal{H}_{\theta
}^{-1}d_{\theta ,t}(\beta _{n}/\left\Vert \beta _{n}\right\Vert ,\pi
_{0})\right) +o_{p,\lambda }(1).  \notag
\end{eqnarray}%
By Lemma \ref{lm:CM_weak}.b, therefore, $\{1/\sqrt{n}\sum_{t=1}^{n}\epsilon
_{t}(\hat{\theta}_{n})F\left( \lambda ^{\prime }\mathcal{W}(x_{t})\right) $ $%
:$ $\Lambda \}$ $\Rightarrow ^{\ast }$ $\{\mathfrak{Z}_{\theta }(\lambda )$ $%
:$ $\Lambda \}$, a zero mean Gaussian process with covariance kernel $E[%
\mathfrak{Z}_{\theta }(\lambda )\mathfrak{Z}_{\theta }(\tilde{\lambda})]$ $=$
$E[\epsilon _{t}^{2}\mathcal{K}_{\theta ,t}(\lambda )\mathcal{K}_{\theta ,t}(%
\tilde{\lambda})]$ where $\mathcal{K}_{\theta ,t}(\lambda )\equiv $ $F\left(
\lambda ^{\prime }\mathcal{W}(x_{t})\right) -\mathfrak{b}_{\theta }(\omega
_{0},\pi _{0},\lambda )^{\prime }\mathcal{H}_{\theta }^{-1}d_{\theta
,t}(\omega _{0},\pi _{0})$.

Now turn to $\hat{v}_{n}(\hat{\theta}_{n},\lambda )$. By Lemma \ref{lm:vn}.b
and $\hat{\theta}_{n}$ $\overset{p}{\rightarrow }$ $\theta _{0}$, $%
\sup_{\lambda \in \Lambda }||\hat{v}_{n}^{2}(\hat{\theta}_{n},\lambda )$ $-$
$v^{2}(\theta _{0},\lambda )||$ $\overset{p}{\rightarrow }0$ where by
construction and Assumption \ref{assum:dgp}.d(vi) $v^{2}(\theta _{0},\lambda
)$ $=$ $E[\epsilon _{t}^{2}\mathcal{K}_{\theta ,t}^{2}(\lambda )]$ $<$ $%
\infty $. By Assumption \ref{assum:nondegen_v2_ae_all} $v^{2}(\theta
_{0},\lambda )$ $>$ $0$ $\forall \lambda $ $\in $ $\Lambda $. Since $E[%
\mathfrak{Z}_{\theta }^{2}(\lambda )]$ $=$ $E[\epsilon _{t}^{2}\mathcal{K}%
_{\theta ,t}^{2}(\lambda )]$, the proof is complete by the mapping theorem. $%
\mathcal{QED}$.\bigskip \newline
\textbf{Proof of Theorem \ref{th:CM_H1}.}\qquad By the proof of Theorem of %
\ref{th:CM_weak}, $\hat{v}_{n}^{2}(\hat{\theta}_{n},\lambda )$ $\rightarrow $
$(0,\infty )$ asymptotically with probability approaching one under any
identification case, and $\forall \lambda $ $\in $ $\Lambda $.

It remains to prove $|1/n\sum_{t=1}^{n}\epsilon _{t}(\hat{\theta}%
_{n})F(\lambda ^{\prime }\mathcal{W}(x_{t}))|$ $\overset{p}{\rightarrow }$ $%
(0,\infty )$ $\forall \lambda $ $\in $ $\Lambda /S$ where $S$ $\subset $ $%
\Lambda $ has Lebesgue measure zero. Consider identification case $\mathcal{C%
}(i,b)$ with $||b||$ $<$ $\infty $. Arguments in the proof of Theorem \ref%
{th:CM_weak}.a imply
\begin{eqnarray*}
\frac{1}{n}\sum_{t=1}^{n}\epsilon _{t}(\hat{\theta}_{n})F\left( \lambda
^{\prime }\mathcal{W}(x_{t})\right) &=&\frac{1}{n}\sum_{t=1}^{n}\epsilon
_{t}(\psi _{n},\hat{\pi}_{n})F\left( \lambda ^{\prime }\mathcal{W}%
(x_{t})\right) -\mathfrak{\hat{b}}_{\psi ,n}(\hat{\pi}_{n},\lambda )^{\prime
}\left( \hat{\psi}_{n}-\psi _{n}\right) \\
&=&\frac{1}{n}\sum_{t=1}^{n}\epsilon _{t}(\psi _{n},\hat{\pi}_{n})F\left(
\lambda ^{\prime }\mathcal{W}(x_{t})\right) +O_{p,\lambda }(1/\sqrt{n}) \\
\text{\ } &=&E\left[ \epsilon _{t}F\left( \lambda ^{\prime }\mathcal{W}%
(x_{t})\right) \right] +\frac{1}{n}\sum_{t=1}^{n}\left\{ \epsilon
_{t}F\left( \lambda ^{\prime }\mathcal{W}(x_{t})\right) -E\left[ \epsilon
_{t}F\left( \lambda ^{\prime }\mathcal{W}(x_{t})\right) \right] \right\} \\
&&-\beta _{n}^{\prime }\frac{1}{n}\sum_{t=1}^{n}\left\{ g(x_{t},\hat{\pi}%
_{n})-g(x_{t},\pi _{0})\right\} F\left( \lambda ^{\prime }\mathcal{W}%
(x_{t})\right) +O_{p,\lambda }(1/\sqrt{n}).
\end{eqnarray*}%
Under $\sqrt{n}||\beta _{n}||$ $\rightarrow $ $[0,\infty )$, by Lemma \ref%
{lm:ulln_xgF_eF} we have $\sup_{\pi \in \Pi ,\lambda \in \Lambda
}||1/n\sum_{t=1}^{n}\epsilon _{t}(\hat{\theta}_{n})F(\lambda ^{\prime }%
\mathcal{W}(x_{t}))$ $-$ \linebreak $E[\epsilon _{t}F(\lambda ^{\prime }%
\mathcal{W}(x_{t}))]||$ $\overset{p}{\rightarrow }$ $0$. The claim now
follows from $E[\epsilon _{t}F(\lambda ^{\prime }\mathcal{W}(x_{t}))]$ $\neq
$ $0$ $\forall \lambda $ $\in $ $\Lambda /S$ for some $S$ $\subset $ $%
\Lambda $ with Lebesgue measure zero by Theorem 2.3 in \cite{StinchWhite1998}%
, cf. \citet[Lemma 1]{Bierens1990}.

Under $\mathcal{C}(ii,\omega _{0})$, and the proof of Theorem \ref%
{th:CM_weak}.b:
\begin{eqnarray*}
\sum_{t=1}^{n}\epsilon _{t}(\hat{\theta}_{n})F(\lambda ^{\prime }\mathcal{W}%
(x_{t}))-E\left[ \epsilon _{t}F(\lambda ^{\prime }\mathcal{W}(x_{t}))\right]
&=&\frac{1}{n}\sum_{t=1}^{n}\left\{ \epsilon _{t}F\left( \lambda ^{\prime }%
\mathcal{W}(x_{t})\right) -E\left[ \epsilon _{t}F\left( \lambda ^{\prime }%
\mathcal{W}(x_{t})\right) \right] \right\} \\
&&+\mathfrak{b}_{\theta }(\omega _{0},\pi _{0},\lambda )^{\prime }\mathcal{H}%
_{\theta }^{-1}\frac{1}{n}\sum_{t=1}^{n}\epsilon _{t}d_{\theta ,t}(\beta
_{n}/\left\Vert \beta _{n}\right\Vert ,\pi _{0})+O_{p,\lambda }(1/\sqrt{n}).
\end{eqnarray*}%
Mixing under Assumption \ref{assum:dgp}.b implies ergodicity. The ergodic
theorem therefore yields \linebreak $1/n\sum_{t=1}^{n}\epsilon _{t}d_{\theta
,t}(\beta _{n}/\left\Vert \beta _{n}\right\Vert ,\pi _{0})$ $\overset{p}{%
\rightarrow }$ $0$ in view of identification Assumption \ref{assum:dgp}%
.a(ii) and stationarity. By Lemma \ref{lm:ulln_xgF_eF}, $1/n\sum_{t=1}^{n}\{%
\epsilon _{t}F(\lambda ^{\prime }\mathcal{W}(x_{t}))$ $-$ $E[\epsilon
_{t}F(\lambda ^{\prime }\mathcal{W}(x_{t}))]\}$ $\overset{p}{\rightarrow }$ $%
0$, and Theorem 2.3 in \cite{StinchWhite1998} yields $E[\epsilon
_{t}F(\lambda ^{\prime }\mathcal{W}(x_{t}))]$ $\neq $ $0$ $\forall \lambda $
$\in $ $\Lambda /S$. $\mathcal{QED}.$\bigskip \newline
\textbf{Proof of Theorem \ref{th:ICS_stat}.}\medskip \newline
\textbf{Claim a.}\qquad Let $\mathcal{C}(i,b)$ hold.\medskip

\textbf{Step 1}.\qquad We first show $||\hat{\Sigma}_{n}$ $-$ $\bar{\Sigma}%
(\pi ^{\ast }(b),b)||$ $\overset{p}{\rightarrow }$ $0$, where $\bar{\Sigma}%
(\pi ,b)$ $\equiv $ $\Sigma (\omega ^{\ast }(\pi ,b),\pi )$ $=$ $\Sigma
(||\beta _{0}||,\omega ^{\ast }(\pi ,b),$ $\zeta _{0},\pi )$. Recall $\Sigma
(||\beta ||,\omega ,\zeta ,\pi )$ $=$ $\Sigma (\theta ^{+})$ $\equiv $ $%
\mathcal{H}_{\theta }(\theta ^{+})^{-1}\mathcal{V}(\theta ^{+})\mathcal{H}%
_{\theta }(\theta ^{+})^{-1}$, cf. (\ref{Sigma_the+}).\ Recall $\hat{\beta}%
_{n}$ $=$ $\hat{\beta}_{n}(\hat{\pi}_{n})$ hence we write $\omega (\hat{\beta%
}_{n})$ $=$ $\omega _{n}(\hat{\pi}_{n})$ as in (\ref{wwn}). Thus joint weak
convergence for $\sqrt{n}(\hat{\psi}_{n}(\hat{\pi}_{n})$ $-$ $\psi _{n}),$ $%
\hat{\pi}_{n}$ and $\omega (\hat{\beta}_{n})$ holds by (\ref{psi_pi_om}).
Hence:%
\begin{equation*}
\left[ \sqrt{n}\left( \hat{\psi}_{n}(\hat{\pi}_{n})-\psi _{n}\right)
^{\prime },\hat{\pi}_{n}^{\prime },\omega (\hat{\beta}_{n})^{\prime }\right]
^{\prime }\overset{d}{\rightarrow }\left[ \tau (\pi ^{\ast }(b),b)^{\prime
},\pi ^{\ast }(b)^{\prime },\frac{\tau _{\beta }(\pi ^{\ast }(b),b)^{\prime }%
}{\left\Vert \tau _{\beta }(\pi ^{\ast }(b),b)\right\Vert }\right] ^{\prime
}.
\end{equation*}%
Joint convergence for $\hat{\pi}_{n}$ and $\omega (\hat{\beta}_{n})$, $%
\sup_{\pi \in \Pi }||\hat{\psi}_{n}(\pi )$ $-$ $\psi _{n}||$ $\overset{p}{%
\rightarrow }$ $0$ by the proof of Theorem \ref{th:ols}, uniform convergence
Lemma \ref{lm:JV} for $\widehat{\mathcal{H}}_{n}(\theta ^{+})$\ and $%
\mathcal{\hat{V}}_{n}(\theta ^{+})$, and the mapping theorem, together yield
$||\hat{\Sigma}_{n}$ $-$ $\bar{\Sigma}(\pi ^{\ast }(b),b)||$ $\overset{p}{%
\rightarrow }$ $0$.\medskip

\textbf{Step 2}.\qquad Now invoke $(\sqrt{n}(\hat{\psi}_{n}(\hat{\pi}_{n})$ $%
-$ $\psi _{n}),\hat{\pi}_{n})$ $\overset{d}{\rightarrow }$ $(\tau (\pi
^{\ast }(b),b),\pi ^{\ast }(b))$ by Theorem \ref{th:ols}.a, and the mapping
theorem, to complete the proof.\medskip \newline
\textbf{Claim b.}\qquad Let $\mathcal{C}(ii,\omega _{0})$ hold, and let $%
\{\kappa _{n}\}$\ be a sequence of positive constants, $\kappa _{n}$ $%
\rightarrow $ $\infty $ and $\kappa _{n}$ $=$ $o(\sqrt{n})$. Write:
\begin{eqnarray*}
&&\kappa _{n}^{-2}\mathcal{A}_{n}^{2}=\frac{1}{k_{\beta }}\frac{n}{\kappa
_{n}^{2}}\left( \hat{\beta}_{n}-\beta _{n}\right) ^{\prime }\hat{\Sigma}%
_{\beta ,\beta ,n}^{-1}\left( \hat{\beta}_{n}-\beta _{n}\right) \\
&&\text{ \ \ \ \ \ \ \ \ \ \ \ \ \ \ \ }+2\frac{1}{k_{\beta }}\frac{\sqrt{n}%
}{\kappa _{n}}\left( \hat{\beta}_{n}-\beta _{n}\right) ^{\prime }\hat{\Sigma}%
_{\beta ,\beta ,n}^{-1}\frac{\sqrt{n}}{\kappa _{n}}\beta _{n}+\frac{1}{%
k_{\beta }}\frac{\sqrt{n}}{\kappa _{n}}\beta _{n}^{\prime }\hat{\Sigma}%
_{\beta ,\beta ,n}^{-1}\frac{\sqrt{n}}{\kappa _{n}}\beta _{n}\equiv \mathcal{%
B}_{n,1}+\mathcal{B}_{n,2}+\mathcal{B}_{n,3}.
\end{eqnarray*}%
Lemma \ref{lm:JV} yields that $\hat{\Sigma}_{\beta ,\beta ,n}$ is positive
definite asymptotically with probability approaching one. Hence $\mathcal{B}%
_{n,1}\overset{p}{\rightarrow }$ $0$ by application of Theorem \ref{th:ols}%
.b. If $\sqrt{n}||\beta _{n}||/\kappa _{n}$ $=$ $O(1)$ then $|\mathcal{B}%
_{n,2}|$ $\overset{p}{\rightarrow }$ $0$ because $\sqrt{n}\kappa _{n}^{-1}(%
\hat{\beta}_{n}$ $-\beta _{n})$ $=$ $O_{p}(1/\kappa _{n})$ $=$ $o_{p}(1)$,
and $\mathcal{B}_{n,3}\overset{p}{\rightarrow }$ $[0,\infty )$, hence $%
\kappa _{n}^{-1}\mathcal{A}_{n}$ $\overset{p}{\rightarrow }$ $[0,\infty )$.

If $\sqrt{n}||\beta _{n}||/\kappa _{n}$ $\rightarrow $ $\infty $ then $|%
\mathcal{B}_{n,2}|$ $=$ $O_{p}(\mathcal{B}_{n,3})$ because $\sqrt{n}\kappa
_{n}^{-1}(\hat{\beta}_{n}$ $-\beta _{n})$ $=$ $o_{p}(1)$, and $\mathcal{B}%
_{n,3}\overset{p}{\rightarrow }$ $\infty $, hence $\kappa _{n}^{-1}\mathcal{A%
}_{n}$ $\overset{p}{\rightarrow }$ $\infty $. An example of this final case
is $\beta _{0}$ $\neq $ $0$: $||\beta _{n}||$ $\rightarrow $ $||\beta _{0}||$
$>$ $0$ while $\kappa _{n}$ $=$ $o(\sqrt{n})$ by supposition hence $\sqrt{n}%
||\beta _{n}||/\kappa _{n}$ $\rightarrow $ $\infty $. $\mathcal{QED}$%
.\bigskip \newline
\textbf{Proof of Theorem \ref{th:pv_weak}.}\qquad \medskip \newline
\textbf{Claim (a).}\qquad Recall $h$ $\equiv $ $(\gamma _{0},b)$ $\in $ $%
\mathfrak{H}$ $\equiv $ $\{h$ $:$ $\gamma _{0}$ $\in $ $\Gamma ^{\ast }$,
and $||b||$ $<$ $\infty $, with $\beta _{0}$ $=$ $0\}$. where $\Gamma ^{\ast
}$ is the true parameter space in (\ref{gamma*}).

Let $F_{\gamma }$ be the distribution function of $W_{t}$ $=$ $%
[y_{t},x_{t}^{\prime }]^{\prime }$ under some $\gamma $ $\in $ $\Gamma
^{\ast }$. Let $P_{\gamma }$ denote probability under $F_{\gamma }$. Recall $%
\mathcal{F}_{\infty }(c)$ $\equiv $ $P(\mathcal{T}(\lambda )$ $\leq $ $c)$,
where $\{\mathcal{T(}\lambda \mathcal{)}$ $:$ $\lambda $ $\in $ $\Lambda
\mathcal{\}}$ is the asymptotic null chi-squared process under strong
identification, and $\mathcal{F}_{\lambda ,h}(c)$ $\equiv $ $P(\mathcal{T}%
_{\psi }(\lambda ,h)$ $\leq $ $c)$ where $\{\mathcal{T}_{\psi }(\lambda ,h)$
$:$ $\lambda $ $\in $ $\Lambda \}$ is the asymptotic null process under weak
identification. Write the finite sample p-values $p_{n}^{\infty }(\lambda )$
$\equiv $ $1$ $-$ $\mathcal{F}_{\infty }(\mathcal{T}_{n}(\lambda ))$ $=$ $%
\mathcal{\bar{F}}_{\infty }(\mathcal{T}_{n}(\lambda ))$ and $p_{n}(\lambda
,h)$ $\equiv $ $1$ $-\mathcal{F}_{\lambda ,h}(\mathcal{T}_{n}(\lambda ))$ $=$
$\mathcal{\bar{F}}_{\lambda ,h}(\mathcal{T}_{n}(\lambda ))$. Recall that $%
\mathcal{F}_{\lambda ,\tilde{h}}(\cdot )$ is continuous by Assumption \ref%
{assum:pvs}.\medskip

\textbf{Step 1 (LF).} The asymptotic size $AsySz$ $\equiv $ $\sup_{\lambda
\in \Lambda }AsySz(\lambda )$ can be written as:%
\begin{eqnarray*}
AsySz &=&\sup_{\lambda \in \Lambda }\limsup_{n\rightarrow \infty
}\sup_{\gamma \in \Gamma ^{\ast }}P_{\gamma }\left( \max \left\{ \sup_{h\in
\mathfrak{H}}\left\{ p_{n}(\lambda ,h)\right\} ,p_{n}^{\infty }(\lambda
)\right\} <\alpha |H_{0}\right) \\
&=&\sup_{\lambda \in \Lambda }\limsup_{n\rightarrow \infty }\sup_{\gamma \in
\Gamma ^{\ast }}P_{\gamma }\left( \max \left\{ \sup_{h\in \mathfrak{H}%
}\left\{ \mathcal{\bar{F}}_{\lambda ,h}(\mathcal{T}_{n}(\lambda ))\right\} ,%
\mathcal{\bar{F}}_{\infty }(\mathcal{T}_{n}(\lambda ))\right\} <\alpha
|H_{0}\right) \equiv \mathfrak{A},
\end{eqnarray*}%
say. By Theorem \ref{th:CM_weak}.a, $\{\mathcal{T}_{n}(\lambda )$ $:$ $%
\Lambda \}$ $\Rightarrow ^{\ast }$ $\{\mathcal{T}_{\psi }(\lambda ,h)$ $:$ $%
\Lambda \}$ under $\mathcal{C}(i,b)$ with $||b||$ $<$ $\infty $. Weak
convergence implies convergence in finite dimensional distributions. By the
definition of distribution convergence, and the mapping theorem, weak
convergence therefore yields:
\begin{eqnarray}
\mathfrak{A} &=&\sup_{\lambda \in \Lambda }\sup_{\tilde{h}\in \mathfrak{H}%
}P\left( \max \left\{ \sup_{h\in \mathfrak{H}}\left\{ \mathcal{\bar{F}}%
_{\lambda ,h}\left( \mathcal{T}_{\psi }(\lambda ,\tilde{h})\right) \right\} ,%
\mathcal{\bar{F}}_{\infty }\left( \mathcal{T}_{\psi }(\lambda ,\tilde{h}%
)\right) \right\} <\alpha \right)  \label{A_lamda1} \\
&\leq &\sup_{\lambda \in \Lambda }\sup_{\tilde{h}\in \mathfrak{H}}P\left(
\sup_{h\in \mathfrak{H}}\left\{ \mathcal{\bar{F}}_{\lambda ,h}\left(
\mathcal{T}_{\psi }(\lambda ,\tilde{h})\right) \right\} <\alpha \right)
\notag \\
&\leq &\sup_{\lambda \in \Lambda }\sup_{\tilde{h}\in \mathfrak{H}}P\left(
\mathcal{\bar{F}}_{\lambda ,\tilde{h}}\left( \mathcal{T}_{\psi }(\lambda ,%
\tilde{h})\right) <\alpha \right) =\alpha .  \notag
\end{eqnarray}%
In the first inequality notice $\sup_{h\in \mathfrak{H}}$ operates only on
the distribution function $\mathcal{\bar{F}}_{\lambda ,h}$, whereas $\sup_{%
\tilde{h}\in \mathfrak{H}}$ operates on the limit process $\mathcal{T}_{\psi
}(\lambda ,\tilde{h})$. The inequalities hold since $\sup_{h\in \mathfrak{H}%
}\{\mathcal{\bar{F}}_{\lambda ,h}(\mathcal{T}_{\psi }(\lambda ,\tilde{h}))\}$
$\leq $ \linebreak $\max \{\sup_{h\in \mathfrak{H}}\{\mathcal{\bar{F}}%
_{\lambda ,h}(\mathcal{T}_{\psi }(\lambda ,\tilde{h}))\},$ $\mathcal{\bar{F}}%
_{\infty }(\mathcal{T}_{\psi }(\lambda ,\tilde{h}))\}$ and $\mathcal{\bar{F}}%
_{\lambda ,\tilde{h}}(\mathcal{T}_{\psi }(\lambda ,\tilde{h}))$ $\leq $ $%
\sup_{h\in \mathfrak{H}}\{\mathcal{\bar{F}}_{\lambda ,h}(\mathcal{T}_{\psi
}(\lambda ,\tilde{h}))\}$. The last equality in (\ref{A_lamda1}) applies
since $\mathcal{T}_{\psi }(\lambda ,\tilde{h})$ is distributed $\mathcal{F}%
_{\lambda ,\tilde{h}}$ which is continuous: hence $P(\mathcal{\bar{F}}%
_{\lambda ,\tilde{h}}(\mathcal{T}_{\psi }(\lambda ,\tilde{h}))$ $<$ $\alpha
) $ $=$ $\alpha $ for any $\tilde{h}$ and $\lambda $.

Under $\mathcal{C}(ii,\omega _{0})$, $\{\mathcal{T}_{n}(\lambda )$ $: $ $%
\Lambda \}$ $\Rightarrow ^{\ast }$ $\{\mathcal{T}(\lambda )$ $:$ $\Lambda \}$
by Theorem \ref{th:CM_weak}.b. Since $\mathcal{T}(\lambda )$ is distributed $%
\mathcal{F}_{\infty }$, which is continuous:
\begin{equation*}
\mathfrak{A}=\sup_{\lambda \in \Lambda }\sup_{\gamma \in \Gamma ^{\ast
}}P\left( \max \left\{ \sup_{h\in \mathfrak{H}}\left\{ \mathcal{\bar{F}}%
_{\lambda ,h}\left( \mathcal{T}(\lambda )\right) \right\} ,\mathcal{\bar{F}}%
_{\infty }\left( \mathcal{T}(\lambda )\right) \right\} <\alpha \right) \leq
\sup_{\lambda \in \Lambda }\sup_{\gamma \in \Gamma ^{\ast }}P\left( \mathcal{%
\bar{F}}_{\infty }\left( \mathcal{T}(\lambda )\right) <\alpha \right)
=\alpha .
\end{equation*}

\textbf{Step 2 (ICS-1).}\qquad\ By Theorem \ref{th:ICS_stat}.a, $\mathcal{A}%
_{n}$ $=$ $O_{p}(1)$\ under $\mathcal{C}(i,b)$ with $||b||$ $<$ $\infty $.
Further $\kappa _{n}$ $\rightarrow $ $\infty $ and $\kappa _{n}$ $= $ $%
o(n^{1/2})$. Hence
\begin{equation*}
p_{n}^{(ICS-1)}(\lambda )\equiv \left\{ p_{n}^{(LF)}(\lambda )\text{ if }%
\mathcal{A}_{n}\leq \kappa _{n},\text{ }p_{n}^{\infty }(\lambda )\text{ if }%
\mathcal{A}_{n}>\kappa _{n}\right\} =p_{n}^{(LF)}(\lambda )
\end{equation*}%
asymptotically with probability approaching one.

Under $\mathcal{C}(ii,\omega _{0})$, $\mathcal{A}_{n}\overset{p}{\rightarrow
}$ $\infty $ by Theorem \ref{th:ICS_stat}.a. If $\kappa _{n}^{-1}\mathcal{A}%
_{n}\overset{p}{\rightarrow }$ $[0,\infty )$ then again $p_{n}^{(ICS-1)}(%
\lambda )=$ $p_{n}^{(LF)}(\lambda )$ asymptotically with probability
approaching one. If $\kappa _{n}^{-1}\mathcal{A}_{n}\overset{p}{\rightarrow }
$ $\infty $, for example when $\beta _{0}=$ $0$ (see Theorem \ref%
{th:ICS_stat}.b,) then $p_{n}^{(ICS-1)}(\lambda )=$ $p_{n}^{\infty }(\lambda
)$ asymptotically with probability approaching one.

In each case, a p-value is chosen that leads to correct asymptotic level $%
AsySz$ $\leq $ $\alpha $\ in view of Step 1.\medskip \newline
\textbf{Claim (b).}\qquad Let $H_{0}$ be false. By Theorem \ref{th:CM_H1} $%
\mathcal{T}_{n}(\lambda )$ $\overset{p}{\rightarrow }$ $\infty $ $\forall
\lambda $ $\in $ $\Lambda /S$ where $S$ $\subset $ $\Lambda $ has Lebesgue
measure zero. Theorem \ref{th:CM_weak}.a states $\sup_{\lambda \in \Lambda
}\{\mathcal{T}_{\psi }(\lambda ,b)\}$ $<$ $\infty $ $a.s$., hence the
distribution of $\mathcal{T}_{\psi }(\lambda ,b)$ has support $[0,\infty )$.
Therefore $\sup_{h\in \mathfrak{H}}\{p_{n}(\lambda ,h)\}$ $\overset{p}{%
\rightarrow }$ $0$ and $p_{n}^{\infty }(\lambda )$ $\overset{p}{\rightarrow }
$ $0$, hence by construction $p_{n}^{(LF)}(\lambda )$ $\overset{p}{%
\rightarrow }$ $0$ and by arguments under Step 2 above $p_{n}^{(ICS-1)}(%
\lambda )$ $\overset{p}{\rightarrow }$ $0$. $\mathcal{QED}.\bigskip $\newline
\textbf{Proof of Theorem \ref{th:boot_pv}.\medskip }\newline
\textbf{Claim (a).\medskip }

\textbf{Step 1.} Operate conditionally on the sample $\mathfrak{W}_{n}$ $%
\equiv $ $\{(y_{t},x_{t})\}_{t=1}^{n}$. In this step we prove the
bootstrapped test statistic converges weakly in probability to the Theorem %
\ref{th:CM_weak}.a null limit process under weak identification:%
\begin{equation}
\left\{ \mathcal{\hat{T}}_{\psi ,n}^{\ast }(\lambda ,h):\lambda \in \Lambda
\right\} \mathcal{\Rightarrow }^{p}\left\{ \left( \frac{\mathfrak{T}_{\psi
}(\pi ^{\ast }(b),\lambda ,b)}{\bar{v}(\pi ^{\ast }(b),\lambda ,b)}\right)
^{2}:\lambda \in \Lambda \right\} =\left\{ \mathcal{T}_{\psi }(\lambda
,h):\lambda \in \Lambda \right\} .  \label{Cn_boot_weak}
\end{equation}%
We then prove the claim in Step 2.\medskip

\textbf{Step 1.1} Recall $\sigma _{0}\widehat{\mathcal{G}}_{\psi ,n}^{\ast
}(\pi )$ $\equiv $ $1/\sqrt{n}\sum_{t=1}^{n}\sigma _{0}z_{t}\widehat{%
\mathcal{H}}_{\psi ,n}^{-1/2}(\pi )d_{\psi ,t}(\pi )$ where $z_{t}$\ is iid $%
N(0,1)$. We will prove $\{\sigma _{0}\widehat{\mathcal{G}}_{\psi ,n}^{\ast
}(\pi )$ $:$ $\pi $ $\in $ $\Pi \}$ $\mathcal{\Rightarrow }^{p}$ $\{\mathcal{%
H}_{\psi }^{-1/2}(\pi )\mathcal{G}_{\psi }(\pi )$ $:$ $\pi $ $\in $ $\Pi \}$%
, where $\mathcal{G}_{\psi }(\pi )$ $=$ $\mathcal{G}_{\psi }(\psi _{0},\pi )$%
, and $\mathcal{G}_{\psi }(\theta )$ $\equiv $ $\mathcal{G}_{\psi }(\psi
,\pi )$ is the Lemma \ref{lm:G_uclt_weak} case $\mathcal{C}(i,b)$ limit
process. We need to prove convergence in finite dimensional distributions,
and demonstrate stochastic equicontinuity.\footnote{%
See Theorem 3.1 in \cite{GineZinn1990}, and see \cite{Dudley1978} and
\citet[Theorem
1.1.3]{GineZinn1986}.}

In order to establish convergence in finite dimensional distributions, we
use Hansen's (\citeyear{Hansen1996}, proof of Theorem 2) argument. Denote $%
E_{\mathfrak{W}_{n}}[\cdot ]$ $=$ $E[\cdot |\mathfrak{W}_{n}]$. By
Gaussianicity of $z_{t}$, $\sigma _{0}\widehat{\mathcal{G}}_{\psi ,n}^{\ast
}(\pi )$ is (conditionally) normally distributed with mean zero and
covariance kernel:%
\begin{equation*}
E_{\mathfrak{W}_{n}}\left[ \sigma _{0}^{2}\widehat{\mathcal{G}}_{\psi
,n}^{\ast }(\pi )\widehat{\mathcal{G}}_{\psi ,n}^{\ast }(\tilde{\pi}%
)^{\prime }\right] =\sigma _{0}^{2}\widehat{\mathcal{H}}_{\psi
,n}^{-1/2}(\pi )\frac{1}{n}\sum_{t=1}^{n}d_{\psi ,t}(\pi )d_{\psi ,t}(\tilde{%
\pi})^{\prime }\widehat{\mathcal{H}}_{\psi ,n}^{-1/2}(\tilde{\pi})=\sigma
_{0}^{2}\widehat{\mathcal{H}}_{\psi ,n}^{-1/2}(\pi )\widehat{\mathcal{H}}%
_{\psi ,n}(\pi ,\tilde{\pi})\widehat{\mathcal{H}}_{\psi ,n}^{-1/2}(\tilde{\pi%
}),
\end{equation*}%
where $\widehat{\mathcal{H}}_{\psi ,n}(\pi ,\tilde{\pi})$\ is implicitly
defined. Let $\mathcal{W}$ be the set of asymptotic sample $%
\{(y_{t},x_{t})\}_{t=1}^{\infty }$ such that
\begin{equation}
\sup_{\pi ,\tilde{\pi}\in \Pi \times \Pi }\left\Vert E_{\mathfrak{W}_{n}}%
\left[ \widehat{\mathcal{G}}_{\psi ,n}^{\ast }(\pi )\widehat{\mathcal{G}}%
_{\psi ,n}^{\ast }(\tilde{\pi})^{\prime }\right] -\mathcal{H}_{\psi
}^{-1/2}(\pi )\mathcal{H}_{\psi }(\pi ,\tilde{\pi})\mathcal{H}_{\psi
}^{-1/2}(\pi )\right\Vert \overset{p}{\rightarrow }0.  \label{supEGG}
\end{equation}%
By Lemma \ref{lm:H_ulln_weak} $\sup_{\pi \in \Pi }||\widehat{\mathcal{H}}%
_{\psi ,n}(\pi )$ $-$ $\mathcal{H}_{\psi }(\pi )||$ $\overset{p}{\rightarrow
}$ $0$, and by Assumption \ref{assum:dgp}.d(iii) $\mathcal{H}_{\psi }(\pi )$
is positive definite uniformly on $\Pi $. By the same argument used to prove
Lemma \ref{lm:H_ulln_weak}, $\sup_{\pi ,\tilde{\pi}\in \Pi }||\widehat{%
\mathcal{H}}_{\psi ,n}(\pi ,\tilde{\pi})$ $-$ $\mathcal{H}_{\psi }(\pi ,%
\tilde{\pi})||$ $\overset{p}{\rightarrow }$ $0$ where $\mathcal{H}_{\psi
}(\pi ,\tilde{\pi})$ $\equiv $ $E[d_{\psi ,t}(\pi )d_{\psi ,t}(\tilde{\pi}%
)^{\prime }]$. This proves $P(\mathfrak{W}_{n}$ $\in $ $\mathcal{W})$ $=$ $1$%
. Therefore $\sigma _{0}\widehat{\mathcal{G}}_{\psi ,n}^{\ast }(\pi )$
converges in finite dimensional distributions to a zero mean Gaussian law
with covariance kernel $\sigma _{0}^{2}\mathcal{H}_{\psi }^{-1/2}(\pi )%
\mathcal{H}_{\psi }(\pi ,\tilde{\pi})\mathcal{H}_{\psi }^{-1/2}(\pi )$.

$\mathcal{H}_{\psi }^{-1/2}(\pi )\mathcal{G}_{\psi ,n}(\pi )$ has the same
limit under the null and case $\mathcal{C}(i,b)$ with $||b||$ $<$ $\infty $.
This follows by Lemma \ref{lm:G_uclt_weak}, the $\mathcal{G}_{\psi ,n}(\psi
_{0,n},\pi )$ identify (\ref{G_psi_n}), and the following moment under the
null and Assumption \ref{assum:dgp}.a(i):
\begin{equation*}
E\left[ \mathcal{H}_{\psi }^{-1/2}(\pi )\frac{1}{\sqrt{n}}%
\sum_{t=1}^{n}\epsilon _{t}d_{\psi ,t}(\pi )\frac{1}{\sqrt{n}}%
\sum_{t=1}^{n}\epsilon _{t}d_{\psi ,t}(\tilde{\pi})^{\prime }\mathcal{H}%
_{\psi }^{-1/2}(\tilde{\pi})\right] =\sigma _{0}^{2}\mathcal{H}_{\psi
}^{-1/2}(\pi )\mathcal{H}_{\psi }(\pi ,\tilde{\pi})\mathcal{H}_{\psi
}^{-1/2}(\tilde{\pi}).
\end{equation*}%
Now apply Lemma \ref{lm:G_uclt_weak} to yield that $\mathcal{H}_{\psi
}^{-1/2}(\pi )\mathcal{G}_{\psi ,n}(\pi )$ converges in finite dimensional
distributions to $\mathcal{H}_{\psi }^{-1/2}(\pi )\mathcal{G}_{\psi }(\pi )$%
, a zero mean Gaussian law with kernel $\sigma _{0}^{2}\mathcal{H}_{\psi
}^{-1/2}(\pi )\mathcal{H}_{\psi }(\pi ,\tilde{\pi})\mathcal{H}_{\psi
}^{-1/2}(\pi )$. Since Gaussian processes are fully characterized by their
mean and covariance functions, $\sigma _{0}\widehat{\mathcal{G}}_{\psi
,n}^{\ast }(\pi )$ therefore converges in finite dimensional distributions
to $\mathcal{H}_{\psi }^{-1/2}(\pi )\mathcal{G}_{\psi }(\pi )$.

Next we establish stochastic equicontinuity. Let $r$ $\in $ $\mathbb{R}%
^{k_{x}+k_{\beta }}$, $r^{\prime }r$ $=$ $1$, be arbitrary. By the mean
value theorem, for some $\mathring{\pi}$ $\in $ $\Pi $, $||\mathring{\pi}$ $%
- $ $\pi ||$ $\leq $ $||\tilde{\pi}$ $-$ $\pi ||$:%
\begin{equation*}
r^{\prime }\widehat{\mathcal{H}}_{\psi ,n}^{-1/2}(\pi )d_{\psi ,t}(\pi
)-r^{\prime }\widehat{\mathcal{H}}_{\psi ,n}^{-1/2}(\tilde{\pi})d_{\psi ,t}(%
\tilde{\pi})=\left( \left[ r^{\prime }\frac{\partial }{\partial \pi _{i}}%
\widehat{\mathcal{H}}_{\psi ,n}^{-1/2}(\mathring{\pi})d_{\psi ,t}(\mathring{%
\pi})+r^{\prime }\widehat{\mathcal{H}}_{\psi ,n}^{-1/2}(\mathring{\pi})\frac{%
\partial }{\partial \pi _{i}}d_{\psi ,t}(\mathring{\pi})\right]
_{i=1}^{k_{\pi }}\right) ^{\prime }\left( \pi -\tilde{\pi}\right) .
\end{equation*}%
The derivatives are%
\begin{eqnarray*}
&&\frac{\partial }{\partial \pi }d_{\psi ,t}(\pi )=\left[
\begin{array}{c}
\frac{\partial }{\partial \pi }g(x_{t},\pi ) \\
0_{k_{x}\times k_{\pi }}%
\end{array}%
\right] \text{ \ and \ }\frac{\partial }{\partial \pi _{i}}\widehat{\mathcal{%
H}}_{\psi ,n}^{-1/2}(\pi )=-\widehat{\mathcal{H}}_{\psi ,n}^{-1/2}(\pi )%
\frac{\partial }{\partial \pi _{i}}\widehat{\mathcal{H}}_{\psi ,n}(\pi )%
\widehat{\mathcal{H}}_{\psi ,n}^{-1/2}(\pi ) \\
&&\frac{\partial }{\partial \pi _{i}}\widehat{\mathcal{H}}_{\psi ,n}(\pi )=%
\frac{1}{n}\sum_{t=1}^{n}\left[
\begin{array}{cc}
\frac{\partial }{\partial \pi _{i}}g(x_{t},\pi )g(x_{t},\pi )^{\prime
}+g(x_{t},\pi )\frac{\partial }{\partial \pi _{i}}g(x_{t},\pi )^{\prime } &
\frac{\partial }{\partial \pi _{i}}g(x_{t},\pi )x_{t}^{\prime } \\
x_{t}\frac{\partial }{\partial \pi _{i}}g(x_{t},\pi )^{\prime } &
0_{k_{x}\times k_{x}}%
\end{array}%
\right] .
\end{eqnarray*}%
Invoke Chebyshev's inequality, and the fact that $z_{t}$ is iid, independent
of $\mathfrak{W}_{n}$, and not a function of $\pi $, to yield
\begin{eqnarray*}
\mathcal{P}_{n}(\eta ) &\equiv &P\left( \sup_{\pi ,\tilde{\pi}\in \Pi :||\pi
-\tilde{\pi}||\leq \delta }\left\vert \frac{1}{\sqrt{n}}\sum_{t=1}^{n}z_{t}%
\left\{ r^{\prime }\widehat{\mathcal{H}}_{\psi ,n}^{-1/2}(\pi )d_{\psi
,t}(\pi )-r^{\prime }\widehat{\mathcal{H}}_{\psi ,n}^{-1/2}(\tilde{\pi}%
)d_{\psi ,t}(\tilde{\pi})\right\} \right\vert >\eta |\mathfrak{W}_{n}\right)
\\
&\leq &\frac{1}{\eta ^{2}}E\left[ \sup_{\pi ,\tilde{\pi}\in \Pi :||\pi -%
\tilde{\pi}||\leq \delta }\left( \frac{1}{\sqrt{n}}\sum_{t=1}^{n}z_{t}\left%
\{ r^{\prime }\widehat{\mathcal{H}}_{\psi ,n}^{-1/2}(\pi )d_{\psi ,t}(\pi
)-r^{\prime }\widehat{\mathcal{H}}_{\psi ,n}^{-1/2}(\tilde{\pi})d_{\psi ,t}(%
\tilde{\pi})\right\} \right) ^{2}|\mathfrak{W}_{n}\right] \\
&=&\frac{1}{\eta ^{2}}\frac{1}{n}\sum_{t=1}^{n}\sup_{\pi ,\tilde{\pi}\in \Pi
:||\pi -\tilde{\pi}||\leq \delta }\left\{ r^{\prime }\widehat{\mathcal{H}}%
_{\psi ,n}^{-1/2}(\pi )d_{\psi ,t}(\pi )-r^{\prime }\widehat{\mathcal{H}}%
_{\psi ,n}^{-1/2}(\tilde{\pi})d_{\psi ,t}(\tilde{\pi})\right\} ^{2} \\
&\leq &\frac{\delta ^{2}}{\eta ^{2}}\frac{1}{n}\sum_{t=1}^{n}\sup_{\pi \in
\Pi }\left\vert \left[ r^{\prime }\frac{\partial }{\partial \pi _{i}}%
\widehat{\mathcal{H}}_{\psi ,n}^{-1/2}(\pi )d_{\psi ,t}(\pi )+r^{\prime }%
\widehat{\mathcal{H}}_{\psi ,n}^{-1/2}(\pi )\frac{\partial }{\partial \pi
_{i}}d_{\psi ,t}(\pi )\right] _{i=1}^{k_{\pi }}\right\vert ^{2}=\frac{\delta
^{2}}{\eta ^{2}}\mathcal{C}_{n},
\end{eqnarray*}%
say. We prove below that $\mathcal{C}_{n}$ $\overset{p}{\rightarrow }$ $%
\mathcal{C}$ a finite non-negative constant. We can therefore choose any $%
\delta $ $>$ $0$ if $\mathcal{C}$ $=$ $0$, and otherwise $0$ $<$ $\delta $ $%
\leq $ $[\epsilon \eta ^{2}/\mathcal{C}]^{1/2}$, such that for each $%
(\epsilon ,\eta )$ $>$ $0$ there exists $\delta $ $>$ $0$ yielding $%
\lim_{n\rightarrow \infty }\mathcal{P}_{n}(\eta )$ $<$ $\epsilon $
asymptotically with probability approaching one with respect to the sample
draw $\mathfrak{W}_{n}$. This establishes stochastic equicontinuity.

We now prove $\mathcal{C}_{n}$ $\overset{p}{\rightarrow }$ $\mathcal{C}$ $%
\in $ $[0,\infty )$.\ Since $\widehat{\mathcal{H}}_{\psi ,n}(\pi )$ $\equiv $
$1/n\sum_{t=1}^{n}d_{\psi ,t}(\pi )d_{\psi ,t}(\pi )^{\prime }$, note that:%
\begin{eqnarray*}
&&\frac{1}{n}\sum_{t=1}^{n}\sup_{\pi \in \Pi }\left\vert \left[ r^{\prime }%
\frac{\partial }{\partial \pi _{i}}\widehat{\mathcal{H}}_{\psi
,n}^{-1/2}(\pi )d_{\psi ,t}(\pi )+r^{\prime }\widehat{\mathcal{H}}_{\psi
,n}^{-1/2}(\pi )\frac{\partial }{\partial \pi _{i}}d_{\psi ,t}(\pi )\right]
_{i=1}^{k_{\pi }}\right\vert ^{2} \\
&&\text{ \ \ \ \ \ \ \ }=\sum_{i,j=1}^{k_{\pi }}\sup_{\pi \in \Pi }r^{\prime
}\frac{\partial }{\partial \pi _{i}}\widehat{\mathcal{H}}_{\psi
,n}^{-1/2}(\pi )\widehat{\mathcal{H}}_{\psi ,n}(\pi )\frac{\partial }{%
\partial \pi _{j}}\widehat{\mathcal{H}}_{\psi ,n}^{-1/2}(\pi )r \\
&&\text{\ \ \ \ \ \ \ \ \ \ \ \ \ \ \ \ \ \ \ \ \ }+\sum_{i,j=1}^{k_{\pi
}}\sup_{\pi \in \Pi }r^{\prime }\widehat{\mathcal{H}}_{\psi ,n}^{-1/2}(\pi )%
\frac{1}{n}\sum_{t=1}^{n}\frac{\partial }{\partial \pi _{i}}d_{\psi ,t}(\pi )%
\frac{\partial }{\partial \pi _{j}}d_{\psi ,t}(\pi )^{\prime }\widehat{%
\mathcal{H}}_{\psi ,n}^{-1/2}(\pi )r \\
&&\text{ \ \ \ \ \ \ \ \ \ \ \ \ \ \ \ \ \ \ \ \ }+2\sum_{i,j=1}^{k_{\pi
}}\sup_{\pi \in \Pi }r^{\prime }\frac{\partial }{\partial \pi _{i}}\widehat{%
\mathcal{H}}_{\psi ,n}^{-1/2}(\pi )\frac{1}{n}\sum_{t=1}^{n}d_{\psi ,t}(\pi )%
\frac{\partial }{\partial \pi _{i}}d_{\psi ,t}(\pi )^{\prime }\widehat{%
\mathcal{H}}_{\psi ,n}^{-1/2}(\pi )r.
\end{eqnarray*}%
The argument used to prove Lemmas \ref{lm:H_ulln_weak} and \ref%
{lm:H_ulln_strong} extend to each component of $(\partial /\partial \pi _{i})%
\widehat{\mathcal{H}}_{\psi ,n}^{-1/2}(\pi )$, $1/n\sum_{t=1}^{n}d_{\psi
,t}(\pi )(\partial /\partial \pi _{i})d_{\psi ,t}(\pi )^{\prime }$ and $%
1/n\sum_{t=1}^{n}(\partial /\partial \pi _{i})d_{\psi ,t}(\pi )(\partial
/\partial \pi _{i})d_{\psi ,t}(\pi )^{\prime }$, in view of the Assumption %
\ref{assum:dgp}.b,c mixing and moment bounds. Thus, each summand has a
uniformly bounded uniform probability limit under any case $\mathcal{C}(i,b)$
or $\mathcal{C}(ii,\omega _{0})$. Hence, by Slutsky's theorem:%
\begin{eqnarray*}
&&\mathcal{C}_{n}\overset{p}{\rightarrow }\sum_{i,j=1}^{k_{\pi }}\sup_{\pi
\in \Pi }\left\{ r^{\prime }\frac{\partial }{\partial \pi _{i}}\mathcal{H}%
_{\psi }^{-1/2}(\pi )\mathcal{H}_{\psi }(\pi )\frac{\partial }{\partial \pi
_{j}}\mathcal{H}_{\psi }^{-1/2}(\pi )r\right\} \\
&&\text{ \ \ \ \ \ \ \ \ \ \ }+\sum_{i,j=1}^{k_{\pi }}\sup_{\pi \in \Pi
}\left\{ r^{\prime }\mathcal{H}_{\psi }^{-1/2}(\pi )E\left[ \frac{\partial }{%
\partial \pi _{i}}d_{\psi ,t}(\pi )\frac{\partial }{\partial \pi _{j}}%
d_{\psi ,t}(\pi )^{\prime }\right] \widehat{\mathcal{H}}_{\psi
,n}^{-1/2}(\pi )r\right\} \\
&&\text{ \ \ \ \ \ \ \ \ \ }+2\sum_{i,j=1}^{k_{\pi }}\sup_{\pi \in \Pi
}\left\{ r^{\prime }\frac{\partial }{\partial \pi _{i}}\mathcal{H}_{\psi
}^{-1/2}(\pi )E\left[ d_{\psi ,t}(\pi )\frac{\partial }{\partial \pi _{i}}%
d_{\psi ,t}(\pi )^{\prime }\right] \mathcal{H}_{\psi }^{-1/2}(\pi )r\right\}
\equiv \mathcal{C}<\infty .
\end{eqnarray*}%
Non-negativity $\mathcal{C}$ $\geq $ $0$ is trivial in view of the quadratic
form of $\mathcal{C}_{n}$.\medskip

\textbf{Step 1.2} Next, we prove the bootstrapped $\hat{\pi}_{n}^{\ast }(\pi
_{0},b)$ $\equiv $ $\argmin_{\pi \in \Pi }\{\hat{\xi}_{\psi ,n}^{\ast }(\pi
,\pi _{0},b)\}$ satisfies:
\begin{equation}
\hat{\pi}_{n}^{\ast }(\pi _{0},b)\mathcal{\Rightarrow }^{p}\argmin_{\pi \in
\Pi }\left\{ -\frac{1}{2}\left( \mathcal{G}_{\psi }(\pi )+\mathcal{D}_{\psi
}(\pi ,\pi _{0})b\right) ^{\prime }\mathcal{H}_{\psi }^{-1}(\pi )\left(
\mathcal{G}_{\psi }(\pi )+\mathcal{D}_{\psi }(\pi ,\pi _{0})b\right)
\right\} \equiv \pi ^{\ast }(b).  \label{L0}
\end{equation}%
$\widehat{\mathcal{H}}_{\psi ,n}(\pi )$ and $\mathcal{\hat{D}}_{\psi ,n}(\pi
,\pi _{0})$ have uniform probability limits $\mathcal{H}_{\psi }(\pi )$ and $%
\mathcal{D}_{\psi }(\pi ,\pi _{0})$ by Lemma \ref{lm:JV}. The Step 1.1
result of weak convergence in probability, the mapping theorem, and $\hat{%
\sigma}_{n}$ $\overset{p}{\rightarrow }$ $\sigma _{0}$ (see Remark \ref%
{rm:sig_hat}), together yield :%
\begin{eqnarray*}
&&\left\{ \hat{\xi}_{\psi ,n}^{\ast }(\pi ,\pi _{0},b):\pi \in \Pi \right\}
\\
&&\text{ \ \ \ }=\left\{ -\frac{1}{2}\left( \hat{\sigma}_{n}\widehat{%
\mathcal{G}}_{\psi ,n}^{\ast }(\pi )+\widehat{\mathcal{H}}_{\psi
,n}^{-1/2}(\pi )\mathcal{\hat{D}}_{\psi ,n}(\pi ,\pi _{0})\times b\right)
^{\prime }\left( \hat{\sigma}_{n}\widehat{\mathcal{G}}_{\psi ,n}^{\ast }(\pi
)+\widehat{\mathcal{H}}_{\psi ,n}^{-1/2}(\pi )\mathcal{\hat{D}}_{\psi
,n}(\pi ,\pi _{0})\times b\right) :\pi \in \Pi \right\} \\
&&\text{ \ \ \ }\Rightarrow ^{p}-\frac{1}{2}\left\{ \mathcal{G}_{\psi }(\pi
)+\mathcal{D}_{\psi }(\pi ,\pi _{0})b\right\} ^{\prime }\mathcal{H}_{\psi
}^{-1}(\pi )\left\{ \mathcal{G}_{\psi }(\pi )+\mathcal{D}_{\psi }(\pi ,\pi
_{0})b:\pi \in \Pi \right\} .
\end{eqnarray*}
Apply the mapping theorem again, and Assumption \ref{assum:pi}.a, to yield (%
\ref{L0}).\medskip

\textbf{Step 1.3} Define $\mathcal{K}_{n,t}(\pi ,\lambda )$ $\equiv $ $%
F(\lambda ^{\prime }\mathcal{W}(x_{t}))$ $-$ $\mathfrak{\hat{b}}_{\psi
,n}(\pi ,\lambda )^{\prime }\widehat{\mathcal{H}}_{\psi ,n}^{-1}(\pi
)d_{\psi ,t}(\pi )$ and recall $\mathcal{K}_{\psi ,t}(\pi ,\lambda )$ $%
\equiv $ $F(\lambda ^{\prime }\mathcal{W}(x_{t}))$ $-$ $\mathfrak{b}_{\psi
}(\pi ,\lambda )^{\prime }\mathcal{H}_{\psi }^{-1}(\pi )d_{\psi ,t}(\pi )$.
We will show $\mathfrak{\hat{T}}_{\psi ,n}^{\ast }(\pi ,\lambda ,\pi _{0},b)$
defined in (\ref{tao_boot}) satisfies:
\begin{equation*}
\left\{ \mathfrak{\hat{T}}_{\psi ,n}^{\ast }(\pi ,\lambda ,\pi _{0},b):\Pi
,\Lambda \right\} \mathcal{\Rightarrow }^{p}\text{\ }\left\{ \mathfrak{T}%
_{\psi }(\pi ,\lambda ,b):\Pi ,\Lambda \right\} ,
\end{equation*}%
where, as in (\ref{T(pi,lam,b)}),%
\begin{eqnarray*}
\mathfrak{T}_{\psi }(\pi ,\lambda ,b) &\equiv &\mathfrak{Z}_{\psi }(\pi
,\lambda )+\mathfrak{b}_{\psi }(\pi ,\lambda )^{\prime }\left( \mathcal{H}%
_{\psi }^{-1}(\pi )\mathcal{D}_{\psi }(\pi )b+\left[ b,0_{k_{\beta
}}^{\prime }\right] ^{\prime }\right) \\
&&+\mathfrak{b}_{\psi }(\pi ,\lambda )^{\prime }\mathcal{H}_{\psi }^{-1}(\pi
)E\left[ d_{\psi ,t}(\pi )\left\{ g(x_{t},\pi _{0})-g(x_{t},\pi )\right\}
^{\prime }\right] b \\
&&+E\left[ \mathcal{K}_{\psi ,t}(\pi ,\lambda )\left\{ g(x_{t},\pi
_{0})-g(x_{t},\pi )\right\} ^{\prime }\right] b,
\end{eqnarray*}%
and $\mathfrak{Z}_{\psi }(\pi ,\lambda )$ is the Lemma \ref{lm:CM_weak} zero
mean Gaussian limit process of $1/\sqrt{n}\sum_{t=1}^{n}\epsilon _{t}%
\mathcal{K}_{\psi ,t}(\pi ,\lambda )$.

Observe that $\widehat{\mathcal{H}}_{\psi ,n}(\pi )$, $\mathcal{\hat{D}}%
_{\psi ,n}(\pi ,\pi _{0})$, $\mathfrak{\hat{b}}_{\psi ,n}(\hat{\psi}_{n},\pi
,\lambda )$, $1/n\sum_{t=1}^{n}d_{\psi ,t}(\pi )g(x_{t},\pi )^{\prime }$ and
$1/n\sum_{t=1}^{n}\mathcal{K}_{n,t}(\pi ,\lambda )g(x_{t},\pi )$ have
uniform probability limits $\mathfrak{b}_{\psi }(\pi ,\lambda )$, $\mathcal{H%
}_{\psi }(\pi )$, $\mathcal{D}_{\psi }(\pi )$, $E[d_{\psi ,t}(\pi
)g(x_{t},\pi )^{\prime }]$, and $E[\mathcal{K}_{\psi ,t}(\pi ,\lambda
)g(x_{t},\pi )]$ by applications of Lemmas \ref{lm:H_ulln_weak}, \ref{lm:bn}
and \ref{lm:ulln_xgF_eF}.

It therefore suffices to prove $\{\sigma _{0}\mathfrak{\hat{Z}}_{\psi
,n}^{\ast }(\pi ,\lambda )$ $:$ $\Pi ,\Lambda \}$ $\mathcal{\Rightarrow }%
^{p} $\ $\{\mathfrak{Z}_{\psi }(\pi ,\lambda )$ $:$ $\Pi ,\Lambda \}$ where $%
\mathfrak{\hat{Z}}_{\psi ,n}^{\ast }(\pi ,\lambda )$ $\equiv $ \linebreak $1/%
\sqrt{n}\sum_{t=1}^{n}\sigma _{0}z_{t}\mathcal{K}_{n,t}(\pi ,\lambda )$.
Note that $\mathfrak{\hat{Z}}_{\psi ,n}^{\ast }(\pi ,\lambda )$ is
(conditionally) normally distributed with zero mean and covariance kernel $%
\sigma _{0}^{2}1/n\sum_{t=1}^{n}\mathcal{K}_{n,t}(\pi ,\lambda )\mathcal{K}%
_{n,t}(\tilde{\pi},\tilde{\lambda})$. Let $\mathcal{\tilde{W}}$ be the set
of asymptotic sample $\{(y_{t},x_{t})\}_{t=1}^{\infty }$ such that
\begin{equation*}
\sup_{\pi ,\tilde{\pi}\in \Pi \times \Pi ,\lambda ,\tilde{\lambda}\in
\Lambda }\left\Vert \frac{1}{n}\sum_{t=1}^{n}\mathcal{K}_{n,t}(\pi ,\lambda )%
\mathcal{K}_{n,t}(\tilde{\pi},\tilde{\lambda})-E\left[ \mathcal{K}_{\psi
,t}(\pi ,\lambda )\mathcal{K}_{\psi ,t}(\tilde{\pi},\tilde{\lambda})\right]
\right\Vert \overset{p}{\rightarrow }0.
\end{equation*}%
By Lemmas \ref{lm:H_ulln_weak} and \ref{lm:bn}: $\sup_{\pi ,\tilde{\pi}\in
\Pi \times \Pi ,\lambda ,\tilde{\lambda}\in \Lambda }|1/n\sum_{t=1}^{n}\{%
\mathcal{K}_{n,t}(\pi ,\lambda )\mathcal{K}_{n,t}(\tilde{\pi},\tilde{\lambda}%
)$ $-$ $\mathcal{K}_{\psi ,t}(\pi ,\lambda )\mathcal{K}_{\psi ,t}(\tilde{\pi}%
,\tilde{\lambda})\}|$ $\overset{p}{\rightarrow }$ $0$ and by the same
arguments used to prove Lemma \ref{lm:bn}: $\sup_{\pi ,\tilde{\pi}\in \Pi
\times \Pi ,\lambda ,\tilde{\lambda}\in \Lambda }|1/n\sum_{t=1}^{n}\{%
\mathcal{K}_{\psi ,t}(\pi ,\lambda )\mathcal{K}_{\psi ,t}(\tilde{\pi},\tilde{%
\lambda})$ $-$ $E[\mathcal{K}_{\psi ,t}(\pi ,\lambda )\mathcal{K}_{\psi ,t}(%
\tilde{\pi},\tilde{\lambda})]\}|$ $\overset{p}{\rightarrow }$ $0$. This
proves $P(\mathfrak{W}_{n}$ $\in $ $\mathcal{\tilde{W}})$ $=$ $1$. Hence, $%
\mathfrak{\hat{Z}}_{\psi ,n}^{\ast }(\pi ,\lambda )$ converges in finite
dimensional distributions to a zero mean Gaussian law with kernel $\sigma
_{0}^{2}E[\mathcal{K}_{\psi ,t}(\pi ,\lambda )\mathcal{K}_{\psi ,t}(\tilde{%
\pi},\tilde{\lambda})]$.

By Lemma \ref{lm:CM_weak}.a, under the null $\{\mathfrak{Z}_{\psi ,n}(\pi
,\lambda ):\Pi ,\Lambda \}\Rightarrow ^{\ast }\{\mathfrak{Z}_{\psi }(\pi
,\lambda )$ $:$ $\Pi ,\Lambda \}$, a zero mean Gaussian process with
covariance kernel $\sigma _{0}^{2}E[\mathcal{K}_{\psi ,t}(\pi ,\lambda )%
\mathcal{K}_{\psi ,t}(\tilde{\pi},\tilde{\lambda})]$. Therefore, the finite
dimensional distributions of $\{\mathfrak{\hat{Z}}_{\psi ,n}^{\ast }(\pi
,\lambda ))$ $:$ $\Pi ,\Lambda \}$ converge to those of $\{\mathfrak{Z}%
_{\psi }(\pi ,\lambda )$ $:$ $\Pi ,\Lambda \}$ under the null.

It remains to prove stochastic equicontinuity for $\mathfrak{\hat{Z}}_{\psi
,n}^{\ast }(\pi ,\lambda )$. By construction, we need only show the sequence
of distributions of $1/\sqrt{n}\sum_{t=1}^{n}z_{t}F(\lambda ^{\prime }%
\mathcal{W}(x_{t}))$ and $1/\sqrt{n}\sum_{t=1}^{n}z_{t}$ $\mathfrak{\hat{b}}%
_{\psi ,n}(\pi ,\lambda )^{\prime }\widehat{\mathcal{H}}_{\psi ,n}^{-1}(\pi
)d_{\psi ,t}(\pi )$ are stochastically equicontinuous, and invoke
probability subadditivity. For $1/\sqrt{n}\sum_{t=1}^{n}z_{t}F(\lambda
^{\prime }\mathcal{W}(x_{t}))$, by the mean value theorem and Chebyshev's
inequality and the fact $\{z_{t}\}_{t=1}^{n}$ is independent of the sample $%
\mathfrak{W}_{n}$:
\begin{eqnarray*}
&&P\left( \sup_{\lambda ,\tilde{\lambda}\in \Lambda :||\lambda -\tilde{%
\lambda}||\leq \delta }\left\vert \frac{1}{\sqrt{n}}\sum_{t=1}^{n}z_{t}\left%
\{ F(\lambda ^{\prime }\mathcal{W}(x_{t}))-F(\tilde{\lambda}^{\prime }%
\mathcal{W}(x_{t}))\right\} \right\vert >\eta |\mathfrak{W}_{n}\right) \\
&\leq &\frac{1}{\eta ^{2}}\frac{1}{n}\sum_{t=1}^{n}\sup_{\lambda ,\tilde{%
\lambda}\in \Lambda :||\lambda -\tilde{\lambda}||\leq \delta }\left\{
F(\lambda ^{\prime }\mathcal{W}(x_{t}))-F(\tilde{\lambda}^{\prime }\mathcal{W%
}(x_{t}))\right\} ^{2}\leq \frac{1}{\eta ^{2}}\frac{1}{n}\sum_{t=1}^{n}%
\sup_{\lambda \in \Lambda }\left\Vert \frac{\partial }{\partial \lambda }%
F(\lambda ^{\prime }\mathcal{W}(x_{t}))\right\Vert ^{2}\times \delta ^{2}.
\end{eqnarray*}%
The Assumption \ref{assum:dgp}.c envelope bounds and ergodicity imply $%
1/n\sum_{t=1}^{n}\sup_{\lambda \in \Lambda }||(\partial /\partial \lambda
)F(\lambda ^{\prime }\mathcal{W}(x_{t}))||^{2}$ $\overset{p}{\rightarrow }$ $%
E[\sup_{\lambda \in \Lambda }||(\partial /\partial \lambda )F(\lambda
^{\prime }\mathcal{W}(x_{t}))||^{2}]$ $\leq $ $K$ $<$ $\infty $. Now pick $0$
$<$ $\delta $ $\leq $ $[\epsilon \eta ^{2}/K]^{1/2}$ to complete the proof
of stochastic equicontinuity asymptotically with probability approaching one
with respect to the draw $\mathfrak{W}_{n}$.

Next, for $1/\sqrt{n}\sum_{t=1}^{n}z_{t}\mathfrak{\hat{b}}_{\psi ,n}(\pi
,\lambda )^{\prime }\widehat{\mathcal{H}}_{\psi ,n}^{-1}(\pi )d_{\psi
,t}(\pi )$\ write
\begin{eqnarray}
&&\frac{1}{\sqrt{n}}\sum_{t=1}^{n}\sigma _{0}z_{t}\left\{ \mathfrak{\hat{b}}%
_{\psi ,n}(\pi ,\lambda )^{\prime }\widehat{\mathcal{H}}_{\psi ,n}^{-1}(\pi
)d_{\psi ,t}(\pi )-\mathfrak{\hat{b}}_{\psi ,n}(\tilde{\pi},\tilde{\lambda}%
)^{\prime }\widehat{\mathcal{H}}_{\psi ,n}^{-1}(\tilde{\pi})d_{\psi ,t}(%
\tilde{\pi})\right\}  \notag \\
&&\text{ \ \ \ \ \ \ \ \ \ \ \ \ \ \ \ }=\mathfrak{\hat{b}}_{\psi ,n}(\tilde{%
\pi},\tilde{\lambda})^{\prime }\widehat{\mathcal{H}}_{\psi ,n}^{-1}(\tilde{%
\pi})\frac{1}{\sqrt{n}}\sum_{t=1}^{n}\sigma _{0}z_{t}\left\{ d_{\psi ,t}(\pi
)-d_{\psi ,t}(\tilde{\pi})\right\}  \label{bbb} \\
&&\text{ \ \ \ \ \ \ \ \ \ \ \ \ \ \ \ \ \ \ \ \ \ \ \ }+\mathfrak{\hat{b}}%
_{\psi ,n}(\tilde{\pi},\tilde{\lambda})^{\prime }\widehat{\mathcal{H}}_{\psi
,n}^{-1}(\pi )\left\{ \widehat{\mathcal{H}}_{\psi ,n}(\tilde{\pi})-\widehat{%
\mathcal{H}}_{\psi ,n}(\pi )\right\} \widehat{\mathcal{H}}_{\psi ,n}^{-1}(%
\tilde{\pi})\frac{1}{\sqrt{n}}\sum_{t=1}^{n}\sigma _{0}z_{t}d_{\psi ,t}(\pi )
\notag \\
&&\text{ \ \ \ \ \ \ \ \ \ \ \ \ \ \ \ \ \ \ \ \ \ \ \ }+\left\{ \mathfrak{%
\hat{b}}_{\psi ,n}(\pi ,\lambda )-\mathfrak{\hat{b}}_{\psi ,n}(\tilde{\pi},%
\tilde{\lambda})\right\} ^{\prime }\frac{1}{\sqrt{n}}\sum_{t=1}^{n}\sigma
_{0}z_{t}\widehat{\mathcal{H}}_{\psi ,n}^{-1}(\pi )d_{\psi ,t}(\pi ).  \notag
\end{eqnarray}%
By Lemmas \ref{lm:H_ulln_weak} and \ref{lm:bn}, $\sup_{\pi \in \Pi }||%
\widehat{\mathcal{H}}_{\psi ,n}(\pi )$ $-$ $\mathcal{H}_{\psi }(\pi )||$ $%
\overset{p}{\rightarrow }$ $0$ and $\sup_{\pi \in \Pi ,\lambda \in \Lambda
}||\mathfrak{\hat{b}}_{\psi ,n}(\pi ,\lambda )$ $-$ $\mathfrak{b}_{\psi
}(\pi ,\lambda )||$ $\overset{p}{\rightarrow }$ $0$. Step 1.1 gives first
order expansions for both $1/\sqrt{n}\sum_{t=1}^{n}\sigma _{0}z_{t}\{d_{\psi
,t}(\pi )$ $-$ $d_{\psi ,t}(\tilde{\pi})\}$ and $\widehat{\mathcal{H}}_{\psi
,n}(\tilde{\pi})$ $-$ $\widehat{\mathcal{H}}_{\psi ,n}(\pi )$ around $\pi $.
Arguments there suffice to prove the first two summands in (\ref{bbb}) are
stochastically equicontinuous asymptotically with probability approaching
one with respect to the sample draw.

Consider the third summand in (\ref{bbb}). By the Step 1.1 argument and
Lemma \ref{lm:H_ulln_weak},\linebreak\ $\sup_{\pi \in \Pi }||1/\sqrt{n}%
\sum_{t=1}^{n}\sigma _{0}z_{t}\widehat{\mathcal{H}}_{\psi ,n}^{-1}(\pi
)d_{\psi ,t}(\pi )||$ $=$ $O_{p}(1)$.

Next, write $\mathfrak{\hat{b}}_{\psi ,n}(\chi )$ $=$ $\mathfrak{\hat{b}}%
_{\psi ,n}(\pi ,\lambda )$ where $\chi $ $=$ $[\pi ^{\prime },\lambda
^{\prime }]^{\prime }$ $\in $ $\mathcal{X}$ $\equiv $ $\Pi $ $\times $ $%
\Lambda $. Two applications of Minkowski's inequality yields:%
\begin{eqnarray*}
\left\vert \sup_{\chi ,\tilde{\chi}\in \mathcal{X}}\left\Vert \mathfrak{\hat{%
b}}_{\psi ,n}(\chi )-\mathfrak{\hat{b}}_{\psi ,n}(\tilde{\chi})\right\Vert
-\sup_{\chi ,\tilde{\chi}\in \mathcal{X}}\left\Vert \mathfrak{b}_{\psi
}(\chi )-\mathfrak{b}_{\psi }(\tilde{\chi})\right\Vert \right\vert &\leq
&\left\vert \sup_{\chi ,\tilde{\chi}\in \mathcal{X}}\left\Vert \mathfrak{%
\hat{b}}_{\psi ,n}(\chi )-\mathfrak{b}_{\psi }(\chi )-\mathfrak{\hat{b}}%
_{\psi ,n}(\tilde{\chi})+\mathfrak{b}_{\psi }(\tilde{\chi})\right\Vert
\right\vert \\
&\leq &\left\vert \sup_{\chi \in \mathcal{X}}\left\Vert \mathfrak{\hat{b}}%
_{\psi ,n}(\chi )-\mathfrak{b}_{\psi }(\chi )\right\Vert \right\vert
+\left\vert \sup_{\tilde{\chi}\in \mathcal{X}}\left\Vert \mathfrak{\hat{b}}%
_{\psi ,n}(\tilde{\chi})-\mathfrak{b}_{\psi }(\tilde{\chi})\right\Vert
\right\vert .
\end{eqnarray*}%
The right hand side is $o_{p}(1)$ by Lemma \ref{lm:bn}. Now apply the mean
value theorem, the Cauchy-Schwartz inequality, and the Assumption \ref%
{assum:dgp}.c envelope bounds to yield:
\begin{eqnarray*}
&&\sup_{\chi ,\tilde{\chi}\in \mathcal{X}:||\chi -\tilde{\chi}||\leq \delta
}\left\Vert \mathfrak{b}_{\psi }(\chi )-\mathfrak{b}_{n}(\tilde{\chi}%
)\right\Vert \leq \sup_{\chi ,\tilde{\chi}\in \mathcal{X}:||\chi -\tilde{\chi%
}||\leq \delta }\left\vert \mathfrak{b}_{\psi }(\chi )-\mathfrak{b}_{n}(%
\tilde{\chi})\right\vert \\
&&\text{ \ \ \ \ \ \ \ \ \ \ \ \ \ }\leq \sup_{\chi ,\tilde{\chi}\in
\mathcal{X}:||\chi -\tilde{\chi}||\leq \delta }\left\{ \left\vert E\left[
\frac{\partial }{\partial \chi }F(\lambda ^{\prime }\mathcal{W}%
(x_{t}))d_{\psi ,t}(\pi )\right] \right\vert \times \left\vert \chi -\tilde{%
\chi}\right\vert \right\} \\
&&\text{ \ \ \ \ \ \ \ \ \ \ \ \ \ }\leq \left\{ \left( \sup_{\lambda \in
\Lambda }E\left[ \left\vert \frac{\partial }{\partial \lambda }F(\lambda
^{\prime }\mathcal{W}(x_{t}))\right\vert ^{2}\right] \times \sup_{\pi \in
\Pi }E\left[ \left\vert d_{\psi ,t}(\pi )\right\vert ^{2}\right] \right)
^{1/2}\right. \\
&&\text{ \ \ \ \ \ \ \ \ \ \ \ \ \ \ \ \ \ \ \ }+\left. \left( \sup_{\lambda
\in \Lambda }E\left[ \left\{ F(\lambda ^{\prime }\mathcal{W}%
(x_{t})^{2}\right\} \right] \times \sup_{\pi \in \Pi }E\left[ \left\vert
\frac{\partial }{\partial \pi }d_{\psi ,t}(\pi )\right\vert ^{2}\right]
\right) ^{1/2}\right\} \times \delta =K\times \delta <\infty .
\end{eqnarray*}%
Stochastic equicontinuity asymptotically with probability approaching one
with respect to the sample draw therefore follows for the third summand in (%
\ref{bbb}) by the Step 1.1 argument.\medskip

\textbf{Step 1.4} We prove joint weak convergence in probability:%
\begin{equation}
\left\{ \mathfrak{\hat{T}}_{\psi ,n}^{\ast }(\pi ,\lambda ,\pi _{0},b),\hat{%
\pi}_{n}^{\ast }(\pi _{0},b):\Pi ,\Lambda \right\} \mathcal{\Rightarrow }%
^{p}\left\{ \mathfrak{T}_{\psi }(\pi ,\lambda ,b),\pi ^{\ast }(b):\Pi
,\Lambda \right\} .  \label{L1}
\end{equation}%
$\hat{\pi}_{n}^{\ast }(\pi _{0},b)$ is a continuous function of $\sigma _{0}%
\widehat{\mathcal{G}}_{\psi ,n}^{\ast }(\pi )$, $\widehat{\mathcal{H}}_{\psi
,n}(\pi )$ and $\mathcal{\hat{D}}_{\psi ,n}(\pi ,\pi _{0})$, where $\widehat{%
\mathcal{H}}_{\psi ,n}(\pi )$ and $\mathcal{\hat{D}}_{\psi ,n}(\pi ,\pi
_{0}) $ have uniform probability limits. Hence, by construction of $%
\mathfrak{\hat{T}}_{\psi ,n}^{\ast }(\pi ,\lambda ,\pi _{0},b)$, and uniform
convergence in probability of key summands in Steps 1.1-1.3, cf. Lemmas \ref%
{lm:H_ulln_weak}, \ref{lm:bn} and \ref{lm:ulln_xgF_eF}, it suffices to show:%
\begin{equation*}
\left\{ \sigma _{0}\mathfrak{\hat{Z}}_{\psi ,n}^{\ast }(\pi ,\lambda
),\sigma _{0}\widehat{\mathcal{G}}_{\psi ,n}^{\ast }(\pi ):\pi \in \Pi
,\lambda \in \Lambda \right\} \mathcal{\Rightarrow }^{p}\left\{ \mathfrak{Z}%
_{\psi }(\pi ,\lambda ),\mathcal{H}_{\psi }^{-1/2}(\pi )\mathcal{G}_{\psi
}(\pi ):\pi \in \Pi ,\lambda \in \Lambda \right\} .
\end{equation*}%
The required result then follows from the mapping theorem.

Let $r$ $=$ $[r_{1},r_{2}^{\prime }]^{\prime }$, $r_{1}$ $\in $ $\mathbb{R}$%
, $r_{2}$ $\in $ $\mathbb{R}^{k_{x}+k_{\beta }}$, $r^{\prime }r$ $=$ $1$,
and define $\mathcal{\hat{L}}_{n,t}(\pi ,\lambda ;r)$ $\equiv $ $r_{1}%
\mathcal{K}_{n,t}(\pi ,\lambda )$ $+$ $r_{2}^{\prime }\widehat{\mathcal{H}}%
_{\psi ,n}^{-1/2}(\pi )$ $\times $ $d_{\psi ,t}(\pi )$ and $\mathcal{L}%
_{n,t}(\pi ,\lambda ;r)$ $\equiv $ $r_{1}\mathcal{K}_{n,t}(\pi ,\lambda )$ $%
+ $ $r_{2}^{\prime }\mathcal{H}_{\psi }^{-1/2}(\pi )$ $\times $ $d_{\psi
,t}(\pi )$. Any linear combination $\sigma _{0}r_{1}\mathfrak{\hat{Z}}_{\psi
,n}^{\ast }(\pi ,\lambda )$ $+$ $\sigma _{0}r_{2}^{\prime }\widehat{\mathcal{%
G}}_{\psi ,n}^{\ast }(\pi )=1/\sqrt{n}\sum_{t=1}^{n}\sigma _{0}z_{t}\mathcal{%
L}_{n,t}(\pi ,\lambda ;r)$ is normally distributed with zero mean and
covariance kernel $\sigma _{0}^{2}n^{-1}\sum_{t=1}^{n}\mathcal{L}_{n,t}(\pi
,\lambda ;r)\mathcal{L}_{n,t}(\tilde{\pi},\tilde{\lambda};r)$. By arguments
used to prove Lemmas \ref{lm:bn} and \ref{lm:ulln_xgF_eF}:
\begin{eqnarray*}
&&\sup_{\pi ,\tilde{\pi}\in \Pi \times \Pi ,\lambda ,\tilde{\lambda}\in
\Lambda }\left\Vert \frac{1}{n}\sum_{t=1}^{n}\left\{ \mathcal{\hat{L}}%
_{n,t}(\pi ,\lambda ;r)\mathcal{\hat{L}}_{n,t}(\tilde{\pi},\tilde{\lambda}%
;r)-\mathcal{L}_{n,t}(\pi ,\lambda ;r)\mathcal{L}_{n,t}(\tilde{\pi},\tilde{%
\lambda};r)\right\} \right\Vert \overset{p}{\rightarrow }0 \\
&&\sup_{\pi ,\tilde{\pi}\in \Pi \times \Pi ,\lambda ,\tilde{\lambda}\in
\Lambda }\left\Vert \frac{1}{n}\sum_{t=1}^{n}\mathcal{L}_{n,t}(\pi ,\lambda
;r)\mathcal{L}_{n,t}(\tilde{\pi},\tilde{\lambda};r)-E\left[ \mathcal{L}%
_{n,t}(\pi ,\lambda ;r)\mathcal{L}_{n,t}(\tilde{\pi},\tilde{\lambda};r)%
\right] \right\Vert \overset{p}{\rightarrow }0.
\end{eqnarray*}%
$(\sigma _{0}\mathfrak{\hat{Z}}_{\psi ,n}^{\ast }(\pi ,\lambda ),\sigma _{0}%
\widehat{\mathcal{G}}_{\psi ,n}^{\ast }(\pi ))$ therefore convergences in
finite dimensional distributions to $(\mathfrak{Z}_{\psi }(\pi ,\lambda ),%
\mathcal{H}_{\psi }^{-1/2}(\pi )\mathcal{G}_{\psi }(\pi ))$. Stochastic
equicontinuity for $1/\sqrt{n}\sum_{t=1}^{n}\sigma _{0}z_{t}\mathcal{L}%
_{n,t}(\pi ,\lambda ;r)$ follows from Step 1.1 for $\widehat{\mathcal{G}}%
_{\psi ,n}^{\ast }(\pi )$, Step 1.3 for $\mathfrak{\hat{Z}}_{\psi ,n}^{\ast
}(\pi ,\lambda )$, and probability subadditivity.\medskip

\textbf{Step 1.5} Finally, consider the test statistic denominator $\hat{v}%
_{n}^{2}(\omega ,\pi ,\lambda )$. By Lemma \ref{lm:vn}:
\begin{equation}
\sup_{\substack{ \{\omega \in \mathbb{R}^{k_{\beta }}:\omega ^{\prime
}\omega =1\}  \\ \times \Pi \times \Lambda }}\left\vert \hat{v}%
_{n}^{2}(\omega ,\pi ,\lambda )-E\left[ \epsilon _{t}^{2}(\psi _{0},\pi
)\left\{ F\left( \lambda ^{\prime }\mathcal{W}(x_{t})\right) -\mathfrak{b}%
_{\theta }(\omega ,\pi ,\lambda )^{\prime }\mathcal{H}_{\theta }^{-1}(\omega
,\pi )d_{\theta ,t}(\omega ,\pi )\right\} ^{2}\right] \right\vert \overset{p}%
{\rightarrow }0.  \label{L2}
\end{equation}

In conclusion, by weak convergence in probability results in Steps 1.1 and
1.2, plus the mapping theorem, $\hat{\tau}_{\beta ,n}^{\ast }(\pi _{0},b)$ $%
\Rightarrow ^{p}$ $-\mathcal{S}_{\beta }\mathcal{H}_{\psi }^{-1}(\pi ^{\ast
}(b))\{\sigma _{0}\mathcal{G}_{\psi }(\pi ^{\ast }(b))$ $+$ $\mathcal{D}%
_{\psi }(\pi ^{\ast }(b),\pi _{0})b\}$ $=$ $\tau _{\beta }(\pi ^{\ast
}(b),b) $. Hence $\hat{\omega}_{n}^{\ast }(\pi _{0},b)$ $\Rightarrow ^{p}$ $%
\tau _{\beta }(\pi ^{\ast }(b),b)/||\tau _{\beta }(\pi ^{\ast }(b),b)||$.
The preceding limit, (\ref{L0})-(\ref{L2}), and the mapping theorem yield (%
\ref{Cn_boot_weak}).\medskip

\textbf{Step 2.} Assume $\mathcal{C}(i,b)$ with $||b||$ $<$ $\infty $. By
independence and the Glivenko-Cantelli theorem $\hat{p}_{n,\mathcal{M}%
}^{\ast }(\lambda ,h)$ $\overset{p}{\rightarrow }$ $P(\mathcal{\hat{T}}%
_{\psi ,n,1}^{\ast }(\lambda ,h)$ $\geq $ $\mathcal{T}_{n}(\lambda )|%
\mathfrak{W}_{n})$ as $\mathcal{M}$ $\rightarrow $ $\infty $. Now define $%
F_{n,\lambda }(c)$ $\equiv $ $P(\mathcal{T}_{n}(\lambda )$ $\leq $ $c)$,
hence $p_{n}(\lambda ,h)$ $\equiv $ $1$ $-$ $\mathcal{F}_{\lambda ,h}(%
\mathcal{T}_{n}(\lambda ))$.

$\{\mathcal{T}_{n}(\lambda )$ $:$ $\lambda $ $\in $ $\Lambda \}$ $%
\Rightarrow ^{\ast }$ $\{\mathcal{T}_{\psi }(\lambda ,b)$ $:$ $\lambda $ $%
\in $ $\Lambda \}$ by Theorem \ref{th:CM_weak}. By (\ref{Cn_boot_weak}),
conditionally on the sample draw $\mathfrak{W}_{n}$ $\equiv $ $%
\{(y_{t},x_{t})\}_{t=1}^{n}$, $\{\mathcal{\hat{T}}_{\psi ,n,j}^{\ast
}(\lambda ,h)\}_{j=1}^{\mathcal{M}}$ is a sequence of iid draws from $\{%
\mathcal{T}_{\psi }(\lambda ,h)$ $:\Lambda \}$, asymptotically with
probability approaching one with respect to $\mathfrak{W}_{n}$. Since $%
\mathcal{T}_{n}(\lambda )$, and $\mathcal{\hat{T}}_{\psi ,n,j}^{\ast
}(\lambda ,h)$ conditionally on $\mathfrak{W}_{n}$, have the same weak
limits in probability under $H_{0}$, uniformly on $\Lambda $, it follows
that $\sup_{c\geq 0}|P(\mathcal{\hat{T}}_{\psi ,n,j}^{\ast }(\lambda ,h)$ $%
\leq $ $c|\mathfrak{W}_{n})$ $-$ $F_{n,\lambda }(c)|$ $\overset{p}{%
\rightarrow }$ $0$ $\forall \lambda $ $\in $ $\Lambda $
\citep[Section 3, eq's (3.4) and
(3.5)]{GineZinn1990}. Therefore, as claimed $\hat{p}_{n,\mathcal{M}%
_{n}}^{\ast }(\lambda ,h)$ $=$ $1$ $-$ $F_{n,\lambda }(\mathcal{T}%
_{n}(\lambda ))$ $+$ $o_{p}(1)$ $=$ $p_{n}(\lambda ,h)$ $+$ $o_{p}(1)$ given
$\mathcal{M}_{n}$ $\rightarrow $ $\infty $.\medskip \newline
\textbf{Claim (b). }Recall $F_{n,\lambda }(c)$ $\equiv $ $P(\mathcal{T}%
_{n}(\lambda )$ $\leq $ $c)$ and $F_{n,\lambda ,h}^{\ast }(c)$ $\equiv $ $P(%
\mathcal{\hat{T}}_{\psi ,n,1}^{\ast }(\lambda ,h)$ $\leq $ $c|\mathfrak{W}%
_{n})$.\medskip

\textbf{Step 1. }In order to prove $\sup_{\lambda \in \Lambda }|\hat{p}_{n,%
\mathcal{M}_{n}}^{\ast }(\lambda ,h)$ $-$ $p_{n}(\lambda ,h)|$ $\overset{p}{%
\rightarrow }$ $0$, it suffices to show:%
\begin{eqnarray}
&&\sup_{\lambda \in \Lambda }\left\vert \hat{p}_{n,\mathcal{M}}^{\ast
}(\lambda ,h)-P\left( \mathcal{\hat{T}}_{\psi ,n,1}^{\ast }(\lambda ,h)\geq
\mathcal{T}_{n}(\lambda )|\mathfrak{W}_{n}\right) \right\vert  \label{sup_pP}
\\
&&\sup_{\lambda \in \Lambda }\left\vert P\left( \mathcal{\hat{T}}_{\psi
,n,1}^{\ast }(\lambda ,h)\geq \mathcal{T}_{n}(\lambda )|\mathfrak{W}%
_{n}\right) -\left\{ 1-F_{n,\lambda }(\mathcal{T}_{n}(\lambda ))\right\}
\right\vert \overset{p}{\rightarrow }0.  \label{sup_PF}
\end{eqnarray}

Consider (\ref{sup_pP}). Under (a) we have pointwise as $\mathcal{M}$ $%
\rightarrow $ $\infty $:
\begin{eqnarray*}
&&\hat{p}_{n,\mathcal{M}}^{\ast }(\lambda ,h)-P\left( \mathcal{\hat{T}}%
_{\psi ,n,1}^{\ast }(\lambda ,h)\geq \mathcal{T}_{n}(\lambda )|\mathfrak{W}%
_{n}\right) \\
&&\text{ \ \ \ \ \ \ \ \ \ \ \ \ \ \ \ }=\frac{1}{\mathcal{M}}\sum_{j=1}^{%
\mathcal{M}}\left\{ I\left( \mathcal{\hat{T}}_{\psi ,n,j}^{\ast }(\lambda
,h)>\mathcal{T}_{n}(\lambda )\right) -E\left[ I\left( \mathcal{\hat{T}}%
_{\psi ,n,j}^{\ast }(\lambda ,h)>\mathcal{T}_{n}(\lambda )\right) |\mathfrak{%
W}_{n}\right] \right\} \overset{p}{\rightarrow }0.
\end{eqnarray*}%
It remains to establish equicontinuity on $\Lambda $ for a uniform
Glivenko-Cantelli theorem (\ref{sup_pP})
\citep[Theorem
2.8.1]{VaartWellner1996}. The $\mathcal{V}(\mathcal{C})$ class satisfies the
required condition \citep[p. 168]{VaartWellner1996}. We therefore need only
demonstrate $\{I(\mathcal{\hat{T}}_{\psi ,n,j}^{\ast }(\lambda ,h)$ $>$ $%
\mathcal{T}_{n}(\lambda ))$ $-$ $E[I(\mathcal{\hat{T}}_{\psi ,n,j}^{\ast
}(\lambda ,h)$ $>$ $\mathcal{T}_{n}(\lambda ))|\mathfrak{W}_{n}]$ $:$ $%
\lambda $ $\in $ $\Lambda \}$ lies in $\mathcal{V}(\mathcal{C})$. In the
following we use properties of $\mathcal{V}(\mathcal{C})$\ functions without
citation. See, e.g., \citet[Chapt. 2.6]{VaartWellner1996}, and see the
discussion above Assumption \ref{assum_pvs_uniform}.

The test weight $F(\cdot )$ is in $\mathcal{V}(\mathcal{C})$ under
Assumption \ref{assum_pvs_uniform}, and trivially $\{\lambda ^{\prime
}W(x_{t})$ $:$ $\lambda $ $\in $ $\Lambda \}$ lies in $\mathcal{V}(\mathcal{C%
})$, hence $\{F(\lambda ^{\prime }W(x_{t}))$ $:$ $\lambda $ $\in $ $\Lambda
\}$ lies in $\mathcal{V}(\mathcal{C})$. Therefore $\{\mathcal{\hat{T}}_{\psi
,n,j}^{\ast }(\lambda ,h),\mathcal{T}_{n}(\lambda )$ $:$ $\lambda $ $\in $ $%
\Lambda \}$ are in $\mathcal{V}(\mathcal{C})$\ because they involve ratios
and products of linearly combined $\{F(\lambda ^{\prime
}W(x_{t}))\}_{t=1}^{n}$. Therefore $\{I(\mathcal{\hat{T}}_{\psi ,n,j}^{\ast
}(\lambda ,h)$ $>$ $\mathcal{T}_{n}(\lambda ))$ $:$ $\lambda $ $\in $ $%
\Lambda \}$ is in $\mathcal{V}(\mathcal{C})$.

Next, by Assumption \ref{assum_pvs_uniform} $\{F_{n,\lambda ,h}^{\ast }(c)$ $%
:$ $\lambda $ $\in $ $\Lambda ,c$ $\in $ $[0,\infty )\}$ belongs to the $%
\mathcal{V}(\mathcal{C})$ class. Hence $\{E[I(\mathcal{\hat{T}}_{\psi
,n,j}^{\ast }(\lambda ,h)$ $>$ $\mathcal{T}_{n}(\lambda ))|\mathfrak{W}_{n}]$
$:$ $\lambda $ $\in $ $\Lambda \}$ is in $\mathcal{V}(\mathcal{C})$, and
therefore $\{I(\mathcal{\hat{T}}_{\psi ,n,j}^{\ast }(\lambda ,h)$ $>$ $%
\mathcal{T}_{n}(\lambda ))$ $-$ $E[I(\mathcal{\hat{T}}_{\psi ,n,j}^{\ast
}(\lambda ,h)$ $>$ $\mathcal{T}_{n}(\lambda ))|\mathfrak{W}_{n}]$ $:$ $%
\lambda $ $\in $ $\Lambda \}$ lies in $\mathcal{V}(\mathcal{C})$ as
required, proving (\ref{sup_pP}).

Now consider (\ref{sup_PF}). Under (a) pointwise $P(\mathcal{\hat{T}}_{\psi
,n,1}^{\ast }(\lambda ,h)$ $\geq $ $\mathcal{T}_{n}(\lambda )|\mathfrak{W}%
_{n})$ $-\ \{1$ $-$ $F_{n,\lambda }(\mathcal{T}_{n}(\lambda ))\}$ $\overset{p%
}{\rightarrow }$ $0$. We have from above that $\mathcal{\hat{T}}_{\psi
,n,1}^{\ast }(\lambda ,h)$ and $\mathcal{T}_{n}(\lambda )$\ are in $\mathcal{%
V}(\mathcal{C})$. Under Assumption \ref{assum_pvs_uniform} $\{F_{n,\lambda
}(c),F_{n,\lambda ,h}^{\ast }(c)$ $:$ $\lambda $ $\in $ $\Lambda ,c$ $\in $ $%
[0,\infty )\}$ lie in $\mathcal{V}(\mathcal{C})$. Hence $\{P(\mathcal{\hat{T}%
}_{\psi ,n,1}^{\ast }(\lambda ,h)$ $\geq $ $\mathcal{T}_{n}(\lambda )|%
\mathfrak{W}_{n})$ $-\ \{1$ $-$ $F_{n,\lambda }(\mathcal{T}_{n}(\lambda ))\}$
$:$ $\lambda $ $\in $ $\Lambda \}$ lies in $\mathcal{V}(\mathcal{C})$,
promoting the required uniform convergence (\ref{sup_PF})
\citep[cf.][Theorem
2.8.1]{VaartWellner1996}.\medskip

\textbf{Step 2. }Finally, we prove $AsySz^{\ast }$ $\leq $ $\alpha $.
Consider the LF p-value $\hat{p}_{n,\mathcal{M}_{n}}^{(LF)}(\lambda )$. $%
AsySz^{\ast }$ can be written as:%
\begin{eqnarray*}
AsySz^{\ast } &=&\sup_{\lambda \in \Lambda }\limsup_{n\rightarrow \infty
}\sup_{\gamma \in \Gamma ^{\ast }}P_{\gamma }\left( \max \left\{ \sup_{h\in
\mathfrak{H}}\left\{ \hat{p}_{n,\mathcal{M}_{n}}^{\ast }(\lambda ,h)\right\}
,p_{n}^{\infty }(\lambda )\right\} <\alpha |H_{0}\right) \\
&=&\sup_{\lambda \in \Lambda }\limsup_{n\rightarrow \infty }\sup_{\gamma \in
\Gamma ^{\ast }}P_{\gamma }\left( \sup_{h\in \mathfrak{H}}\left\{ \hat{p}_{n,%
\mathcal{M}_{n}}^{\ast }(\lambda ,h),\mathcal{\bar{F}}_{\infty }(\mathcal{T}%
_{n}(\lambda ))\right\} <\alpha |H_{0}\right) .
\end{eqnarray*}%
By Step 1 $\sup_{\lambda \in \Lambda }|\hat{p}_{n,\mathcal{M}_{n}}^{\ast
}(\lambda ,h)$ $-$ $p_{n}(\lambda ,h)|$ $\overset{p}{\rightarrow }$ $0$,
hence
\begin{equation*}
AsySz^{\ast }=\sup_{\lambda \in \Lambda }\limsup_{n\rightarrow \infty
}\sup_{\gamma \in \Gamma ^{\ast }}P_{\gamma }\left( \sup_{h\in \mathfrak{H}%
}\left\{ \mathcal{\bar{F}}_{\lambda ,h}(\mathcal{T}_{n}(\lambda
))+o_{p,\lambda }(1),\mathcal{\bar{F}}_{\infty }(\mathcal{T}_{n}(\lambda
))\right\} <\alpha |H_{0}\right) \equiv \mathfrak{A}^{\ast },
\end{equation*}%
say. By Theorem \ref{th:CM_weak}.a, $\{\mathcal{T}_{n}(\lambda )$ $:$ $%
\Lambda \}$ $\Rightarrow ^{\ast }$ $\{\mathcal{T}_{\psi }(\lambda ,h)$ $:$ $%
\Lambda \}$ under $\mathcal{C}(i,b)$ with $||b||$ $<$ $\infty $. Weak
convergence implies convergence in finite dimensional distributions. By the
definition of distribution convergence, and the mapping theorem, weak
convergence therefore yields:%
\begin{equation*}
\mathfrak{A}^{\ast }=\sup_{\lambda \in \Lambda }\sup_{\tilde{h}\in \mathfrak{%
H}}P\left( \max \left\{ \sup_{h\in \mathfrak{H}}\left\{ \mathcal{\bar{F}}%
_{\lambda ,h}\left( \mathcal{T}_{\psi }(\lambda ,\tilde{h})\right) \right\} ,%
\mathcal{\bar{F}}_{\infty }\left( \mathcal{T}_{\psi }(\lambda ,\tilde{h}%
)\right) \right\} <\alpha \right) .
\end{equation*}

The remainder of the proof under $\mathcal{C}(i,b)$ and $\mathcal{C}%
(ii,\omega _{0})$, and for $\hat{p}_{n,\mathcal{M}_{n}}^{(ICS-1)}(\lambda )$%
, now follows directly from the proof of Theorem \ref{th:pv_weak}. $\mathcal{%
QED}$.\bigskip \newline
\textbf{Proof of Theorem \ref{th:pvot}.}\qquad Recall $\hat{p}_{n,\mathcal{M}%
}^{(\cdot )}(\lambda )$\ is the LF or ICS-1 p-value computed with the weak
identification p-value approximation $\hat{p}_{n,\mathcal{M}}^{\ast
}(\lambda ,h)$. Write $\mathcal{\hat{P}}_{n}(\alpha )$ $=$ $\mathcal{\hat{P}}%
_{n,\mathcal{M}_{n}}(\alpha )$ $\equiv $ $\int_{\Lambda }I(\hat{p}_{n,%
\mathcal{M}_{n}}^{(\cdot )}(\lambda )$ $<$ $\alpha )d\lambda $. Define the
infeasible PVOT $\mathcal{P}_{n}^{(\cdot )}(\alpha )$ $\equiv $ $%
\int_{\Lambda }I(p_{n}^{(\cdot )}(\lambda )$ $<$ $\alpha )d\lambda $, where $%
p_{n}^{(\cdot )}(\lambda )$ is the infeasible LF or ICS-1 p-value (see
Section \ref{sec:CM_pv_construct}).

By Theorem \ref{th:boot_pv}.a with $\mathcal{M}$ $=$ $\mathcal{M}_{n}$ $%
\rightarrow $ $\infty $ as $n$ $\rightarrow $ $\infty $, and Lebesgue's
dominated convergence theorem: $\mathcal{\hat{P}}_{n}(\alpha )$ $-$ $%
\mathcal{P}_{n}(\alpha )$ $=$ $\int_{\Lambda }\{I(\hat{p}_{n,\mathcal{M}%
_{n}}^{(\cdot )}(\lambda )$ $<$ $\alpha )$ $-$ $I(p_{n}^{(\cdot )}(\lambda )$
$<$ $\alpha )\}d\lambda $ $\overset{p}{\rightarrow }$ $0$. It therefore
suffices to prove the claim for $\mathcal{P}_{n}(\alpha )$.\medskip \newline
\textbf{Step 1 (}$H_{0}$\textbf{).} Consider identification case $\mathcal{C}%
(i,b)$. The LF p-value satisfies $p_{n}^{(LF)}(\lambda )$ $\geq $ $%
p_{n}(\lambda ,h)$ $=$ $\mathcal{\bar{F}}_{\lambda ,h}(\mathcal{T}%
_{n}(\lambda ))$ for any fixed $h$, hence in the LF case:%
\begin{equation}
P\left( \mathcal{P}_{n}(\alpha )>\alpha \right) \leq P\left( \int_{\Lambda
}I\left( \mathcal{\bar{F}}_{\lambda ,h}(\mathcal{T}_{n}(\lambda ))<\alpha
\right) d\lambda >\alpha \right) .  \label{PLF}
\end{equation}%
By Theorem \ref{th:CM_weak}.a and Assumption \ref{assum:pvs}, and using the
notation of Section \ref{sec:CM_pv_construct}, $\{\mathcal{T}_{n}(\lambda )$
$:$ $\lambda $ $\in $ $\Lambda \}$ $\Rightarrow ^{\ast }$ $\{\mathcal{T}%
_{\psi }(\lambda ,h)$ $:$ $\lambda $ $\in $ $\Lambda \}$. The limit process $%
\mathcal{T}_{\psi }(\lambda ,h)$ satisfies Assumption 1.a in \cite{Hill2018}%
, and $\mathcal{\bar{F}}_{\lambda ,h}(\mathcal{T}_{n}(\lambda ))$ trivially
satisfies Assumption 1.b in \cite{Hill2018}. Hence, $\lim_{n\rightarrow
\infty }P(\int_{\Lambda }I(\mathcal{\bar{F}}_{\lambda ,h}(\mathcal{T}%
_{n}(\lambda ))$ $<$ $\alpha )d\lambda $ $>$ $\alpha )$ $\leq $ $\alpha $ by
Theorem 3.1 in \cite{Hill2018}. In view of (\ref{PLF}), this proves $%
\lim_{n\rightarrow \infty }P(\mathcal{P}_{n}(\alpha )$ $>$ $\alpha )$ $\leq $
$\alpha $.

The ICS-1 p-value satisfies $p_{n}^{(ICS-1)}(\lambda )$ $\geq $ $%
p_{n}^{(LF)}(\lambda )I(\mathcal{A}_{n}$ $\leq $ $\kappa _{n})$ $\geq $ $%
\mathcal{\bar{F}}_{\lambda ,h}(\mathcal{T}_{n}(\lambda ))I(\mathcal{A}_{n}$ $%
\leq $ $\kappa _{n})$, hence in the ICS-1 case:%
\begin{equation*}
P\left( \mathcal{P}_{n}(\alpha )>\alpha \right) \leq P\left( \int_{\Lambda
}I\left( \mathcal{\bar{F}}_{\lambda ,h}(\mathcal{T}_{n}(\lambda ))I\left(
\mathcal{A}_{n}\leq \kappa _{n}\right) <\alpha \right) d\lambda >\alpha
\right) .
\end{equation*}%
Assumption 1.a in \cite{Hill2018} holds by the above argument. By the
Cauchy-Schwartz inequality, Theorem \ref{th:ICS_stat}.a and $\kappa _{n}$ $%
\rightarrow $ $\infty $:
\begin{eqnarray*}
E\left\vert \mathcal{\bar{F}}_{\lambda ,h}(\mathcal{T}_{n}(\lambda ))I\left(
\mathcal{A}_{n}\leq \kappa _{n}\right) -\mathcal{\bar{F}}_{\lambda ,h}(%
\mathcal{T}_{n}(\lambda ))\right\vert &=&E\left\vert \mathcal{\bar{F}}%
_{\lambda ,h}(\mathcal{T}_{n}(\lambda ))I\left( \mathcal{A}_{n}>\kappa
_{n}\right) \right\vert \\
&\leq &\left( E\left[ \mathcal{\bar{F}}_{\lambda ,h}^{2}(\mathcal{T}%
_{n}(\lambda ))\right] \right) ^{1/2}P\left( \mathcal{A}_{n}>\kappa
_{n}\right) ^{1/2}\rightarrow 0.
\end{eqnarray*}%
Hence $\mathcal{\bar{F}}_{\lambda ,h}(\mathcal{T}_{n}(\lambda ))I(\mathcal{A}%
_{n}\leq \kappa _{n})$ $-$ $\mathcal{\bar{F}}_{\lambda ,h}(\mathcal{T}%
_{n}(\lambda ))$ $\overset{p}{\rightarrow }$ $0$ by Markov's inequality,
which verifies Assumption 1.b in \cite{Hill2018}. Theorem 3.1 in \cite%
{Hill2018} now yields $\lim_{n\rightarrow \infty }P(\int_{\Lambda }I(%
\mathcal{\bar{F}}_{\lambda ,h}(\mathcal{T}_{n}(\lambda ))I(\mathcal{A}_{n}$ $%
\leq $ $\kappa _{n})$ $<$ $\alpha )d\lambda $ $>$ $\alpha )$ $\leq $ $\alpha
$, proving $\lim_{n\rightarrow \infty }P(\mathcal{\hat{P}}_{n}(\alpha )$ $>$
$\alpha )$ $\leq $ $\alpha $.

Under identification case $\mathcal{C}(ii,\omega _{0})$ the claim follows
from Theorem \ref{th:CM_weak}.b, and Theorem 3.1 in \cite{Hill2018}.\medskip
\newline
\textbf{Step 2 (}$H_{1}$\textbf{).} By Theorem \ref{th:pv_weak}.b, $%
p_{n}^{(\cdot )}(\lambda )$ $\overset{p}{\rightarrow }$ $0$ $\forall \lambda
$ $\in $ $\Lambda /S$ where $S$\ has Lebesgue measure zero. The claim now
follows from Theorem 2.2.b in \cite{Hill2018}. $\mathcal{QED}$.\bigskip
\newline
\textbf{Proof of Theorem \ref{th:pvot_size}.} Consider the LF case. By
Theorem \ref{th:boot_pv} $\sup_{\lambda \in \Lambda }|\hat{p}_{n,\mathcal{M}%
}^{\ast }(\lambda ,h)$ $-$ $p_{n}(\lambda ,h)|$ $\overset{p}{\rightarrow }$ $%
0$, and by Theorem \ref{th:CM_weak}.a, $\{\mathcal{T}_{n}(\lambda )$ $:$ $%
\Lambda \}$ $\Rightarrow ^{\ast }$ $\{\mathcal{T}_{\psi }(\lambda ,h)$ $:$ $%
\Lambda \}$ under $\mathcal{C}(i,b)$ with $||b||$ $<$ $\infty $. An
application of the mapping theorem, compactness of $\Lambda $ and dominated
convergence yield:%
\begin{eqnarray*}
AsySz(pvot) &=&\limsup_{n\rightarrow \infty }\sup_{\gamma \in \Gamma ^{\ast
}}P_{\gamma }\left( \int_{\Lambda }I\left( \max \left\{ \sup_{h\in \mathfrak{%
H}}\left\{ \hat{p}_{n,\mathcal{M}_{n}}^{\ast }(\lambda ,h)\right\}
,p_{n}^{\infty }(\lambda )\right\} <\alpha \right) d\lambda >\alpha
|H_{0}\right) \\
&=&\limsup_{n\rightarrow \infty }\sup_{\gamma \in \Gamma ^{\ast }}P_{\gamma
}\left( \int_{\Lambda }I\left( \max \left\{ \sup_{h\in \mathfrak{H}}\left\{
\mathcal{\bar{F}}_{\lambda ,h}(\mathcal{T}_{n}(\lambda ))\right\}
,p_{n}^{\infty }(\lambda )\right\} <\alpha \right) d\lambda +o_{p}(1)>\alpha
|H_{0}\right) \\
&=&\sup_{\tilde{h}\in \mathfrak{H}}P\left( \int_{\Lambda }I\left( \max
\left\{ \sup_{h\in \mathfrak{H}}\left\{ \mathcal{\bar{F}}_{\lambda ,h}(%
\mathcal{T}_{\psi }(\lambda ,\tilde{h}))\right\} ,p_{n}^{\infty }(\lambda
)\right\} <\alpha \right) d\lambda >\alpha |H_{0}\right) \\
&\leq &\sup_{h\in \mathfrak{H}}P\left( \int_{\Lambda }I\left( \mathcal{\bar{F%
}}_{\lambda ,h}(\mathcal{T}_{\psi }(\lambda ,h))<\alpha \right) d\lambda
>\alpha |H_{0}\right) .
\end{eqnarray*}%
Since $\mathcal{T}_{\psi }(\lambda ,h)$ is distributed $\mathcal{\bar{F}}%
_{\lambda ,h}$ we have $P(\int_{\Lambda }I(\mathcal{\bar{F}}_{\lambda ,h}(%
\mathcal{T}_{\psi }(\lambda ,h))$ $<$ $\alpha )d\lambda $ $>$ $\alpha
|H_{0}) $ $\leq $ $\alpha $ for all $h$ $\in $ $\mathfrak{H}$
\citep[see the
proof of Theorem 3.1 in][]{Hill2018}, hence $AsySz(pvot)$ $\leq $ $\alpha $.
The ICS-1 case, and under $\mathcal{C}(ii,\omega _{0})$, follow similarly. $%
\mathcal{QED}$.

\setcounter{equation}{0}

\section{APPENDIX : SUPPORTING RESULTS\label{app:lemmas}}

The following supporting results are proved in the supplemental material %
\citet[Appendix B]{Supp_Mat_2020}. Let \underline{$\iota $}$(A)$ and $\bar{%
\iota}(A)$\ denote the minimum and maximum eigenvalues of matrix $A$.

All subsequent Gaussian processes have \emph{almost surely} uniformly
continuous and bounded sample paths, hence in many cases we just say \textit{%
Gaussian process}. The first two lemmas cover uniform asymptotics for $%
\mathcal{G}_{\psi ,n}(\theta )$\ and $\widehat{\mathcal{H}}_{\psi ,n}(\pi )$%
\ under non- and weak identification.

\begin{lemma}
\label{lm:G_uclt_weak}Under $\mathcal{C}(i,b)$ and Assumption \ref{assum:dgp}%
, $\{\mathcal{G}_{\psi ,n}(\theta )$ $:$ $\theta $ $\in $ $\Theta \}$ $%
\Rightarrow ^{\ast }$ $\{\mathcal{G}_{\psi }(\theta )$ $:$ $\theta $ $\in $ $%
\Theta \}$, a zero mean Gaussian process with \emph{almost surely} uniformly
continuous and bounded sample paths and covariance $E[\mathcal{G}_{\psi
}(\theta )\mathcal{G}_{\psi }(\tilde{\theta})^{\prime }]$, $||E[\mathcal{G}%
_{\psi }(\theta )\mathcal{G}_{\psi }(\theta )^{\prime }]||$ $<$ $\infty $.
\end{lemma}

\begin{lemma}
\label{lm:H_ulln_weak}Under $\mathcal{C}(i,b)$ and Assumption \ref{assum:dgp}%
, $\sup_{\pi \in \Pi }||\widehat{\mathcal{H}}_{\psi ,n}(\pi )$ $- $ $%
\mathcal{H}_{\psi }(\pi )||$ $\overset{p}{\rightarrow }$ $0$, where
\underline{$\iota $}$(\mathcal{H}_{\psi }(\pi ))$ $>$ $0$ and $\bar{\iota}(%
\mathcal{H}_{\psi }(\pi ))$ $<$ $\infty $ for each $\pi $ $\in $ $\Pi $.
\end{lemma}

The next two lemmas are uniform asymptotics under semi-strong and strong
identification.

\begin{lemma}
\label{lm:G_uclt_strong}Under $\mathcal{C}(ii,\omega _{0})$ and Assumption %
\ref{assum:dgp}, $\{\mathcal{G}_{\theta ,n}(\theta )$ $:$ $\theta $ $\in $ $%
\Theta \}$ $\Rightarrow ^{\ast }$ $\{\mathcal{G}_{\theta }(\theta ) $ $:$ $%
\theta $ $\in $ $\Theta \}$, a zero mean Gaussian process with \emph{almost
surely} uniformly continuous and bounded sample paths.
\end{lemma}

\begin{corollary}
\label{cor:G_uclt_strong}Let $\theta _{n}$ $\equiv $ $[\beta _{n}^{\prime
},\zeta _{0}^{\prime },\pi _{0}^{\prime }]^{\prime }$ be the sequence of
true values under local drift $\{\beta _{n}\}$. Under $\mathcal{C}(ii,\omega
_{0})$ and Assumption \ref{assum:dgp}, $\sqrt{n}\mathfrak{B}(\beta
_{n})^{-1}(\partial /\partial \theta )Q_{n}(\theta _{n})$ $\overset{d}{%
\rightarrow }$ $\mathcal{G}_{\theta }$, a zero mean Gaussian law with a
finite, positive definite covariance $E[\mathcal{G}_{\theta }\mathcal{G}%
_{\theta }^{\prime }]$, and has a version that has \emph{almost surely}
uniformly continuous and bounded sample paths. Moreover, $E[\mathcal{G}%
_{\theta }\mathcal{G}_{\theta }^{\prime }]$ $=$ $\sigma _{0}^{2}E[d_{\theta
,t}d_{\theta ,t}^{\prime }]$ under $H_{0}$.
\end{corollary}

\begin{lemma}
\label{lm:H_ulln_strong}Under $\mathcal{C}(ii,\omega _{0})$ and Assumption %
\ref{assum:dgp}, $\widehat{\mathcal{H}}_{n}$ $\overset{p}{\rightarrow }$ $%
\mathcal{H}_{\theta }$, and \underline{$\iota $}$(\mathcal{H}_{\theta })$ $>$
$0$ and $\bar{\iota}(\mathcal{H}_{\theta })$ $<$ $\infty $.
\end{lemma}

Next, we tackle general versions of $\widehat{\mathcal{H}}_{n}$ and $%
\mathcal{\hat{V}}_{n}$ in (\ref{Hn}) and (\ref{JV_hat}) that are required
for uniform asymptotics. Recall we use $\mathcal{\hat{V}}_{n}$ for the
Identification Category Selection statistic $\mathcal{A}_{n}$ $\equiv $ $%
(k_{\beta }^{-1}n\hat{\beta}_{n}^{\prime }\hat{\Sigma}_{\beta ,\beta ,n}^{-1}%
\hat{\beta}_{n})^{1/2}$ in (\ref{A_ICS}), where $\hat{\Sigma}_{\beta ,\beta
,n}$ is the upper $k_{\beta }$ $\times $ $k_{\beta }$ block of $\hat{\Sigma}%
_{n}$ $\equiv $ $\widehat{\mathcal{H}}_{n}^{-1}\mathcal{\hat{V}}_{n}\widehat{%
\mathcal{H}}_{n}^{-1}$.

Define the augmented parameter and space $\theta ^{+}$ $\equiv $ $[||\beta
||,\omega ^{\prime },\zeta ^{\prime },\pi ^{\prime }]^{\prime }$ $\in $ $%
\Theta ^{+}$ $\equiv $ $\{\theta ^{+}$ $\in $ $\mathbb{R}^{k_{x}+k_{\beta
}+k_{\pi }+1}$ $:$ $\theta ^{+}$ $=$ $[||\beta ||,\omega (\beta ),\zeta ,\pi
]^{\prime }$ $:$ $\beta $ $\in $ $\mathcal{B},$ $\zeta $ $\in $ $\mathcal{Z}%
(\beta ),$ $\pi $ $\in $ $\Pi \}$. Define $\epsilon _{t}(\theta ^{+})$ $%
\equiv $ $y_{t}$ $-$ $\zeta ^{\prime }x_{t}$ $-$ $||\beta ||\omega ^{\prime
}g(x_{t},\pi )$, and:%
\begin{equation*}
\widehat{\mathcal{H}}_{n}(\theta ^{+})\equiv \frac{1}{n}\sum_{t=1}^{n}d_{%
\theta ,t}(\omega \left( \beta \right) ,\pi )d_{\theta ,t}(\omega \left(
\beta \right) ,\pi )^{\prime }\text{, \ }\mathcal{\hat{V}}_{n}(\theta
^{+})\equiv \frac{1}{n}\sum_{t=1}^{n}\epsilon _{t}^{2}(\theta ^{+})d_{\theta
,t}(\omega \left( \beta \right) ,\pi )d_{\theta ,t}(\omega \left( \beta
\right) ,\pi )^{\prime }.
\end{equation*}%
Hence $\widehat{\mathcal{H}}_{n}(\hat{\theta}_{n}^{+})$ $=$ $\widehat{%
\mathcal{H}}_{n}$ and $\mathcal{\hat{V}}_{n}(\hat{\theta}_{n}^{+})$ $=$ $%
\mathcal{\hat{V}}_{n}$. Define $\mathcal{H}_{\theta }(\theta ^{+})$ $\equiv $
$E[d_{\theta ,t}(\omega ,\pi )d_{\theta ,t}(\omega ,\pi )^{\prime }]$ and $%
\mathcal{V}(\theta ^{+})$ $\equiv $ \linebreak $E[\epsilon _{t}^{2}(\theta
^{+})d_{\theta ,t}(\omega ,\pi )d_{\theta ,t}(\omega ,\pi )^{\prime }]$. In
the interest of decreasing notation we use the same argument $\theta ^{+}$
for both $\widehat{\mathcal{H}}_{n}(\theta ^{+})$ and $\mathcal{\hat{V}}%
_{n}(\theta ^{+})$, although $\widehat{\mathcal{H}}_{n}(\theta ^{+})$ only
depends on $(\omega (\beta ),\pi )$.

\begin{lemma}
\label{lm:JV}Under Assumption \ref{assum:dgp} $\sup_{\theta ^{+}\in \Theta
^{+}:}||\widehat{\mathcal{H}}_{n}(\theta ^{+})$ $-$ $\mathcal{H}_{\theta
}(\theta ^{+})||$ $\overset{p}{\rightarrow }$ $0$, $\sup_{\pi \in \Pi }||%
\mathcal{\hat{D}}_{\psi ,n}(\pi ,\pi _{0})$ $-$ $\mathcal{D}_{\psi }(\pi )||$
$\overset{p}{\rightarrow }$ $0$, and $\sup_{\theta ^{+}\in \Theta ^{+}:}||%
\mathcal{\hat{V}}_{n}(\theta ^{+})$ $-$ $\mathcal{V}(\theta ^{+})||$ $%
\overset{p}{\rightarrow }$ $0$, where $\inf_{\theta ^{+}\in \Theta ^{+}:}$%
\underline{$\iota $}$(\mathcal{H}_{\theta }(\theta ^{+}))$ $>$ $0$, $\bar{%
\iota}(\mathcal{H}_{\theta })$ $<$ $\infty $, $\inf_{\theta ^{+}\in \Theta
^{+}:}$\underline{$\iota $}$(\mathcal{V}(\theta ^{+}))$ $>$ $0$, and $\bar{%
\iota}(\mathcal{V}_{\theta })$ $<$ $\infty $.
\end{lemma}

In order to characterize the weak limit of $\hat{\pi}_{n}$, we need the
following results. \cite{AndrewsCheng2012} exploit the following normalizing
constants for the criterion derivative under $\mathcal{C}(i,b)$:
\begin{equation*}
a_{n}=\left\{
\begin{array}{ll}
\sqrt{n} & \text{if }\mathcal{C}(i,b)\text{ and }\left\Vert b\right\Vert
<\infty \\
\left\Vert \beta _{n}\right\Vert ^{-1} & \text{if }\mathcal{C}(i,b)\text{
and }\left\Vert b\right\Vert =\infty%
\end{array}%
\right.
\end{equation*}%
By construction $||\beta _{n}||^{-1}$ $\leq $ $\sqrt{n}$ for large $n$ when $%
\sqrt{n}||\beta _{n}||$ $\rightarrow $ $\infty $ hence $a_{n}\leq $ $\sqrt{n}
$ for large $n$ when $||b||$ $=$ $\infty $. Recall $\psi _{0,n}$ $\equiv $ $%
[0_{k_{\beta }}^{\prime },\zeta _{0}^{\prime }]^{\prime }$ hence $%
Q_{0,n}\equiv $ $Q_{n}(\psi _{0,n},\pi )$ does not depend on $\pi $. Now
define:%
\begin{equation*}
\mathcal{Z}_{n}(\pi )=-a_{n}\widehat{\mathcal{H}}_{\psi ,n}^{-1}(\pi )\frac{%
\partial }{\partial \psi }Q_{n}(\psi _{0,n},\pi ).
\end{equation*}%
Under $\mathcal{C}(i,b)$, $\widehat{\mathcal{H}}_{\psi ,n}(\pi )$ is
positive definite uniformly on $\Pi $, asymptotically with probability
approaching one. See Lemma \ref{lm:H_ulln_weak}. Write $Q_{n}^{c}(\pi
)\equiv Q_{n}(\hat{\psi}_{n}(\pi ),\pi )$.

\begin{lemma}
\label{lm:aQQ}Let drift case $\mathcal{C}(i,b)$ and Assumption \ref%
{assum:dgp} hold.\medskip \newline
$a.$ In general $a_{n}(\hat{\psi}_{n}(\pi )$ $-$ $\psi _{0,n})$ $=$ $%
\mathcal{Z}_{n}(\pi )$.$\medskip $\newline
$b$. $a_{n}^{2}\{Q_{n}^{c}(\pi )$ $-$ $Q_{0,n}\}$ $=$ $-2^{-1}\mathcal{Z}%
_{n}(\pi )^{\prime }\widehat{\mathcal{H}}_{\psi ,n}(\pi )\mathcal{Z}_{n}(\pi
)$ where $Q_{0,n}\equiv $ $Q_{n}(\psi _{0,n},\pi )$.
\end{lemma}

Define $\vartheta _{\psi }(\pi ,\omega _{0})$ $\equiv $ $-2^{-2}\omega
_{0}^{\prime }\mathcal{D}_{\psi }(\pi )^{\prime }\mathcal{H}_{\psi
}^{-1}(\pi )\mathcal{D}_{\psi }(\pi )\omega _{0}$. The following is a key
result for characterizing the asymptotic properties of $\hat{\pi}_{n}$ under
weak identification.

\begin{lemma}
\label{lm:nQQ}Let drift case $\mathcal{C}(i,b)$ and Assumption \ref%
{assum:dgp} hold.\medskip \newline
$a.$ If $||b||$ $<$ $\infty $ then $\{n(Q_{n}^{c}(\pi )$ $-$ $Q_{0,n})$ $:$ $%
\pi $ $\in $ $\Pi \}$ $\Rightarrow ^{\ast }$ $\{\xi _{\psi }(\pi ,b)$ $:$ $%
\pi $ $\in $ $\Pi \}$.\medskip \newline
$b.$ If $||b||$ $=$ $\infty $ and $\beta _{n}/||\beta _{n}||$ $\rightarrow $
$\omega _{0}$ for some $\omega _{0}\in $ $\mathbb{R}^{k_{\beta }}$, $%
||\omega _{0}||$ $=$ $1$, then $\sup_{\pi \in \Pi }|(1/||\beta
_{n}||^{2})(Q_{n}^{c}(\pi )$ $-$ $Q_{0,n})-\vartheta _{\psi }(\pi ,\omega
_{0})|$ $\overset{p}{\rightarrow }$ $0$.
\end{lemma}

Write $\epsilon _{t}(\psi ,\pi )$ $=$ $y_{t}$ $-$ $\zeta ^{\prime }x_{t}$ $-$
$\beta ^{\prime }g(x_{t},\pi )$. Recall $\psi _{n}$ is the (possibly
drifting) true value under $H_{0}$.

\begin{lemma}
\label{lm:CM_weak}Let Assumption \ref{assum:dgp} hold.\medskip \newline
$a.$ Under $\mathcal{C}(i,b)$ with $||b||$ $<$ $\infty $:%
\begin{equation*}
\left\{ \frac{1}{\sqrt{n}}\sum_{t=1}^{n}\left\{ \epsilon _{t}(\psi _{n},\pi )%
\mathcal{K}_{\psi ,t}(\pi ,\lambda )-E\left[ \epsilon _{t}(\psi _{n},\pi )%
\mathcal{K}_{\psi ,t}(\pi ,\lambda )\right] \right\} :\Pi ,\Lambda \right\}
\Rightarrow ^{\ast }\left\{ \mathfrak{Z}_{\psi }(\pi ,\lambda ):\Pi ,\Lambda
\right\} ,
\end{equation*}%
a zero mean Gaussian process with covariance kernel $E[\mathfrak{Z}_{\psi
}(\pi ,\lambda )\mathfrak{Z}_{\psi }(\tilde{\pi},\tilde{\lambda})]$. Under $%
H_{0}$,
\begin{equation*}
\sup_{\pi \in \Pi ,\lambda \in \Lambda }\left\vert \frac{1}{\sqrt{n}}%
\sum_{t=1}^{n}\left\{ \epsilon _{t}(\psi _{n},\pi )\mathcal{K}_{\psi ,t}(\pi
,\lambda )-E\left[ \epsilon _{t}(\psi _{n},\pi )\mathcal{K}_{\psi ,t}(\pi
,\lambda )\right] \right\} -\frac{1}{\sqrt{n}}\sum_{t=1}^{n}\epsilon _{t}%
\mathcal{K}_{\psi ,t}(\pi ,\lambda )\right\vert \overset{p}{\rightarrow }0,
\end{equation*}%
and $E[\mathfrak{Z}_{\psi }(\pi ,\lambda )\mathfrak{Z}_{\psi }(\tilde{\pi},%
\tilde{\lambda})]$ $=$ $\sigma _{0}^{2}E[\mathcal{K}_{\psi ,t}(\pi ,\lambda )%
\mathcal{K}_{\psi ,t}(\tilde{\pi},\tilde{\lambda})]$.\medskip \newline
$b.$ Under $\mathcal{C}(i,\omega _{0})$, $\{1/\sqrt{n}\sum_{t=1}^{n}\epsilon
_{t}\mathcal{K}_{\theta ,t}(\lambda )$ $:$ $\lambda $ $\in $ $\Lambda \}$ $%
\Rightarrow ^{\ast }$ $\{\mathfrak{Z}_{\theta }$ $:$ $\lambda $ $\in $ $%
\Lambda \}$, a zero mean Gaussian process with covariance $E[\mathfrak{Z}%
_{\theta }(\lambda )\mathfrak{Z}_{\theta }\tilde{(}\lambda )]$ $=$ $%
E[\epsilon _{t}^{2}\mathcal{K}_{\theta ,t}(\lambda )\mathcal{K}_{\theta ,t}(%
\tilde{\lambda})]$.
\end{lemma}

\begin{lemma}
\label{lm:bn}Under Assumption \ref{assum:dgp}, $\sup_{\omega \in \mathbb{R}%
^{k_{\beta }}:||\omega ||=1,\pi \in \Pi ,\lambda \in \Lambda }||\mathfrak{%
\hat{b}}_{\theta ,n}(\omega ,\pi ,\lambda )$ $-$ $\mathfrak{b}_{\theta
}(\omega ,\pi ,\lambda )||$ $\overset{p}{\rightarrow }$ $0$ and $\sup_{\pi
\in \Pi ,\lambda \in \Lambda }||\mathfrak{\hat{b}}_{\psi ,n}(\pi ,\lambda )$
$-$ $\mathfrak{b}_{\psi }(\pi ,\lambda )||$ $\overset{p}{\rightarrow }$ $0$.
\end{lemma}

\begin{lemma}
\label{lm:vn}Under Assumption \ref{assum:dgp} $\sup_{\theta ^{+}\in \Theta
^{+},\lambda \in \Lambda }||\hat{v}_{n}^{2}(\theta ^{+},\lambda )-$ $%
v^{2}(\theta ^{+},\lambda )||$ $\overset{p}{\rightarrow }$ $0$ and $%
\sup_{\theta \in \Theta ,\lambda \in \Lambda }||\hat{v}_{n}^{2}(\theta
,\lambda )-$ $v^{2}(\theta ,\lambda )||$ $\overset{p}{\rightarrow }$ $0.$
\end{lemma}

Define $v^{2}(\lambda )$ $=$ $v^{2}(\omega _{0},\pi _{0},\lambda )$ where:%
\begin{equation*}
v^{2}(\omega ,\pi ,\lambda )\equiv E\left[ \epsilon _{t}^{2}(\psi _{0},\pi
)\left\{ F\left( \lambda ^{\prime }\mathcal{W}(x_{t})\right) -\mathfrak{b}%
_{\theta }(\omega ,\pi ,\lambda )^{\prime }\mathcal{H}_{\theta }^{-1}(\omega
,\pi )d_{\theta ,t}(\omega ,\pi )\right\} ^{2}\right] .
\end{equation*}

\begin{lemma}
\label{lm:v2_nondegen}Let Assumptions \ref{assum:dgp}.a(i) and \ref%
{assum:nondegen_v2_ae} hold. Under $\mathcal{C}(i,b)$ with $||b||$ $<$ $%
\infty $, the set $\{\lambda $ $\in $ $\Lambda $ $:$ \linebreak $%
\inf_{\omega ^{\prime }\omega =1,\pi \in \Pi }v^{2}(\omega ,\pi ,\lambda )$ $%
=$ $0\}$ has Lebesgue measure zero. Under $\mathcal{C}(ii,\omega _{0})$, the
set $\{\lambda $ $\in $ $\Lambda :v^{2}(\lambda )$ $=$ $0\}$ has Lebesgue
measure zero.
\end{lemma}

Define $\mathcal{M}_{t}(\pi ,\lambda )$ $\equiv $ $x_{t}\{g(x_{t},\pi _{0})$
$-$ $g(x_{t},\pi )\}F(\lambda ^{\prime }\mathcal{W}(x_{t}))$ and $\mathcal{%
\tilde{M}}_{t}(\pi )$ $\equiv $ $x_{t}\{g(x_{t},\pi ) $ $-$ $g(x_{t},\pi
_{0})\}d_{\psi ,t}(\pi )$.

\begin{lemma}
\label{lm:ulln_xgF_eF}Under Assumption \ref{assum:dgp}, $\sup_{\pi \in \Pi
,\lambda \in \Lambda }|1/n\sum_{t=1}^{n}\epsilon _{t}F(\lambda ^{\prime }%
\mathcal{W}(x_{t}))$ $-$ $E[\epsilon _{t}F(\lambda ^{\prime }\mathcal{W}%
(x_{t}))]|$ $\overset{p}{\rightarrow }$ $0$, $\sup_{\pi \in \Pi ,\lambda \in
\Lambda }||1/n\sum_{t=1}^{n}\mathcal{M}_{t}(\pi ,\lambda )$ $-$ $E[\mathcal{M%
}_{t}(\pi ,\lambda )]||$ $\overset{p}{\rightarrow }$ $0$ and $\sup_{\pi \in
\Pi ,\lambda \in \Lambda }||E[\mathcal{M}_{t}(\pi ,\lambda )]|| $ $<$ $%
\infty $, and \linebreak $\sup_{\pi \in \Pi }||1/n\sum_{t=1}^{n}\mathcal{%
\tilde{M}}_{t}(\pi )$ $-$ $E[\mathcal{\tilde{M}}_{t}(\pi )]||$ $\overset{p}{%
\rightarrow }$ $0$ and $\sup_{\pi \in \Pi }||E[\mathcal{\tilde{M}}_{t}(\pi
)]||$ $<$ $\infty $.
\end{lemma}

\begin{table}[h]
\caption{STAR Test Rejection Frequencies: Sample Size $n = 100$}
\label{tbl:starn100}
\begin{center}
\begin{tabular}{l|ccc|c|ccc|c|ccc}
\hline\hline
& \multicolumn{3}{|c|}{$H_{0}$: LSTAR} &  & \multicolumn{3}{|c|}{$H_{1}$-weak
} &  & \multicolumn{3}{|c}{$H_{1}$-strong} \\ \hline
& 1\% & 5\% & 10\% &  & 1\% & 5\% & 10\% &  & 1\% & 5\% & 10\% \\
\hline\hline
& \multicolumn{11}{c}{Strong Identification: $\beta _{n}=.3$} \\ \hline\hline
sup $\mathcal{T}_{n}$ & \multicolumn{1}{|l}{.025} & \multicolumn{1}{l}{.094}
& \multicolumn{1}{l|}{.163} & \multicolumn{1}{|l|}{} & \multicolumn{1}{|l}{
.147} & \multicolumn{1}{l}{.280} & \multicolumn{1}{l|}{.365} &
\multicolumn{1}{|l|}{} & \multicolumn{1}{|l}{.757} & \multicolumn{1}{l}{.872}
& \multicolumn{1}{l}{.907} \\
aver $\mathcal{T}_{n}$ & \multicolumn{1}{|l}{.025} & \multicolumn{1}{l}{.078}
& \multicolumn{1}{l|}{.135} & \multicolumn{1}{|l|}{} & \multicolumn{1}{|l}{
.087} & \multicolumn{1}{l}{.209} & \multicolumn{1}{l|}{.289} &
\multicolumn{1}{|l|}{} & \multicolumn{1}{|l}{.552} & \multicolumn{1}{l}{.726}
& \multicolumn{1}{l}{.804} \\ \hline
rand $\mathcal{T}_{n}$ & \multicolumn{1}{|l}{.011} & \multicolumn{1}{l}{.052}
& \multicolumn{1}{l|}{.096} & \multicolumn{1}{|l|}{} & \multicolumn{1}{|l}{
.053} & \multicolumn{1}{l}{.143} & \multicolumn{1}{l|}{.232} &
\multicolumn{1}{|l|}{} & \multicolumn{1}{|l}{.446} & \multicolumn{1}{l}{.635}
& \multicolumn{1}{l}{.732} \\
rand LF & \multicolumn{1}{|l}{\textbf{.007}} & \multicolumn{1}{l}{\textbf{%
.015}} & \multicolumn{1}{l|}{\textbf{.038}} & \multicolumn{1}{|l|}{} &
\multicolumn{1}{|l}{\textbf{.013}} & \multicolumn{1}{l}{\textbf{.066}} &
\multicolumn{1}{l|}{\textbf{.141}} & \multicolumn{1}{|l|}{} &
\multicolumn{1}{|l}{\textbf{.442}} & \multicolumn{1}{l}{\textbf{.553}} &
\multicolumn{1}{l}{\textbf{.661}} \\
rand ICS-1 & \multicolumn{1}{|l}{\textbf{.013}} & \multicolumn{1}{l}{\textbf{%
.050}} & \multicolumn{1}{l|}{\textbf{.089}} & \multicolumn{1}{|l|}{} &
\multicolumn{1}{|l}{\textbf{.028}} & \multicolumn{1}{l}{\textbf{.089}} &
\multicolumn{1}{l|}{\textbf{.170}} & \multicolumn{1}{|l|}{} &
\multicolumn{1}{|l}{\textbf{.379}} & \multicolumn{1}{l}{\textbf{.593}} &
\multicolumn{1}{l}{\textbf{.692}} \\ \hline
sup $p_{n}$ & .009 & .039 & .068 &  & .036 & .118 & .209 &  & .378 & .554 &
.656 \\
sup $p_{n}$ LF & \textbf{.006} & \textbf{.009} & \textbf{.032} &  & \textbf{%
.012} & \textbf{.057} & \textbf{.120} &  & \textbf{.262} & \textbf{.457} &
\textbf{.572} \\
sup $p_{n}$ ICS-1 & \textbf{.006} & \textbf{.036} & \textbf{.061} &  &
\textbf{.020} & \textbf{.081} & \textbf{.138} &  & \textbf{.310} & \textbf{%
.506} & \textbf{.617} \\ \hline
PVOT & \multicolumn{1}{|l}{.015} & \multicolumn{1}{l}{.065} &
\multicolumn{1}{l|}{.124} & \multicolumn{1}{|l|}{} & \multicolumn{1}{|l}{.101
} & \multicolumn{1}{l}{.257} & \multicolumn{1}{l|}{.335} &
\multicolumn{1}{|l|}{} & \multicolumn{1}{|l}{.727} & \multicolumn{1}{l}{.859}
& \multicolumn{1}{l}{.883} \\
PVOT LF & \multicolumn{1}{|l}{\textbf{.007}} & \multicolumn{1}{l}{\textbf{%
.014}} & \multicolumn{1}{l|}{\textbf{.052}} & \multicolumn{1}{|l|}{} &
\multicolumn{1}{|l}{\textbf{.026}} & \multicolumn{1}{l}{\textbf{.121}} &
\multicolumn{1}{l|}{\textbf{.208}} & \multicolumn{1}{|l|}{} &
\multicolumn{1}{|l}{\textbf{.552}} & \multicolumn{1}{l}{\textbf{.781}} &
\multicolumn{1}{l}{\textbf{.817}} \\
PVOT ICS-1 & \multicolumn{1}{|l}{\textbf{.007}} & \multicolumn{1}{l}{\textbf{%
.043}} & \multicolumn{1}{l|}{\textbf{.073}} & \multicolumn{1}{|l|}{} &
\multicolumn{1}{|l}{\textbf{.042}} & \multicolumn{1}{l}{\textbf{.153}} &
\multicolumn{1}{l|}{\textbf{.237}} & \multicolumn{1}{|l|}{} &
\multicolumn{1}{|l}{\textbf{.622}} & \multicolumn{1}{l}{\textbf{.815}} &
\multicolumn{1}{l}{\textbf{.842}} \\ \hline\hline
& \multicolumn{11}{c}{Weak Identification: $\beta _{n}=.3/\sqrt{n}$} \\
\hline\hline
sup $\mathcal{T}_{n}$ & \multicolumn{1}{|l}{.064} & \multicolumn{1}{l}{.155}
& \multicolumn{1}{l|}{.239} &  & \multicolumn{1}{|l}{.337} &
\multicolumn{1}{l}{.574} & \multicolumn{1}{l|}{.681} &  &
\multicolumn{1}{|l}{.929} & \multicolumn{1}{l}{.978} & \multicolumn{1}{l}{
.993} \\
aver $\mathcal{T}_{n}$ & \multicolumn{1}{|l}{.057} & \multicolumn{1}{l}{.146}
& \multicolumn{1}{l|}{.219} &  & \multicolumn{1}{|l}{.215} &
\multicolumn{1}{l}{.430} & \multicolumn{1}{l|}{.554} &  &
\multicolumn{1}{|l}{.739} & \multicolumn{1}{l}{.888} & \multicolumn{1}{l}{
.932} \\ \hline
rand $\mathcal{T}_{n}$ & \multicolumn{1}{|l}{.027} & \multicolumn{1}{l}{.083}
& \multicolumn{1}{l|}{.175} &  & \multicolumn{1}{|l}{.164} &
\multicolumn{1}{l}{.343} & \multicolumn{1}{l|}{.474} &  &
\multicolumn{1}{|l}{.604} & \multicolumn{1}{l}{.810} & \multicolumn{1}{l}{
.870} \\
rand LF & \textbf{.012} & \textbf{.042} & \textbf{.093} &  & \textbf{.060} &
\textbf{.161} & \textbf{.308} &  & \textbf{.467} & \textbf{.685} & \textbf{%
.794} \\
rand ICS-1 & \textbf{.012} & \textbf{.046} & \textbf{.104} &  & \textbf{.116}
& \textbf{.261} & \textbf{.382} &  & \textbf{.545} & \textbf{.749} & \textbf{%
.841} \\ \hline
sup $p_{n}$ & .019 & .087 & .145 &  & .107 & .253 & .411 &  & .493 & .700 &
.785 \\
sup $p_{n}$ LF & \textbf{.001} & \textbf{.061} & \textbf{.084} &  & \textbf{%
.036} & \textbf{.124} & \textbf{.230} &  & \textbf{.351} & \textbf{.598} &
\textbf{.698} \\
sup $p_{n}$ ICS-1 & \textbf{.001} & \textbf{.065} & \textbf{.085} &  &
\textbf{.088} & \textbf{.193} & \textbf{.335} &  & \textbf{.454} & \textbf{%
.663} & \textbf{.756} \\ \hline
PVOT & \multicolumn{1}{|l}{.038} & \multicolumn{1}{l}{.127} &
\multicolumn{1}{l|}{.196} &  & \multicolumn{1}{|l}{.328} &
\multicolumn{1}{l}{.542} & \multicolumn{1}{l|}{.591} &  &
\multicolumn{1}{|l}{.893} & \multicolumn{1}{l}{.968} & \multicolumn{1}{l}{
.950} \\
PVOT LF & \textbf{.015} & \textbf{.049} & \textbf{.108} &  &
\multicolumn{1}{|l}{\textbf{.108}} & \multicolumn{1}{l}{\textbf{.320}} &
\multicolumn{1}{l|}{\textbf{.398}} &  & \multicolumn{1}{|l}{\textbf{.710}} &
\multicolumn{1}{l}{\textbf{.911}} & \multicolumn{1}{l}{\textbf{.916}} \\
PVOT ICS-1 & \multicolumn{1}{|l}{\textbf{.014}} & \multicolumn{1}{l}{\textbf{%
.049}} & \multicolumn{1}{l|}{\textbf{.107}} &  & \multicolumn{1}{|l}{\textbf{%
.221}} & \multicolumn{1}{l}{\textbf{.435}} & \multicolumn{1}{l|}{\textbf{.486%
}} &  & \multicolumn{1}{|l}{\textbf{.830}} & \multicolumn{1}{l}{\textbf{.942}
} & \multicolumn{1}{l}{\textbf{.932}} \\ \hline\hline
& \multicolumn{11}{c}{Non-Identification: $\beta _{n}=\beta _{0}=0$} \\
\hline\hline
sup $\mathcal{T}_{n}$ & \multicolumn{1}{|l}{.066} & \multicolumn{1}{l}{.164}
& \multicolumn{1}{l|}{.249} &  & \multicolumn{1}{|l}{.358} &
\multicolumn{1}{l}{.584} & \multicolumn{1}{l|}{.696} &  &
\multicolumn{1}{|l}{.902} & \multicolumn{1}{l}{.970} & \multicolumn{1}{l}{
.983} \\
aver $\mathcal{T}_{n}$ & \multicolumn{1}{|l}{.062} & \multicolumn{1}{l}{.148}
& \multicolumn{1}{l|}{.226} &  & \multicolumn{1}{|l}{.233} &
\multicolumn{1}{l}{.438} & \multicolumn{1}{l|}{.548} &  &
\multicolumn{1}{|l}{.716} & \multicolumn{1}{l}{.872} & \multicolumn{1}{l}{
.911} \\ \hline
rand $\mathcal{T}_{n}$ & \multicolumn{1}{|l}{.044} & \multicolumn{1}{l}{.107}
& \multicolumn{1}{l|}{.186} &  & \multicolumn{1}{|l}{.184} &
\multicolumn{1}{l}{.380} & \multicolumn{1}{l|}{.505} &  &
\multicolumn{1}{|l}{.634} & \multicolumn{1}{l}{.793} & \multicolumn{1}{l}{
.864} \\
rand LF & \textbf{.013} & \textbf{.046} & \textbf{.115} &  & \textbf{.069} &
\textbf{.191} & \textbf{.327} &  & \textbf{.498} & \textbf{.725} & \textbf{%
.818} \\
rand ICS-1 & \textbf{.013} & \textbf{.047} & \textbf{.116} &  & \textbf{.137}
& \textbf{.298} & \textbf{.481} &  & \textbf{.583} & \textbf{.769} & \textbf{%
.847} \\ \hline
sup $p_{n}$ & .018 & .080 & .167 &  & .117 & .272 & .363 &  & .514 & .710 &
.807 \\
sup $p_{n}$ LF & \textbf{.011} & \textbf{.043} & \textbf{.083} &  & \textbf{%
.042} & \textbf{.122} & \textbf{.221} &  & \textbf{.383} & \textbf{.612} &
\textbf{.740} \\
sup $p_{n}$ ICS-1 & \textbf{.011} & \textbf{.044} & \textbf{.086} &  &
\textbf{.093} & \textbf{.205} & \textbf{.293} &  & \textbf{.464} & \textbf{%
.683} & \textbf{.783} \\ \hline
PVOT & \multicolumn{1}{|l}{.049} & \multicolumn{1}{l}{.134} &
\multicolumn{1}{l|}{.190} &  & \multicolumn{1}{|l}{.322} &
\multicolumn{1}{l}{.554} & \multicolumn{1}{l|}{.624} &  &
\multicolumn{1}{|l}{.890} & \multicolumn{1}{l}{.962} & \multicolumn{1}{l}{
.957} \\
PVOT LF & \textbf{.015} & \textbf{.061} & \textbf{.117} &  & \textbf{.122} &
\textbf{.322} & \textbf{.415} &  & \textbf{.740} & \textbf{.911} & \textbf{%
.936} \\
PVOT ICS-1 & \textbf{.015} & \textbf{.057} & \textbf{.116} &  & \textbf{.253}
& \textbf{.464} & \textbf{.570} &  & \textbf{.847} & \textbf{.939} & \textbf{%
.954} \\ \hline\hline
\end{tabular}%
\end{center}
\par
{\small Numerical values are rejection frequency at the given level. LSTAR
is Logistic STAR. Empirical power is not size-adjusted. \textit{sup} $%
\mathcal{T}_{n}$ and \textit{ave} $\mathcal{T}_{n}$ tests are based on a
wild bootstrapped p-value. \textit{rand} $\mathcal{T}_{n}$: $\mathcal{T}%
_{n}(\lambda )$ with uniformly randomly chosen $\lambda $\ on [1,5]. \textit{%
sup} ${p}_{n}$ is the supremum p-value test where p-values are computed from
the chi-squared distribution. PVOT uses the chi-squared distribution. LF
implies the \emph{least favorable} p-value is used, and ICS-$1$ implies the
type 1 \emph{identification category selection} p-value is used with
threshold $\kappa _{n} $ $=$ $\ln (\ln (n))$. }
\end{table}

\clearpage
\begin{table}[h]
\caption{STAR Test Rejection Frequencies: Sample Size $n = 250$}
\label{tbl:starn250}
\begin{center}
\begin{tabular}{l|ccc|c|ccc|c|ccc}
\hline\hline
& \multicolumn{3}{|c|}{$H_{0}$: LSTAR} &  & \multicolumn{3}{|c|}{$H_{1}$-weak
} &  & \multicolumn{3}{|c}{$H_{1}$-strong} \\ \hline
& 1\% & 5\% & 10\% &  & 1\% & 5\% & 10\% &  & 1\% & 5\% & 10\% \\
\hline\hline
& \multicolumn{11}{c}{Strong Identification: $\beta _{n}=.3$} \\ \hline\hline
\multicolumn{1}{l|}{sup $\mathcal{T}_{n}$} & \multicolumn{1}{|l}{.018} &
\multicolumn{1}{l}{.088} & \multicolumn{1}{l|}{.163} &  &
\multicolumn{1}{|l}{.359} & \multicolumn{1}{l}{.468} & \multicolumn{1}{l|}{
.551} &  & \multicolumn{1}{|l}{.953} & \multicolumn{1}{l}{.984} &
\multicolumn{1}{l}{.990} \\
\multicolumn{1}{l|}{aver $\mathcal{T}_{n}$} & \multicolumn{1}{|l}{.014} &
\multicolumn{1}{l}{.077} & \multicolumn{1}{l|}{.133} &  &
\multicolumn{1}{|l}{.262} & \multicolumn{1}{l}{.387} & \multicolumn{1}{l|}{
.468} &  & \multicolumn{1}{|l}{.873} & \multicolumn{1}{l}{.949} &
\multicolumn{1}{l}{.975} \\ \hline
\multicolumn{1}{l|}{rand $\mathcal{T}_{n}$} & \multicolumn{1}{|l}{.014} &
\multicolumn{1}{l}{.064} & \multicolumn{1}{l|}{.126} &  &
\multicolumn{1}{|l}{.165} & \multicolumn{1}{l}{.299} & \multicolumn{1}{l|}{
.396} &  & \multicolumn{1}{|l}{.793} & \multicolumn{1}{l}{.912} &
\multicolumn{1}{l}{.952} \\
\multicolumn{1}{l|}{rand LF} & \textbf{.001} & \textbf{.010} & \textbf{.025}
&  & \textbf{.067} & \textbf{.235} & \textbf{.368} &  & \textbf{.688} &
\textbf{.888} & \textbf{.936} \\
\multicolumn{1}{l|}{rand ICS-1} & \textbf{.008} & \textbf{.031} & \textbf{%
.077} &  & \textbf{.076} & \textbf{.244} & \textbf{.375} &  & \textbf{.762}
& \textbf{.902} & \textbf{.947} \\ \hline
sup $p_{n}$ & .003 & .039 & .066 &  & .103 & .264 & .358 &  & .743 & .876 &
.917 \\
sup $p_{n}$ LF & \textbf{.000} & \textbf{.007} & \textbf{.021} &  & \textbf{%
.032} & \textbf{.214} & \textbf{.303} &  & \textbf{.605} & \textbf{.838} &
\textbf{.899} \\
\multicolumn{1}{l|}{sup $p_{n}$ ICS-1} & \textbf{.003} & \textbf{.035} &
\textbf{.063} &  & \textbf{.038} & \textbf{.217} & \textbf{.316} &  &
\textbf{.714} & \textbf{.870} & \textbf{.912} \\ \hline
\multicolumn{1}{l|}{PVOT} & \multicolumn{1}{|l}{.016} & \multicolumn{1}{l}{
.067} & \multicolumn{1}{l|}{.125} &  & \multicolumn{1}{|l}{.328} &
\multicolumn{1}{l}{.437} & \multicolumn{1}{l|}{.517} &  &
\multicolumn{1}{|l}{.952} & \multicolumn{1}{l}{.983} & \multicolumn{1}{l}{
.991} \\
\multicolumn{1}{l|}{PVOT LF} & \textbf{.004} & \textbf{.020} & \textbf{.041}
&  & \textbf{.132} & \textbf{.348} & \textbf{.417} &  & \textbf{.938} &
\textbf{.972} & \textbf{.976} \\
\multicolumn{1}{l|}{PVOT ICS-1} & \textbf{.011} & \textbf{.051} & \textbf{%
.108} &  & \textbf{.147} & \textbf{.370} & \textbf{.433} &  & \textbf{.947}
& \textbf{.978} & \textbf{.985} \\ \hline\hline
& \multicolumn{11}{c}{Weak Identification: $\beta _{n}=.3/\sqrt{n}$} \\
\hline\hline
\multicolumn{1}{l|}{sup $\mathcal{T}_{n}$} & \multicolumn{1}{|l}{.051} &
\multicolumn{1}{l}{.139} & \multicolumn{1}{l|}{.224} &  &
\multicolumn{1}{|l}{.764} & \multicolumn{1}{l}{.922} & \multicolumn{1}{l|}{
.957} &  & \multicolumn{1}{|l}{.992} & \multicolumn{1}{l}{1.00} &
\multicolumn{1}{l}{1.00} \\
\multicolumn{1}{l|}{aver $\mathcal{T}_{n}$} & \multicolumn{1}{|l}{.046} &
\multicolumn{1}{l}{.118} & \multicolumn{1}{l|}{.215} &  &
\multicolumn{1}{|l}{.539} & \multicolumn{1}{l}{.779} & \multicolumn{1}{l|}{
.853} &  & \multicolumn{1}{|l}{.969} & \multicolumn{1}{l}{.992} &
\multicolumn{1}{l}{.998} \\ \hline
\multicolumn{1}{l|}{rand $\mathcal{T}_{n}$} & \multicolumn{1}{|l}{.027} &
\multicolumn{1}{l}{.086} & \multicolumn{1}{l|}{.169} &  &
\multicolumn{1}{|l}{.451} & \multicolumn{1}{l}{.695} & \multicolumn{1}{l|}{
.785} &  & \multicolumn{1}{|l}{.911} & \multicolumn{1}{l}{.979} &
\multicolumn{1}{l}{.993} \\
\multicolumn{1}{l|}{rand LF} & \textbf{.018} & \textbf{.060} & \textbf{.097}
&  & \textbf{.180} & \textbf{.481} & \textbf{.641} &  & \textbf{.851} &
\textbf{.961} & \textbf{.980} \\
\multicolumn{1}{l|}{rand ICS-1} & \textbf{.018} & \textbf{.058} & \textbf{%
.098} &  & \textbf{.298} & \textbf{.633} & \textbf{.770} &  & \textbf{.926}
& \textbf{.975} & \textbf{.991} \\ \hline
sup $p_{n}$ & .017 & .056 & .097 &  & .330 & .615 & .712 &  & .858 & .975 &
.991 \\
sup $p_{n}$ LF & \textbf{.008} & \textbf{.026} & \textbf{.067} &  & \textbf{%
.115} & \textbf{.416} & \textbf{.587} &  & \textbf{.698} & \textbf{.926} &
\textbf{.978} \\
\multicolumn{1}{l|}{sup $p_{n}$ ICS-1} & \textbf{.008} & \textbf{.030} &
\textbf{.072} &  & \textbf{.294} & \textbf{.580} & \textbf{.687} &  &
\textbf{.852} & \textbf{.975} & \textbf{.991} \\ \hline
\multicolumn{1}{l|}{PVOT} & \multicolumn{1}{|l}{.051} & \multicolumn{1}{l}{
.122} & \multicolumn{1}{l|}{.201} &  & \multicolumn{1}{|l}{.740} &
\multicolumn{1}{l}{.894} & \multicolumn{1}{l|}{.934} &  &
\multicolumn{1}{|l}{1.00} & \multicolumn{1}{l}{1.00} & \multicolumn{1}{l}{
1.00} \\
\multicolumn{1}{l|}{PVOT LF} & \textbf{.014} & \textbf{.061} & \textbf{.110}
&  & \textbf{.380} & \textbf{.708} & \textbf{.805} &  & \textbf{.990} &
\textbf{1.00} & \textbf{1.00} \\
\multicolumn{1}{l|}{PVOT ICS-1} & \textbf{.015} & \textbf{.060} & \textbf{%
.111} &  & \textbf{.618} & \textbf{.848} & \textbf{.878} &  & \textbf{.999}
& \textbf{1.00} & \textbf{1.00} \\ \hline\hline
& \multicolumn{11}{c}{Non-Identification: $\beta _{n}=\beta _{0}=0$} \\
\hline\hline
sup $\mathcal{T}_{n}$ & \multicolumn{1}{|l}{.061} & \multicolumn{1}{l}{.152}
& \multicolumn{1}{l|}{.223} &  & \multicolumn{1}{|l}{.751} &
\multicolumn{1}{l}{.922} & \multicolumn{1}{l|}{.956} &  &
\multicolumn{1}{|l}{1.00} & \multicolumn{1}{l}{1.00} & \multicolumn{1}{l}{
1.00} \\
aver $\mathcal{T}_{n}$ & \multicolumn{1}{|l}{.054} & \multicolumn{1}{l}{.145}
& \multicolumn{1}{l|}{.200} &  & \multicolumn{1}{|l}{.526} &
\multicolumn{1}{l}{.765} & \multicolumn{1}{l|}{.849} &  &
\multicolumn{1}{|l}{.975} & \multicolumn{1}{l}{.996} & \multicolumn{1}{l}{
.999} \\ \hline
rand $\mathcal{T}_{n}$ & \multicolumn{1}{|l}{.036} & \multicolumn{1}{l}{.123}
& \multicolumn{1}{l|}{.184} &  & \multicolumn{1}{|l}{.417} &
\multicolumn{1}{l}{.696} & \multicolumn{1}{l|}{.803} &  &
\multicolumn{1}{|l}{.025} & \multicolumn{1}{l}{.976} & \multicolumn{1}{l}{
.988} \\
rand LF & \textbf{.008} & \textbf{.047} & \textbf{.108} &  & \textbf{.205} &
\textbf{.504} & \textbf{.655} &  & \textbf{.838} & \textbf{.955} & \textbf{%
.973} \\
rand ICS-1 & \textbf{.008} & \textbf{.049} & \textbf{.109} &  & \textbf{.411}
& \textbf{.653} & \textbf{.770} &  & \textbf{.923} & \textbf{.977} & \textbf{%
.989} \\ \hline
sup $p_{n}$ & .026 & .068 & .123 &  & .380 & .650 & .772 &  & .850 & .946 &
.968 \\
sup $p_{n}$ LF & \textbf{.008} & \textbf{.038} & \textbf{.079} &  & \textbf{%
.132} & \textbf{.430} & \textbf{.592} &  & \textbf{.728} & \textbf{.915} &
\textbf{.946} \\
sup $p_{n}$ ICS-1 & \textbf{.008} & \textbf{.004} & \textbf{.081} &  &
\textbf{.340} & \textbf{.629} & \textbf{.750} &  & \textbf{.842} & \textbf{%
.945} & \textbf{.968} \\ \hline
PVOT & \multicolumn{1}{|l}{.036} & \multicolumn{1}{l}{.145} &
\multicolumn{1}{l|}{.211} &  & \multicolumn{1}{|l}{.732} &
\multicolumn{1}{l}{.885} & \multicolumn{1}{l|}{.930} &  &
\multicolumn{1}{|l}{1.00} & \multicolumn{1}{l}{1.00} & \multicolumn{1}{l}{
1.00} \\
PVOT LF & \textbf{.010} & \textbf{.058} & \textbf{.114} &  & \textbf{.373} &
\textbf{.717} & \textbf{.806} &  & \textbf{.990} & \textbf{1.00} & \textbf{%
1.00} \\
PVOT ICS-1 & \textbf{.010} & \textbf{.059} & \textbf{.116} &  & \textbf{.682}
& \textbf{.853} & \textbf{.898} &  & \textbf{1.00} & \textbf{1.00} & \textbf{%
1.00} \\ \hline\hline
\end{tabular}%
\end{center}
\par
{\small Numerical values are rejection frequency at the given level. LSTAR
is Logistic STAR. Empirical power is not size-adjusted. \textit{sup} $%
\mathcal{T}_{n}$ and \textit{ave} $\mathcal{T}_{n}$ tests are based on a
wild bootstrapped p-value. \textit{rand} $\mathcal{T}_{n}$: $\mathcal{T}%
_{n}(\lambda )$ with uniformly randomly chosen $\lambda $\ on [1,5]. \textit{%
sup} ${p}_{n}$ is the supremum p-value test where p-values are computed from
the chi-squared distribution. PVOT uses the chi-squared distribution. LF
implies the \emph{least favorable} p-value is used, and ICS-$1$ implies the
type 1 \emph{identification category selection} p-value is used with
threshold $\kappa _{n} $ $=$ $\ln (\ln (n))$. }
\end{table}

\begin{table}[h]
\caption{STAR Test Rejection Frequencies: Sample Size $n = 500$}
\label{tbl:starn500}
\begin{center}
\begin{tabular}{l|ccc|c|ccc|c|ccc}
\hline\hline
& \multicolumn{3}{|c|}{$H_{0}$: LSTAR} &  & \multicolumn{3}{|c|}{$H_{1}$-weak
} &  & \multicolumn{3}{|c}{$H_{1}$-strong} \\ \hline
& 1\% & 5\% & 10\% &  & 1\% & 5\% & 10\% &  & 1\% & 5\% & 10\% \\
\hline\hline
& \multicolumn{11}{c}{Strong Identification: $\beta _{n}=.3$} \\ \hline\hline
sup $\mathcal{T}_{n}$ & \multicolumn{1}{|l}{.029} & \multicolumn{1}{l}{.069}
& \multicolumn{1}{l|}{.153} &  & \multicolumn{1}{|l}{.441} &
\multicolumn{1}{l}{.590} & \multicolumn{1}{l|}{.676} &  &
\multicolumn{1}{|l}{.997} & \multicolumn{1}{l}{.999} & \multicolumn{1}{l}{
.999} \\
aver $\mathcal{T}_{n}$ & \multicolumn{1}{|l}{.022} & \multicolumn{1}{l}{.055}
& \multicolumn{1}{l|}{.120} &  & \multicolumn{1}{|l}{.382} &
\multicolumn{1}{l}{.546} & \multicolumn{1}{l|}{.624} &  &
\multicolumn{1}{|l}{.988} & \multicolumn{1}{l}{.996} & \multicolumn{1}{l}{
.997} \\ \hline
rand $\mathcal{T}_{n}$ & \multicolumn{1}{|l}{.008} & \multicolumn{1}{l}{.049}
& \multicolumn{1}{l|}{.098} &  & \multicolumn{1}{|l}{.328} &
\multicolumn{1}{l}{.488} & \multicolumn{1}{l|}{.598} &  &
\multicolumn{1}{|l}{.976} & \multicolumn{1}{l}{.999} & \multicolumn{1}{l}{
.996} \\
rand LF & \textbf{.001} & \textbf{.018} & \textbf{.042} &  & \textbf{.227} &
\textbf{.450} & \textbf{.565} &  & \textbf{.967} & \textbf{.989} & \textbf{%
.998} \\
rand ICS-1 & \textbf{.009} & \textbf{.046} & \textbf{.096} &  & \textbf{.230}
& \textbf{.449} & \textbf{.565} &  & \textbf{.974} & \textbf{.990} & \textbf{%
.998} \\ \hline
sup $p_{n}$ & .005 & .039 & .078 &  & .295 & .457 & .536 &  & .961 & .990 &
.997 \\
sup $p_{n}$ LF & \textbf{.002} & \textbf{.010} & \textbf{.033} &  & \textbf{%
.223} & \textbf{.427} & \textbf{.528} &  & \textbf{.949} & \textbf{.985} &
\textbf{.997} \\
sup $p_{n}$ ICS-1 & \textbf{.005} & \textbf{.039} & \textbf{.077} &  &
\textbf{.228} & \textbf{.432} & \textbf{.528} &  & \textbf{.962} & \textbf{%
.990} & \textbf{.997} \\ \hline
PVOT & \multicolumn{1}{|l}{.014} & \multicolumn{1}{l}{.055} &
\multicolumn{1}{l|}{.115} &  & \multicolumn{1}{|l}{.423} &
\multicolumn{1}{l}{.568} & \multicolumn{1}{l|}{.655} &  &
\multicolumn{1}{|l}{.996} & \multicolumn{1}{l}{.999} & \multicolumn{1}{l}{
.999} \\
PVOT LF & \textbf{.002} & \textbf{.023} & \textbf{.051} &  & \textbf{.311} &
\textbf{.509} & \textbf{.618} &  & \textbf{.995} & \textbf{.998} & \textbf{%
1.00} \\
PVOT ICS-1 & \textbf{.013} & \textbf{.058} & \textbf{.106} &  & \textbf{.314}
& \textbf{.510} & \textbf{.618} &  & \textbf{.995} & \textbf{.998} & \textbf{%
1.00} \\ \hline\hline
& \multicolumn{11}{c}{Weak Identification: $\beta _{n}=.3/\sqrt{n}$} \\
\hline\hline
sup $\mathcal{T}_{n}$ & \multicolumn{1}{|l}{.044} & \multicolumn{1}{l}{.134}
& \multicolumn{1}{l|}{.184} &  & \multicolumn{1}{|l}{.984} &
\multicolumn{1}{l}{.998} & \multicolumn{1}{l|}{1.00} &  &
\multicolumn{1}{|l}{1.00} & \multicolumn{1}{l}{1.00} & \multicolumn{1}{l}{
1.00} \\
aver $\mathcal{T}_{n}$ & \multicolumn{1}{|l}{.029} & \multicolumn{1}{l}{.125}
& \multicolumn{1}{l|}{.176} &  & \multicolumn{1}{|l}{.883} &
\multicolumn{1}{l}{.968} & \multicolumn{1}{l|}{/989} &  &
\multicolumn{1}{|l}{1.00} & \multicolumn{1}{l}{1.00} & \multicolumn{1}{l}{
1.00} \\ \hline
rand $\mathcal{T}_{n}$ & \multicolumn{1}{|l}{.032} & \multicolumn{1}{l}{.096}
& \multicolumn{1}{l|}{.162} &  & \multicolumn{1}{|l}{.817} &
\multicolumn{1}{l}{.929} & \multicolumn{1}{l|}{.970} &  &
\multicolumn{1}{|l}{.995} & \multicolumn{1}{l}{.998} & \multicolumn{1}{l}{
.998} \\
rand LF & \textbf{.009} & \textbf{.051} & \textbf{.108} &  & \textbf{.519} &
\textbf{.835} & \textbf{.914} &  & \textbf{.984} & \textbf{.996} & \textbf{%
.998} \\
rand ICS-1 & \textbf{.009} & \textbf{.051} & \textbf{.120} &  & \textbf{.785}
& \textbf{.921} & \textbf{.954} &  & \textbf{.990} & \textbf{.998} & \textbf{%
1.00} \\ \hline
sup $p_{n}$ & .020 & .047 & .093 &  & .721 & .892 & .943 &  & .985 & .998 &
1.00 \\
sup $p_{n}$ LF & \textbf{.015} & \textbf{.025} & \textbf{.054} &  & \textbf{%
.451} & \textbf{.772} & \textbf{.883} &  & \textbf{.961} & \textbf{.992} &
\textbf{1.00} \\
sup $p_{n}$ ICS-1 & \textbf{.014} & \textbf{.026} & \textbf{.056} &  &
\textbf{.710} & \textbf{.890} & \textbf{.940} &  & \textbf{.986} & \textbf{%
.998} & \textbf{1.00} \\ \hline
PVOT & \multicolumn{1}{|l}{.050} & \multicolumn{1}{l}{.118} &
\multicolumn{1}{l|}{.194} &  & \multicolumn{1}{|l}{.981} &
\multicolumn{1}{l}{.995} & \multicolumn{1}{l|}{1.00} &  &
\multicolumn{1}{|l}{1.00} & \multicolumn{1}{l}{1.00} & \multicolumn{1}{l}{
1.00} \\
PVOT LF & \textbf{.012} & \textbf{.053} & \textbf{.109} &  & \textbf{.823} &
\textbf{.965} & \textbf{.975} &  & \textbf{1.00} & \textbf{1.00} & \textbf{%
1.00} \\
PVOT ICS-1 & \textbf{.012} & \textbf{.054} & \textbf{.109} &  & \textbf{.958}
& \textbf{.987} & \textbf{.993} &  & \textbf{1.00} & \textbf{1.00} & \textbf{%
1.00} \\ \hline\hline
& \multicolumn{11}{c}{Non-Identification: $\beta _{n}=\beta _{0}=0$} \\
\hline\hline
sup $\mathcal{T}_{n}$ & \multicolumn{1}{|l}{.051} & \multicolumn{1}{l}{.151}
& \multicolumn{1}{l|}{.196} &  & \multicolumn{1}{|l}{.981} &
\multicolumn{1}{l}{.998} & \multicolumn{1}{l|}{.998} &  &
\multicolumn{1}{|l}{1.00} & \multicolumn{1}{l}{1.00} & \multicolumn{1}{l}{
1.00} \\
aver $\mathcal{T}_{n}$ & \multicolumn{1}{|l}{.043} & \multicolumn{1}{l}{.136}
& \multicolumn{1}{l|}{.189} &  & \multicolumn{1}{|l}{.886} &
\multicolumn{1}{l}{.968} & \multicolumn{1}{l|}{.984} &  &
\multicolumn{1}{|l}{1.00} & \multicolumn{1}{l}{1.00} & \multicolumn{1}{l}{
1.00} \\ \hline
rand $\mathcal{T}_{n}$ & \multicolumn{1}{|l}{.047} & \multicolumn{1}{l}{.111}
& \multicolumn{1}{l|}{.177} &  & \multicolumn{1}{|l}{.826} &
\multicolumn{1}{l}{.938} & \multicolumn{1}{l|}{.967} &  &
\multicolumn{1}{|l}{.997} & \multicolumn{1}{l}{1.00} & \multicolumn{1}{l}{
1.00} \\
rand LF & \textbf{.006} & \textbf{.058} & \textbf{.110} &  & \textbf{.549} &
\textbf{.859} & \textbf{.926} &  & \textbf{1.00} & \textbf{1.00} & \textbf{%
1.00} \\
rand ICS-1 & \textbf{.006} & \textbf{.058} & \textbf{.109} &  & \textbf{.827}
& \textbf{.940} & \textbf{.973} &  & \textbf{1.00} & \textbf{1.00} & \textbf{%
1.00} \\ \hline
sup $p_{n}$ & .032 & .081 & .126 &  & .718 & .904 & .934 &  & .995 & .999 &
.999 \\
sup $p_{n}$ LF & \textbf{.013} & \textbf{.051} & \textbf{.085} &  & \textbf{%
.414} & \textbf{.778} & \textbf{.875} &  & \textbf{.965} & \textbf{.999} &
\textbf{1.00} \\
sup $p_{n}$ ICS-1 & \textbf{.013} & \textbf{.051} & \textbf{.086} &  &
\textbf{.704} & \textbf{.903} & \textbf{.934} &  & \textbf{.995} & \textbf{%
.999} & \textbf{1.00} \\ \hline
PVOT & \multicolumn{1}{|l}{.061} & \multicolumn{1}{l}{.148} &
\multicolumn{1}{l|}{.208} &  & \multicolumn{1}{|l}{.977} &
\multicolumn{1}{l}{.993} & \multicolumn{1}{l|}{.996} &  &
\multicolumn{1}{|l}{1.00} & \multicolumn{1}{l}{1.00} & \multicolumn{1}{l}{
1.00} \\
PVOT LF & \textbf{.014} & \textbf{.058} & \textbf{.108} &  & \textbf{.853} &
\textbf{.970} & \textbf{.989} &  & \textbf{1.00} & \textbf{1.00} & \textbf{%
1.00} \\
PVOT ICS-1 & \textbf{.013} & \textbf{.057} & \textbf{.107} &  & \textbf{.978}
& \textbf{.996} & \textbf{.998} &  & \textbf{1.00} & \textbf{1.00} & \textbf{%
1.00} \\ \hline\hline
\end{tabular}%
\end{center}
\par
{\small Numerical values are rejection frequency at the given level. LSTAR
is Logistic STAR. Empirical power is not size-adjusted. \textit{sup} $%
\mathcal{T}_{n}$ and \textit{ave} $\mathcal{T}_{n}$ tests are based on a
wild bootstrapped p-value. \textit{rand} $\mathcal{T}_{n}$: $\mathcal{T}%
_{n}(\lambda )$ with uniformly randomly chosen $\lambda $\ on [1,5]. \textit{%
sup} ${p}_{n}$is the supremum p-value test where p-values are computed from
the chi-squared distribution. PVOT uses the chi-squared distribution. LF
implies the \emph{least favorable} p-value is used, and ICS-$1$ implies the
type 1 \emph{identification category selection} p-value is used with
threshold $\kappa _{n} $ $=$ $\ln (\ln (n))$. }
\end{table}

\end{document}